\documentclass[11pt,a4paper,usenames,dvipsnames]{article}

\usepackage{url,amsmath,amssymb,latexsym,mathrsfs,comment,amsthm,enumerate,geometry,calc,ifthen,mathdots,textcomp,multicol}
\usepackage{xfrac}
\usepackage[symbol]{footmisc}
\usepackage{cases}

\usepackage[noadjust]{cite}

\usepackage[colorlinks]{hyperref}

\usepackage[margin=10pt,font=small,labelfont=bf, labelsep=period]{caption}

\usepackage{makecell}

\usepackage[all,cmtip,2cell]{xy}

\usepackage[x11names, svgnames, rgb,table]{xcolor}
\usepackage{hhline}
\usepackage{multirow}

\usepackage{pst-node}
\usepackage{tikz-cd}

\usepackage{enumitem}

\usepackage[utf8]{inputenc}
\usepackage[T1]{fontenc}

\usepackage{tikz}
\usepgflibrary{snakes,arrows,shapes}
\pgfdeclarelayer{background layer}
\pgfsetlayers{background layer,main}

\usetikzlibrary{decorations.markings}
\usetikzlibrary{arrows,matrix}
\usepgflibrary{arrows}
  
\tikzset{->-/.style={decoration={
  markings,
  mark=at position #1 with {\arrow{{latex}}}},postaction={decorate}}}

\newcommand\nc\newcommand
\nc\rnc\renewcommand

\nc\trans[1]{\left(\begin{smallmatrix}#1\end{smallmatrix}\right)}
\nc\tmat[1]{\left[\begin{smallmatrix}#1\end{smallmatrix}\right]}
\nc\tmatx[1]{\big[\begin{smallmatrix}#1\end{smallmatrix}\big]}

\nc\RR{\mathbb R}
\nc\CCC{\mathbb C}
\nc\MnR{M_n(\RR)}
\nc\MnC{M_n(\CCC)}
\nc\MnK{M_n(\K)}
\nc\+\dagger
\nc\tr{{\operatorname{T}}}

\nc\FIG{\operatorname{\mathsf{IG}}}
\nc\FPG{\operatorname{\mathsf{PG}}}
\nc\bbB{\mathbb B}
\nc\N{\mathbb N}
\nc\Z{\mathbb Z}
\nc\TL{\mathcal T\!\mathcal L}
\rnc\S{\mathcal S}
\nc{\lmap}[1]{\mapstochar \xrightarrow {\ #1\ }}
\nc\spc{\;\;\!\!}
\nc{\ssim}{\mathrel{\raise0.25 ex\hbox{\oalign{$\approx$\crcr\noalign{\kern-0.84 ex}$\approx$}}}}

\nc\FP{F_P}
\nc\FPd{F_{P'}}
\nc\MP{M_P}

\newcommand{\darcx}[3]{\draw(#1,0)arc(180:90:#3) (#1+#3,#3)--(#2-#3,#3) (#2-#3,#3) arc(90:0:#3);}
\newcommand{\darc}[2]{\darcx{#1}{#2}{.4}}
\newcommand{\uarcx}[3]{\draw(#1,2)arc(180:270:#3) (#1+#3,2-#3)--(#2-#3,2-#3) (#2-#3,2-#3) arc(270:360:#3);}
\newcommand{\uarc}[2]{\uarcx{#1}{#2}{.4}}

\nc{\uarcs}[1]{{\foreach \x/\y in {#1}{ \uarc{\x}{\y} }}}
\nc{\darcs}[1]{{\foreach \x/\y in {#1}{ \darc{\x}{\y} }}}

\newcommand{\darcxcol}[4]{\draw[#4](#1,0)arc(180:90:#3) (#1+#3,#3)--(#2-#3,#3) (#2-#3,#3) arc(90:0:#3);}
\newcommand{\darccol}[3]{\darcxcol{#1}{#2}{.4}{#3}}
\nc{\darccols}[2]{{\foreach \x/\y in {#1}{ \darccol{\x}{\y}{#2} }}}
\newcommand{\uarcxcol}[4]{\draw[#4](#1,2)arc(180:270:#3) (#1+#3,2-#3)--(#2-#3,2-#3) (#2-#3,2-#3) arc(270:360:#3);}
\newcommand{\uarccol}[3]{\uarcxcol{#1}{#2}{.4}{#3}}
\nc{\uarccols}[2]{{\foreach \x/\y in {#1}{ \uarccol{\x}{\y}{#2} }}}

\nc{\uvert}[1]{\fill (#1,2)circle(.2);}
\rnc{\lvert}[1]{\fill (#1,0)circle(.2);}
\nc{\custpartn}[3]{{\lower1.4 ex\hbox{
\begin{tikzpicture}[scale=.3]
\foreach \x in {#1}
{ \uvert{\x}  }
\foreach \x in {#2}
{ \lvert{\x}  }
#3 \end{tikzpicture}
}}}

\newcounter{ncols}
\newcounter{incols}

\nc\la\langle
\nc\ra\rangle

\nc\ol\overline
\nc\ul\underline
\nc\wh\widehat
\nc\da\downarrow
\nc\pr\bullet
\nc\ccirc{\circ'}
\nc\bc\odot
\nc\sm\setminus

\nc\bP{{\bf P}}
\nc\Grp{{\bf Grp}}
\nc\SL{{\bf SL}}
\rnc\SS{{\bf SS}}
\nc\SRS{{\bf SRS}}
\nc\RSS{{\bf RSS}}
\nc\OG{{\bf OG}}
\nc\BC{{\bf BC}}
\nc\DRC{{\bf DRC}}
\nc\DRCR{{\bf DRCR}}
\nc\IGRSS{{\bf IGRSS}}
\nc\IGRSSP{\IGRSS_P}
\nc\PA{{\bf PA}}
\nc\PG{{\bf PG}}
\nc\WCPG{{\bf WCPG}}
\nc\TCPG{{\bf TCPG}}
\nc\CPG{{\bf CPG}}
\nc\PC{{\bf PC}}
\nc\WCPC{{\bf WCPC}}
\nc\TCPC{{\bf TCPC}}
\nc\CPC{{\bf CPC}}
\nc\Set{{\bf Set}}
\nc\IS{{\bf IS}}
\nc\IG{{\bf IG}}
\nc\bG{{\bf G}}
\nc\bO{{\bf O}}
\nc\bS{{\bf S}}
\nc\bC{{\bf C}}
\nc\bD{{\bf D}}
\nc\bI{{\bf I}}
\nc\bT{{\bf T}}
\nc\bbG{\mathbb G}
\nc\bbS{\mathbb S}

\nc\ds\displaystyle
\nc\ts\textstyle

\nc\bn{{\bf n}}

\nc\al\alpha
\nc\be\beta
\nc\ga\gamma
\nc\de\delta
\nc\ve\varepsilon
\nc\lam\lambda
\nc\om\omega
\nc\ze\zeta
\nc\si\sigma
\nc\Si\Sigma
\nc\Ga\Gamma
\nc\De\Delta
\nc\Om\Omega
\nc\io\iota

\nc\nab\nabla

\nc\vt\vartheta
\nc\vd\partial
\rnc\th\theta
\nc\Th\Theta

\nc\BY{\qquad\text{by}\qquad}
\nc\GIVENBY{\qquad\text{given by}\qquad}
\nc\AND{\qquad\text{and}\qquad}
\nc\ANDSIM{,\qquad\text{and similarly}\qquad}
\nc\ANDSO{\qquad\text{and so}\qquad}
\nc\ANd{\quad\text{and}\quad}
\nc\COMMA{,\qquad}
\nc\COMMa{,\quad}
\nc\WHERE{\qquad\text{where}\qquad}
\nc\IE{\qquad\text{i.e.}\qquad}

\rnc\iff{\ \Leftrightarrow\ }
\nc\IFf{\quad \Leftrightarrow\quad }
\nc\Iff{\ \ \Leftrightarrow\ \ }
\nc\IFF{\qquad \Leftrightarrow\qquad }
\rnc\implies{\ \Rightarrow\ }
\nc\IMPLIES{\qquad \Rightarrow\qquad }

\nc\set[2]{\{#1:#2\}}
\nc\bigset[2]{\big\{#1:#2\big\}}
\nc\pres[2]{\la#1:#2\ra}

\nc\bit{\begin{itemize}[label=\textbullet, leftmargin=5mm]}
\nc\eit{\end{itemize}}
\nc\ben{\begin{enumerate}[label=\textup{(\roman*)},leftmargin=10mm]}
\nc\bena{\begin{enumerate}[label=\textup{(\alph*)},leftmargin=10mm]}
\nc\een{\end{enumerate}}
\nc\bmc{\begin{multicols}}
\nc\emc{\end{multicols}}

\nc\pf{\begin{proof}}
\nc\epf{\end{proof}}
\nc\pfclaim{\begin{quote}\begin{proof}}
\nc\epfclaim{\end{proof}\end{quote}}
\nc\epfres{\hfill\qed}
\nc\epfreseq{\tag*{\qed}}
\let\oldproofname=\proofname
\renewcommand{\proofname}{\rm\bf{\oldproofname}}
\nc{\pfitem}[1]{\medskip \noindent #1.}
\nc{\firstpfitem}[1]{#1.}
\nc{\pfcase}[1]{\medskip\noindent {\bf Case #1.}}
\nc\aftercases{\medskip}

\renewcommand{\H}{\mathrel{\mathscr H}}
\renewcommand{\L}{\mathrel{\mathscr L}}
\newcommand{\R}{\mathrel{\mathscr R}}
\newcommand{\D}{\mathrel{\mathscr D}}
\newcommand{\J}{\mathrel{\mathscr J}}
\newcommand{\K}{\mathbb K}
\nc\rH{\mathrel{\H}}
\nc\rL{\mathrel{\L}}
\nc\rR{\mathrel{\R}}
\nc\rD{\mathrel{\D}}
\nc\rJ{\mathrel{\J}}
\nc\rK{\mathrel{\K}}
\nc\rsi{\mathrel{\si}}

\nc\leql{\leq_\ell}
\nc\leqr{\leq_r}
\nc\leqe{\leq_e}
\nc\leqL{\leq_{\L}}
\nc\leqR{\leq_{\R}}
\nc\leqJ{\leq_\J}
\nc\leqH{\leq_\H}
\nc\leqF{\leq_{\F}}
\nc\geqF{\geq_{\F}}

\newcommand{\opp}{{\operatorname{op}}}
\newcommand{\id}{\operatorname{id}}
\newcommand{\dom}{\operatorname{dom}}

\newcommand{\rank}{\operatorname{rank}}
\rnc\Row{\operatorname{Row}}
\rnc\Col{\operatorname{Col}}

\nc\mr\mathrel
\nc\pc[2]{(#1,#2)^\sharp}

\nc\U{\mathcal U}
\nc\V{\mathcal V}
\nc\G{\mathcal G}
\rnc\iff{\ \Leftrightarrow\ }
\rnc\implies{\ \Rightarrow\ }
\nc\Implies{\quad \Rightarrow\quad }
\nc\F{\mathrel{\mathscr F}}
\nc\C{\mathcal C}
\nc\M{\mathcal M}
\nc\CC{\mathscr C}
\nc\DD{\mathcal D}
\nc\bF{{\bf F}}
\nc\I{\mathcal I}
\rnc\O{\mathcal O}
\rnc\P{\mathscr P}
\nc\PP{\mathcal P}
\nc\T{\mathcal T}
\nc\p{\mathfrak p}
\nc\q{\mathfrak q}
\rnc\r{\mathfrak r}
\nc\s{\mathfrak s}
\rnc\t{\mathfrak t}
\nc\bd{{\bf d}}
\nc\br{{\bf r}}
\nc\op\oplus
\nc\lra{\mathrel\leftrightarrow}
\nc\es\varnothing
\nc\rev{\textup{rev}}
\nc\corestt{{\upharpoonleft}}
\nc\restt{{\upharpoonright}}
\nc\corest{{\downharpoonleft}}
\nc\rest{{\downharpoonright}}
\nc\WHERe{\quad\text{where}\quad}
\nc\ot\leftarrow
\rnc\a{\mathfrak a}
\rnc\b{\mathfrak b}
\rnc\c{\mathfrak c}
\rnc\d{\mathfrak d}
\nc\sub\subseteq
\nc\mt\mapsto
\nc\im{\operatorname{im}}
\nc\B{\mathcal B}
\nc\E{\mathbb E}

\numberwithin{equation}{section}

\newtheorem{thm}[equation]{Theorem}
\newtheorem{lemma}[equation]{Lemma}
\newtheorem{cor}[equation]{Corollary}
\newtheorem{prop}[equation]{Proposition}

\theoremstyle{definition}

\newtheorem{defn}[equation]{Definition}
\newtheorem{rem}[equation]{Remark}

\newtheorem{que}[equation]{Question}
\newtheorem{prob}[equation]{Problem}

\newcounter{caseco}

\newcounter{subcaseco}

\newcounter{stepco}

\newcounter{stageco}

\begin{document}

\title{\vspace{-1.5cm}Categorical representation of DRC-semigroups}

\date{}
\author{}

\maketitle

\vspace{-15mm}

\begin{center}
{\large 
James East,%
\hspace{-.25em}\footnote[1]{\label{fn:JE}Centre for Research in Mathematics and Data Science, Western Sydney University, Locked Bag 1797, Penrith NSW 2751, Australia. {\it Emails:} {\tt J.East@westernsydney.edu.au}, {\tt M.Fresacher@westernsydney.edu.au}, {\tt A.ParayilAjmal@westernsydney.edu.au}.}
Matthias Fresacher,%
\hspace{-.25em}\textsuperscript{\ref{fn:JE}}
P.A.~Azeef Muhammed,%
\hspace{-.25em}\textsuperscript{\ref{fn:JE}}
Timothy Stokes\footnote[2]{Department of Mathematics, University of Waikato, Hamilton 3216, New Zealand. {\it Email:} {\tt Tim.Stokes@waikato.ac.nz}.}
}
\end{center}

\vspace{1mm}

\begin{abstract}
DRC-semigroups model associative systems with domain and range operations, and contain many important classes, such as inverse, restriction, Ehresmann, regular $*$-, and $*$-regular semigroups.  In this paper we show that the category of DRC-semigroups is isomorphic to a category of certain biordered categories whose object sets are projection algebras in the sense of Jones.  This extends the recent groupoid approach to regular $*$-semigroups of the first and third authors.  We also establish the existence of free DRC-semigroups by constructing a left adjoint to the forgetful functor into the category of projection algebras.

\medskip

\medskip


\end{abstract}

\tableofcontents

\section{Introduction}\label{sect:intro}

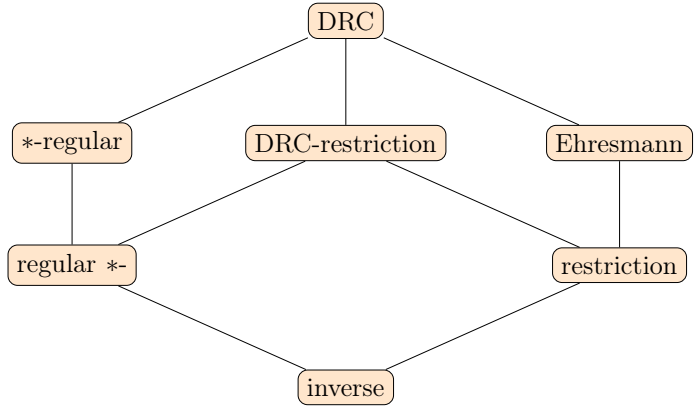
\begin{figure}[ht]
\begin{center}
\scalebox{.9}{
\begin{tikzpicture}[xscale=0.8]
\nc\six5
\nc\two{1.8}
\tikzstyle{vertex}=[circle,draw=black, fill=white, inner sep = 0.07cm]
\node[rounded corners,rectangle,draw,fill=orange!20] (I) at (0,0){inverse};
\node[rounded corners,rectangle,draw,fill=orange!20] (RS) at (-\six,\two){regular $*$-};
\node[rounded corners,rectangle,draw,fill=orange!20] (R) at (\six,\two){restriction};
\node[rounded corners,rectangle,draw,fill=orange!20] (DRCR) at (0,2*\two){DRC-restriction};
\node[rounded corners,rectangle,draw,fill=orange!20] (E) at (\six,2*\two){Ehresmann};
\node[rounded corners,rectangle,draw,fill=orange!20] (SR) at (-\six,2*\two){$*$-regular};
\node[rounded corners,rectangle,draw,fill=orange!20] (DRC) at (0,3*\two){DRC};
\draw (I)--(RS);
\draw (I)--(R);
\draw (RS)--(SR);
\draw (RS)--(DRCR);
\draw (R)--(DRCR);
\draw (R)--(E);
\draw (SR)--(DRC);
\draw (DRCR)--(DRC);
\draw (E)--(DRC);
\end{tikzpicture}
}
\caption{Certain (sub)classes of DRC-semigroups, and the containments among them.}
\label{fig:sgps}
\end{center}
\end{figure}

DRC-semigroups were introduced in \cite{Lawson1991}, and first studied systematically in \cite{Stokes2015} where they were used to model associative systems with domain and range operations.\footnote[3]{DRC-semigroups were referred to in \cite{Lawson1991} as \emph{reduced $U$-semiabundant semigroups satisfying the congruence conditions}.  The current terminology stems from usage in \cite{Stokes2015,Jones2021}.}  Such a semigroup $S$ has additional unary operations $D$ and $R$ satisfying the laws
\begin{align*}
D(a)a &= a ,&
D(ab) &= D(aD(b)) ,&
D(ab) &= D(a)D(ab)D(a) ,& 
R(D(a)) &= D(a),\\
aR(a) &= a ,&
R(ab) &= R(R(a)b) ,&
R(ab) &= R(b)R(ab)R(b) ,&
D(R(a)) &= R(a).
\end{align*}
DRC-semigroups contain many important classes, such as inverse \cite{Lawson1998}, restriction~\cite{Gould_notes}, Ehresmann~\cite{Lawson1991}, regular $*$- \cite{NS1978}, and $*$-regular semigroups \cite{Drazin1979}.  For example, any inverse semigroup is a DRC-semigroup under the operations $D(a) = aa^{-1}$ and $R(a) = a^{-1}a$, or more generally any $*$-regular semigroup under $D(a) = aa^\+$ and $R(a) = a^\+ a$, in terms of the Moore--Penrose inverse.  Prototypical examples of $*$-regular semigroups are the multiplicative monoids~$\MnR$ and~$\MnC$ of real or complex $n\times n$ matrices \cite{Penrose1955}.  The above classes are shown in Figure \ref{fig:sgps}; we refer to \cite{Stokes2015,DW2024} for more on these and other classes, as well as more examples and  historical background.

One of the strong motivations for studying classes of semigroups such as those above stems from the fact that they can be understood categorically, going back to the celebrated Ehresmann--Nambooripad--Schein (ESN) Theorem, as formulated by Lawson in \cite{Lawson1998}.  This result states that the category of inverse semigroups is isomorphic to the category of inductive groupoids, i.e.~the ordered groupoids whose object sets are semilattices.  The main significance of this result is that one can capture the entire semigroup by limited knowledge of the product (the groupoid only `remembers' $ab$ when $R(a) = D(b)$) and the natural partial order (${a\leq b \iff a=D(a)b \iff a=bR(a)}$).  The ESN Theorem has far-reaching consequences and applications, for example to representation theory \cite{Steinberg2006,Steinberg2010,Steinberg2008} and C*-algebras~\cite{Lawson2020,Paterson1999,KS2002}.

One of the first generalisations of the ESN Theorem, to the important class of Ehresmann semigroups and categories, was due to Lawson himself \cite{Lawson1991}.  These include many non-inverse semigroups, such as monoids of binary relations (which are not even regular) and partition monoids \cite{EG2021,MS2021}.  See also \cite{Lawson2021} for a recent alternative approach, and \cite{Wang2016,Lawson2004,DP2018,FitzGerald2019,GW2012,FitzGerald2010,Armstrong1988,Wang2020,Wang2019 ,Nambooripad1979 ,Hollings2012,GH2010,Gould2012,Hollings2010, EPA2024,
Wang2022,FK2021,Stokes2017,Stokes2022,Stokes2022b} for other generalisations.  Of course the main challenge in obtaining an `ESN-type' theorem for a class of semigroups is to identify and axiomatise the appropriate categorical structures used in the representation.  A major part of this is to ascertain the precise structure to place on the object sets, which play the role of the semilattices of Lawson's inductive groupoids.  Typically these are (proper) subsets of the idempotents of the semigroup in question, but they need not be commutative or even closed under multiplication, and in particular they need not be semilattices.

Recent work of the first and third authors \cite{EPA2024} represents regular $*$-semigroups by certain ordered groupoids whose object sets are the projection algebras of Imaoka~\cite{Imaoka1983}.  These are special cases of Jones' projection algebras, defined in \cite{Jones2021} to obtain transformation representations of fundamental DRC-semigroups.  These were subsequently used by Wang in \cite{Wang2022} to obtain an ESN-type theorem for DRC-semigroups, although he did not use categories for his representation, but rather generalisations in which many more products/compositions are required to exist.

The main purpose of the current work is to extend the approaches of \cite{EPA2024} and \cite{Lawson1991} to obtain a purely categorical representation for the class of DRC-semigroups.  Our main result, Theorem~\ref{thm:iso}, obtains an isomorphism between the category of DRC-semigroups and what we call chained projection categories.  These are natural biordered categories whose object sets form projection algebras in the sense of Jones \cite{Jones2021}.  As an application, we establish the existence of free DRC-semigroups, thereby also extending the results of \cite{EGPAR2024}.  We note that special cases of these results, on so-called DRC-restriction semigroups, have been independently obtained by Die and Wang~\cite{DW2024}.  This is the class of DRC-semigroups for which our biordered categories have only a single order, and represents the limit to which the techniques of \cite{EPA2024} generalise without significant modification.  

The paper is organised as follows.  
\bit
\item
We begin in Section \ref{sect:prelim} with the preliminary definitions and results we need on biordered categories.  
\item
In Section \ref{sect:DRC} we recall the definitions and basic properties of DRC semigroups, and construct a functor $\DRC\to\BC$ (in Proposition \ref{prop:Cfunc}), where these are the categories of DRC-semigroups and biordered categories.  
\item
We then turn to projection algebras in Section \ref{sect:P}.  After establishing several preliminary facts, we construct a forgetful functor $\bP:\DRC\to\PA$ (in Proposition \ref{prop:Pfunc}) into the category of projection algebras.
At the object level, this maps a DRC-semigroup $S$ to its underlying projection algebra $\bP(S)$, as defined by Jones \cite{Jones2021}.
We also define the chain category $\CC(P)$ associated to a projection algebra $P$ (see Definition \ref{defn:CP}), which will play an important role in subsequent constructions.
\item
In Section \ref{sect:CPC} we introduce the chained projection categories, and the category $\CPC$ formed by them.  Such a category is in fact a triple $(P,\C,\ve)$, where $\C$ is a biordered category whose object set $P$ is a projection algebra (with close structural ties to $\C$), and where $\ve:\CC(P)\to\C$ is a certain functor from the chain category of $P$.  The main result here is Theorem \ref{thm:Cfunctor}, which constructs a functor $\bC:\DRC\to\CPC$.
\item
In Section \ref{sect:CtoS} we then define a functor $\bS:\CPC\to\DRC$ in the opposite direction; see Theorem~\ref{thm:Sfunctor}.  The main work involves the construction of a DRC-semigroup $\bS(P,\C,\ve)$ associated to a chained projection category $(P,\C,\ve)$; see Definition \ref{defn:pr} and Theorem \ref{thm:S}.
\item
We then prove in Section \ref{sect:iso} that the functors $\bC$ and $\bS$ are in fact mutually inverse isomorphisms between the categories $\DRC$ and $\CPC$, thereby establishing our main result, Theorem \ref{thm:iso}.
\item
In Section \ref{sect:PGP}, we construct a left adjoint to the forgetful functor $\bP:\DRC\to\PA$ from Section~\ref{sect:P}, thereby establishing the existence of free projection-generated DRC-semigroups.  These free semigroups are defined by presentations in terms of generators and defining relations.  Among other things, this has the consequence that every (abstract) projection algebra is the algebra of projections of some DRC-semigroup, a fact that was first established by Jones in \cite{Jones2021} by entirely different methods.  Along the way, we also construct the unique (up to isomorphism) fundamental projection-generated DRC-semigroup with projection algebra~$P$.
\item
Finally, in Sections \ref{sect:*} and \ref{sect:other} we examine the various simplifications that occur for the important special cases of regular $*$-semigroups \cite{EPA2024}, $*$-regular semigroups \cite{Drazin1979}, Ehresmann semigroups~\cite{Lawson1991} and DRC-restriction semigroups \cite{DW2024}.
\eit

\section{Preliminaries on biordered categories}\label{sect:prelim}

Throughout the paper, we identify a small category $\C$ with its morphism set, and identify the objects of $\C$ with the identities, the set of which is denoted $v\C$.  Thus, we have domain and range maps $\bd,\br:\C\to v\C$, for which the following hold, for all $a,b,c\in\C$:
\bit
\item $a\circ b$ exists if and only if $\br(a) = \bd(b)$, in which case $\bd(a\circ b) = \bd(a)$ and $\br(a\circ b) = \br(b)$.  (So we are composing morphisms `left to right'.)
\item $\bd(a)\circ a = a = a\circ \br(a)$.
\item If $a\circ b$ and $b\circ c$ exist, then $(a\circ b)\circ c = a\circ (b\circ c)$.
\item $\bd(p) = \br(p) = p$ for all $p\in v\C$.
\eit
A \emph{left-ordered category} is a pair $(\C,\leq)$, where $\C$ is a small category and $\leq$ a partial order on~$\C$ satisfying the following, for all $a,b,c,d\in\C$:
\bit
\item If $a\leq b$, then $\bd(a)\leq\bd(b)$ and $\br(a)\leq\br(b)$.
\item If $a\leq b$ and $c\leq d$, and if $\br(a)=\bd(c)$ and $\br(b)=\bd(d)$, then $a\circ c\leq b\circ d$.
\item For any $p\in v\C$ with $p\leq\bd(a)$, there exists a unique $u\leq a$ with $\bd(u)=p$.  This element is denoted $u = {}_p\corest a$, and called the \emph{left restriction} of $a$ to $p$.
\eit
(These are the \emph{$\Omega$-structured categories with restrictions}, in the terminology of \cite{Lawson1991}.)
\emph{Right-ordered categories} are defined analogously, where the third item above is replaced by:
\bit
\item For any $q\in v\C$ with $q\leq\br(a)$, there exists a unique $v\leq a$ with $\bd(v)=q$. This element is denoted $v = a\rest_q$, and called the \emph{right restriction} of $a$ to $q$.
\eit
A \emph{biordered category} is a triple $(\C,\leql,\leqr)$, where:
\bit
\item $(\C,\leql)$ is a left-ordered category, and $(\C,\leqr)$ a right-ordered category, and 
\item $\leql$ and $\leqr$ restrict to the same order on $v\C$, meaning that $p\leql q \iff p\leqr q$ for all $p,q\in v\C$.
\eit
(The Ehresmann categories considered in \cite{Lawson1991} are certainly biordered categories in the above sense.)
Unless there is a chance of confusion, we typically speak of `a biordered category $\C$', and assume the orders are named $\leql$ and $\leqr$.  It is important to note that left and right restrictions in a biordered category are quite independent of each other.  For example, if $\br({}_p\corest a) = q$, then we need not have ${}_p\corest a = a\rest_q$; in fact, $a\rest_q$ need not even have domain $p$.  (Some concrete examples are discussed in Section \ref{sect:*}.)

We write $\BC$ for the (large) category of biordered categories.  Morphisms in $\BC$ are the \emph{biordered morphisms}.  Such a morphism $(\C,\leql,\leqr) \to (\C',\leql',\leqr')$ is a functor $\phi:\C\to\C'$ preserving both orders, in the sense that
\begin{equation}\label{eq:bimor}
a\leql b \implies a\phi \leql' b\phi \AND a\leqr b\implies a\phi\leqr' b\phi \qquad\text{for all $a,b\in\C$.}
\end{equation}

Although a biordered category need not be particularly symmetric in its own right, there is a certain symmetry/duality in the category $\BC$ itself.  Specifically, we have an isomorphism $\bO:\BC\to\BC$, given at the object level by $\bO(\C,\leq_l,\leq_r) = (\C^\opp,\leq_r,\leq_l)$.  Here~$\C^\opp$ denotes the \emph{opposite category} to $\C$, in which domains and codomains are swapped, and composition is reversed.  (At the morphism level we simply have $\bO(\phi)=\phi$.)

Throughout the paper, we will typically construct biordered categories by first defining an ordering on objects, and then defining suitable restrictions, as formalised in the next result.  We omit the proof, as it is directly analogous to that of \cite[Lemma 2.3]{EPA2024}.

\begin{lemma}\label{lem:BC}
Suppose $\C$ is a small category for which the following two conditions hold:
\ben
\item \label{*C1} There is a partial order $\leq$ on the object set $v\C$.
\item \label{*C2} For all $a\in\C$, and for all $p\leq\bd(a)$ and $q\leq\br(a)$, there exist morphisms ${}_p\corest a$ and $a\rest_q$ in $\C$ such that the following hold, for all $a,b\in\C$ and $p,q,r,s\in v\C$:
\begin{enumerate}[label=\textup{\textsf{(O\arabic*)}},leftmargin=10mm]
\item \label{O1} If $p\leq\bd(a)$, then $\bd({}_p\corest a)=p$ and $\br({}_p\corest a)\leq\br(a)$.
\item[] If $q\leq\br(a)$, then $\br(a\rest_q)=q$ and $\bd(a\rest_q)\leq\bd(a)$.  \\[-2mm]

\item \label{O2} ${}_{\bd(a)}\corest a = a = a\rest_{\br(a)}$.  \\[-2mm]

\item \label{O3} If $p\leq q\leq\bd(a)$, then ${}_p\corest{}_q\corest a = {}_p\corest a$.
\item[] If $r\leq s\leq\br(a)$, then $a\rest_s\rest_r = a\rest_r$.  \\[-2mm]

\item \label{O4} If $p\leq\bd(a)$ and $\br(a)=\bd(b)$, and if $q=\br({}_p\corest a)$, then ${}_p\corest(a\circ b) = {}_p\corest a\circ{}_q\corest b$.
\item[] If $r\leq\br(b)$ and $\br(a)=\bd(b)$, and if $s=\bd(b\rest_r)$, then $(a\circ b)\rest_r = a\rest_s\circ b\rest_r$.
\end{enumerate}
\een
Then $(\C,\leql,\leqr)$ is a biordered category with orders given by
\begin{equation}\label{eq:aleqb}
\begin{array}{ccccc}
a\leql b &\qquad\Leftrightarrow & \qquad a={}_p\corest b \ \ \text{for some $p\leq\bd(b)$} & \qquad\Leftrightarrow & \qquad a = {}_{\bd(a)}\corest b,\\[2mm]
a\leqr b &\qquad\Leftrightarrow & \qquad a= b\rest_q \ \ \text{for some $q\leq\br(b)$} & \qquad\Leftrightarrow & \qquad a = b\rest_{\br(a)}.
\end{array}
\end{equation}
Moreover, any biordered category has the above form.  \qed
\end{lemma}

For biordered categories $\C$ and $\C'$, it is easy to check that a functor $\phi:\C\to\C'$ is a biordered morphism, i.e.~satisfies \eqref{eq:bimor}, if and only if it satisfies
\begin{equation}\label{eq:bm}
({}_p\corest a)\phi = {}_{p\phi}\corest a\phi \AND (a\rest_q)\phi = a\phi\rest_{q\phi} \qquad\text{for all $a\in\C$, and all $p\leq \bd(a)$ and $q\leq\br(a)$.}
\end{equation}
Here, and henceforth, we denote by $\leq$ the common restriction of $\leql$ and $\leqr$ to $v\C$.

Consider a biordered category $\C$, and let $a\in\C$.  We then have two maps
\[
\vt_a:\bd(a)^\da \to \br(a)^\da \AND \vd_a:\br(a)^\da\to\bd(a)^\da,
\]
given by
\begin{equation}\label{eq:vtvd}
p\vt_a = \br({}_p\corest a) \AND q\vd_a = \bd(a\rest_q) \qquad\text{for $p\leq\bd(a)$ and $q\leq\br(a)$.}
\end{equation}
Here for a partially ordered set $(P,\leq)$, we write $t^\da = \set{s\in P}{s\leq t}$ for the \emph{down-set} of $t\in P$.  It follows quickly from \ref{O3} that 
\begin{equation}\label{eq:vtvd2}
s\vt_{{}_p\corest a} = s\vt_a \AND t\vd_{a\rest_q} = t\vd_a \qquad\text{for all $a\in\C$, and all $s\leq p\leq\bd(a)$ and $t\leq q\leq\br(a)$.}
\end{equation}

A \emph{$v$-congruence} on a small category $\C$ is an equivalence relation $\approx$ on $\C$ satisfying the following, for all $a,b,u,v\in\C$:
\bit
\item $a\approx b \implies [\bd(a)=\bd(b)$ and $\br(a)=\br(b)]$,
\item $a\approx b \implies u\circ a\approx u\circ b$, whenever the stated compositions are defined,
\item $a\approx b \implies a\circ v\approx b\circ v$, whenever the stated compositions are defined.
\eit
If $\C$ is a biordered category, we say that the $v$-congruence $\approx$ is a \emph{biordered congruence} if it additionally satisfies the following, for all $a,b\in\C$:
\bit
\item $a\approx b \implies [{}_p\corest a \approx {}_p\corest b$ and $a\rest_q\approx b\rest_q$] for all $p\leq\bd(a)$ and $q\leq\br(a)$.
\eit
In this case the quotient $\C/{\approx}$ is a biordered category with orders given by
\begin{align*}
\al\leql\be &\IFF a\leql b \ \ \text{for some $a\in \al$ and $b\in\be$,}\\
\ga\leqr\de &\IFF c\leqr d \ \ \text{for some $c\in \ga$ and $d\in\de$.}
\end{align*}

\section{DRC-semigroups}\label{sect:DRC}

We now come to the main object of our study, the class of DRC-semigroups.  We begin in Section~\ref{subsect:DRC} by recalling the definitions, and listing some basic properties; while many of these preliminary lemmas exist in the literature (see for example \cite{Jones2021,Wang2022,Stokes2015}), we give brief proofs to keep the paper self-contained.  In Section \ref{subsect:CS} we show how to construct a biordered category from a DRC-semigroup, leading to a functor $\DRC\to\BC$.

\subsection{Definitions and basic properties}\label{subsect:DRC}

\begin{defn}\label{defn:DRC}
A \emph{DRC-semigroup} is an algebra $(S,\cdot,D,R)$, where $(S,\cdot)$ is a semigroup, and where $D$ and $R$ are unary operations $S\to S$ satisfying the following laws, for all $a,b\in S$:
\begin{enumerate}[label=\textup{\textsf{(DRC\arabic*)}},leftmargin=14mm]\bmc2
\item \label{DRC1} \quad $D(a)a = a$,
\item \label{DRC2} \quad $D(ab) = D(aD(b))$,
\item \label{DRC3} \quad $D(ab) = D(a)D(ab)D(a)$,
\item \label{DRC4} \quad $R(D(a)) = D(a)$,
\item[] $aR(a) = a$,
\item[] $R(ab) = R(R(a)b)$,
\item[] $R(ab) = R(b)R(ab)R(b)$,
\item[] $D(R(a)) = R(a)$.
\emc
\end{enumerate}
Throughout the article, we will almost always identify an algebra with its underlying set, including for a DRC-semigroup $S\equiv(S,\cdot,D,R)$.  Thus, we typically speak of `a DRC-semigroup $S$', write multiplication as juxtaposition, and assume the unary operations are denoted $D$ and $R$.

We denote by $\DRC$ the category of all DRC-semigroups, with DRC-morphisms.  Such a morphism $(S,\cdot,D,R) \to (S',\cdot',D',R')$ in $\DRC$ is a function $\phi:S\to S'$ preserving all the operations, meaning that
\[
(a\cdot b)\phi = (a\phi)\cdot'(b\phi) \COMMA D(a)\phi = D'(a\phi) \AND R(a)\phi = R'(a\phi) \qquad\text{for all $a,b\in S$.}
\]
\end{defn}

As noted in the introduction, Lawson first introduced DRC-semigroups in \cite{Lawson1991}, where he referred to them as \emph{reduced $U$-semiabundant semigroups satisfying the congruence conditions}.  (The laws in axiom \ref{DRC2} are typically called the \emph{congruence conditions} as, for example, the identity $D(ab)=D(aD(b))$ is equivalent, in the presence of the other laws, to the relation $\set{(a,b)\in S\times S}{D(a)=D(b)}$ being a left congruence on $S$.)  Lawson was mainly interested in Ehresmann semigroups, which additionally satisfy the law ${D(a)D(b)=D(b)D(a)}$, and hence also $R(a)R(b) = R(b)R(a)$, which makes \ref{DRC3} redundant in the presence of the other laws.

As with $\BC$, the category $\DRC$ posesses a natural symmetry/duality.  Specifically, if $(S,\cdot,D,R)$ is a DRC-semigroup, and if $(S,\star)$ is the \emph{opposite semigroup} of $S$, where $a\star b = b\cdot a$ for $a,b\in S$, then $(S,\star,R,D)$ is a DRC-semigroup.  Thus, any equality that holds in $(S,\cdot,D,R)$ can be converted to a dual equality by reversing the order of all products, and interchanging $D$ and $R$.

\begin{rem}
Some authors use superscripts for the unary operations on a DRC-semigroup.  For example, Jones \cite{Jones2021} and Wang \cite{Wang2022} write $a^+ = D(a)$ and $a^* = R(a)$.  Here we prefer the functional notation, as it emphasises the connection with domains and ranges in the corresponding categories.  
\end{rem}

Note that \ref{DRC4} says $RD=D$ and $DR=R$.  It follows from these that
\[
D = RD = (DR)D = D(RD) = DD \ANDSIM R=RR.
\]
That is, we have the following consequence of \ref{DRC1}--\ref{DRC4}:
\begin{enumerate}[label=\textup{\textsf{(DRC\arabic*)}},leftmargin=14mm]\bmc2 \addtocounter{enumi}4
\item \label{DRC5} \quad $D(D(a)) = D(a)$,
\item[] $R(R(a)) = R(a)$.
\emc
\end{enumerate}
In calculations to follow, we use $=_1$ to indicate equality in a DRC-semigroup by one or more applications of~\ref{DRC1}, and similarly for \ref{DRC2}--\ref{DRC5}.

It follows from \ref{DRC4} that $\im(D) = \im(R)$.  We denote this common image by $\bP(S)$, so
\[
\bP(S) = \set{D(a)}{a\in S} = \set{R(a)}{a\in S}.
\]
The elements of $\bP(S)$ are called \emph{projections}.  We will typically use the next result without explicit reference.

\begin{lemma}[cf.~{\cite[Proposition 1.2]{Stokes2015}}]
We have $\bP(S) = \set{p\in S}{p^2 = p = D(p) = R(p)}$.
\end{lemma}

\pf
Certainly every element of the stated set belongs to $\im(D) = \bP(S)$.  Conversely, let $p\in\bP(S)$, so that $p = D(a)$ for some $a\in S$.  Then $p = D(a) =_5 D(D(a)) = D(p)$, and similarly $p = R(p)$.  From $p = D(p)$, it also follows that $p^2 = D(p) p =_1 p$.
\epf

Since every projection is an idempotent, it follows that \ref{DRC3} is equivalent (in the presence of the other axioms) to:
\begin{enumerate}[label=\textup{\textsf{(DRC\arabic*)$'$}},leftmargin=15mm] \addtocounter{enumi}2
\item \label{DRC3'} \quad $D(ab) = D(a)D(ab) = D(ab)D(a)$,
\qquad\qquad  $R(ab) = R(b)R(ab) = R(ab)R(b)$.
\end{enumerate}

\begin{lemma}\label{lem:Dpa}
For any $a\in S$ and $p\in\bP(S)$ we have
\[
D(pa) = pD(pa) = D(pa)p \AND R(ap) = pR(ap) = R(ap)p.
\]
\end{lemma}

\pf
The first follows from \ref{DRC3'} and $p = D(p)$; the second is dual.
\epf

Since projections are idempotents, the natural partial order on idempotents restricts to a partial order on $\bP(S)$:
\begin{equation}\label{eq:leqPS}
p\leq q \IFf p = qpq \IFf p = pq = qp \qquad\text{for $p,q\in\bP(S)$.}
\end{equation}
The next result shows that this is equivalent to either of $p=pq$ or $p=qp$.  This is referred to as the \emph{reduced} property in the literature \cite{Stokes2015}.

\begin{cor}
For $p,q\in\bP(S)$, we have $p \leq q \iff p = pq \iff p = qp$. 
\end{cor}

\pf
In light of \eqref{eq:leqPS}, and by symmetry, it suffices to check that $p = pq \implies p = qp$.  But if $p = pq$, then using Lemma \ref{lem:Dpa} we have $p = R(p) = R(pq) = qR(pq) = qR(p) = qp$.
\epf

\begin{rem}
An equivalent definition of DRC-semigroups stems from this order-theoretic perspective.  Thus, a DR-semigroup is defined in \cite{Stokes2015} to be a semigroup $S$ in which a (sub)set of idempotents $P$ is distinguished such that for each $a\in S$ there exist smallest $p,q\in P$ under the natural partial order for which $a=pa=aq$, called $p=D(a)$ and $q=R(a)$.  This defines unary operations $D$ and $R$ on such $S$ (both with image $P$), and it follows that DRC-semigroups are precisely the DR-semigroups satisfying the congruence conditions \ref{DRC2}.
\end{rem}

The set $\bP(S)$ can be given the structure of a \emph{projection algebra}, as defined by Jones in \cite{Jones2021}.  We will recall the formal definition in Section \ref{sect:P}, but the next result will be used to verify the axioms.  
For $p\in\bP(S)$ we define the functions $\th_p,\de_p:\bP(S)\to\bP(S)$ by
\[
q\th_p = R(qp) \AND q\de_p = D(pq) \qquad\text{for $q\in\bP(S)$.}
\]

\begin{lemma}[cf.~{\cite[Proposition 7.2]{Jones2021}}]
\label{lem:PS}
For any $p,q,\in\bP(S)$, we have
\ben
\item \label{PS1} $p\th_p = p$, \qquad
$p\th_{q\th_p} = q\th_p$, \qquad
$\th_q\th_{q\th_p} = \th_q\th_p$, \qquad
$\th_p\de_p = \th_p$, \qquad
$\th_{p\de_q}\th_p = \th_q\th_p$,
\item \label{PS2} $p\de_p = p$, \qquad
$p\de_{q\de_p} = q\de_p$, \qquad
$\de_q\de_{q\de_p} = \de_q\de_p$, \qquad
$\de_p\th_p = \de_p$, \qquad
$\de_{p\th_q}\de_p = \de_q\de_p$.
\een
\end{lemma}

\pf
We just prove \ref{PS1}, as \ref{PS2} is dual.  In what follows, we write $=_*$ to indicate an application of Lemma \ref{lem:Dpa} (and continue to use $=_1$, $=_2$, etc.).  For the first item we have $p\th_p = R(pp) = R(p) = p$, and for the second,
\[
p\th_{q\th_p} = R(pR(qp)) =_* R(R(qp)) =_5 R(qp) = q\th_p.
\]
For the remaining three items, fix $t\in P$, and first note that $t\th_q\th_p = R(R(tq)p) =_2 R(tqp)$.  The third and fifth items follow by combining this with 
\[
t\th_q\th_{q\th_p} = R(R(tq) R(qp)) =_2 R(tq R(qp)) =_* R(tq  pR(qp)) =_1 R(tqp)
\]
and
\[
t\th_{p\de_q}\th_p = R(R(tD(qp))p) =_2 R(tD(qp)p) =_* R(tD(qp)q p) =_1 R(tqp).
\]
For the fourth we have $t\th_p\de_p = D(pR(tp)) =_* D(R(tp)) =_4 R(tp) = t\th_p$.
\epf

\subsection{From DRC-semigroups to biordered categories}\label{subsect:CS}

We now wish to show that a DRC-semigroup $S$ naturally induces a biordered category $\C(S)$ with morphism set $S$, and object/identity set $\bP(S)$.  For this we need the following.

\begin{lemma}\label{lem:Dab}
If $S$ is a DRC-semigroup, and if $a,b\in S$ are such that $R(a)=D(b)$, then $D(ab) = D(a)$ and $R(ab) = R(b)$.
\end{lemma}

\pf
We have $D(ab) =_2 D(aD(b)) = D(aR(a)) =_1 D(a)$; the other is dual.
\epf

As a consequence, we can make the following definition.

\begin{defn}\label{defn:CS}
Given a DRC-semigroup $S$, we define the (small) category $\C=\C(S)$ as follows:
\bit
\item The underlying (morphism) set of $\C$ is $S$.
\item The object/identity set is $v\C = P = \bP(S)$.
\item For $a\in\C$ we have $\bd(a) = D(a)$ and $\br(a) = R(a)$.
\item For $a,b\in\C$ with $\br(a)=\bd(b)$ we have $a\circ b = ab$.
\eit
For $a\in\C$, and for $p\leq\bd(a)$ and $q\leq\br(a)$ we define ${}_p\corest a = pa$ and $a\rest_q = aq$.  (Here $\leq$ is the partial order on $P$ in \eqref{eq:leqPS}.)
\end{defn}

\begin{prop}\label{prop:CS}
For any DRC-semigroup $S$, the category $\C=\C(S)$ is biordered, with~$\leql$ and~$\leqr$ as in \eqref{eq:aleqb}.
\end{prop}

\pf
We just need to check conditions \ref{O1}--\ref{O4} from Lemma \ref{lem:BC}.  By symmetry, we just need to check the properties of left restrictions.

\pfitem{\ref{O1}}  For $p\leq\bd(a)=D(a)$ we have $p = pD(a)$, and then
\[
\bd({}_p\corest a) = D(pa) =_2 D(pD(a)) = D(p) = p.
\]
We also have $R(pa) =_3 R(a)R(pa)R(a)$, which says that $\br({}_p\corest a) = R(pa) \leq R(a) = \br(a)$.

\pfitem{\ref{O2}}  We have ${}_{\bd(a)}\corest a = D(a)a =_1 a$.

\pfitem{\ref{O3}}  From $p\leq q$ we have $p = pq$, and so ${}_p\corest {}_q\corest a = p(qa) = pa = {}_p\corest a$.

\pfitem{\ref{O4}}  Here we have $q = \br({}_p\corest a) = R(pa)$, and also $q = \bd({}_q\corest b) = \bd(qb)$, so
\[
{}_p\corest(a\circ b) = pab =_1paR(pa)b = paqb = pa\circ qb = {}_p\corest a \circ {}_q\corest b.  \qedhere
\]
\epf

\begin{rem}
Note that our orders $\leql$ and $\leqr$ were denoted $\leqr$ and $\leql$ (the other way round) by Lawson in the context of Ehresmann semigroups \cite{Lawson1991}.  We have chosen the current terminology, as our $\leql$ is contained in Green's $\leqL$ preorder, which is defined by $a\leqL b \iff a \in S^1b$.  An analogous statement applies to $\leqr$ and $\leqR$.  These containments are strict in general.  For example, if $M$ is a DRC-monoid with identity $1$, then every element is $\leqL$ and $\leqR$ below $1$, but only the projections are $\leql$ or $\leqr$ below $1$.
\end{rem}

The construction of $\C(S)$ from $S$ is an object map from $\DRC$ to $\BC$.  For a DRC-morphism~$\phi$, we define $\C(\phi) = \phi$.  (Since the underlying set of $\C(S)$ is just $S$, this is well defined.)

\begin{prop}\label{prop:Cfunc}
$\C$ is a functor $\DRC\to\BC$.  
\end{prop}

\pf
It remains to check that any DRC-morphism $\phi:S\to S'$ is a biordered morphism $\C\to\C'$, where $\C=\C(S)$ and $\C'=\C(S')$.  This amounts to checking (cf.~\eqref{eq:bm}) that
\ben
\item \label{Cfunc1} $(a\circ b)\phi = (a\phi)\circ'(b\phi)$ for all $a,b\in\C$ with $\br(a)=\bd(b)$,
\item \label{Cfunc2} $({}_p\corest a)\phi = {}_{p\phi}\corest a\phi$ for all $a\in\C$ and $p\leq\bd(a)$, and
\item \label{Cfunc3} $(a\rest_q)\phi = a\phi\rest_{q\phi}$ for all $a\in\C$ and $q\leq\br(a)$.
\een
Here, we use dashes to distinguish the parameters and operations associated to $\C$ and $\C'$.  We will, however, denote multiplication in both $S$ and $S'$ by juxtaposition.  Next we note that for any $a\in\C$,
\[
\bd'(a\phi) = D'(a\phi) = D(a)\phi = \bd(a)\phi \ANDSIM \br'(a\phi) = \br(a)\phi.
\]
\firstpfitem{\ref{Cfunc1}}  If $\br(a)=\bd(b)$, then $\br'(a\phi) = \br(a)\phi = \bd(b)\phi = \bd'(b\phi)$, and so
\[
(a\circ b)\phi = (a b)\phi = (a\phi)  (b\phi) = (a\phi)\circ' (b\phi).
\]
\firstpfitem{\ref{Cfunc2}}  If $p\leq q$, where for simplicity we write $q=\bd(a)$, then $p = pq$, and so $p\phi = (p\phi)(q\phi)$.  This gives $p\phi\leq q\phi = \bd'(a\phi)$, and it follows that $({}_p\corest a)\phi = (pa)\phi = (p\phi)(a\phi) = {}_{p\phi}\corestt a\phi$.

\pfitem{\ref{Cfunc3}}  This is dual to \ref{Cfunc2}.
\epf

\section{Projection algebras}\label{sect:P}

The set $\bP(S)$ of projections of a DRC-semigroup $S$ need not be a subsemigroup of $S$.  Nevertheless, it can be given the structure of a \emph{projection algebra}, as axiomatised by Jones \cite{Jones2021}.  These algebras will play a crucial role in all that follows, and here we gather the information we need.  In Section \ref{subsect:P} we recall Jones' definition, albeit in an alternative format, and prove some preliminary results (some of which can again be found in the existing literature).  In Section \ref{subsect:PA} we consider the category $\PA$ of projection algebras, and show how the construction of $\bP(S)$ from a DRC-semigroup~$S$ leads to a functor $\DRC\to\PA$.  In Section~\ref{subsect:CP} we show how a projection algebra $P$ induces two biordered categories, the \emph{path category} $\P(P)$ and the \emph{chain category}~$\CC(P)$.  These will be crucial ingredients in constructions of subsequent sections.

\subsection{Definitions and basic properties}\label{subsect:P}

\begin{defn}\label{defn:P}
A \emph{projection algebra} is an algebra $(P,\th,\de)$, where $\th=\set{\th_p}{p\in P}$ and $\de=\set{\de_p}{p\in P}$ are families of unary operations satisfying the following laws, for all $p,q\in P$:
\begin{enumerate}[label=\textup{\textsf{(P\arabic*)}},leftmargin=9mm]\bmc2
\item \label{P1} \quad $p\th_p = p$,
\item \label{P2} \quad $p\th_{q\th_p} = q\th_p$,
\item \label{P3} \quad $\th_q\th_{q\th_p} = \th_q\th_p$,
\item \label{P4} \quad $\th_p\de_p = \th_p$,
\item \label{P5} \quad $\th_{p\de_q}\th_p = \th_q\th_p$,
\item[] $p\de_p = p$,
\item[] $p\de_{q\de_p} = q\de_p$,
\item[] $\de_q\de_{q\de_p} = \de_q\de_p$,
\item[] $\de_p\th_p = \de_p$,
\item[] $\de_{p\th_q}\de_p = \de_q\de_p$.
\emc
\end{enumerate}
The elements of a projection algebra are called \emph{projections}.  We typically speak of `a projection algebra $P$', and assume that its operations are denoted $\th$ and $\de$.
\end{defn}

Again there is a natural symmetry/duality, as whenever $(P,\th,\de)$ is a projection algebra, so too is $(P,\de,\th)$.  As usual, this allows us to shorten proofs, by omitting the proof of any statement whose dual has already been proved.  This will often be done without explicit comment.

\begin{rem}\label{rem:unary}
Jones' original formulation of projection algebras in \cite{Jones2021} utilised two \emph{binary} operations, denoted $\times$ and $\star$, which are related to our unary operations by:
\[
p \times q = q\de_p \AND p \star q = p\th_q \qquad\text{for $p,q\in P$.}
\]
Jones' axioms (with his labelling) are as follows:
\bmc2
\begin{enumerate}[label=\textup{\textsf{(LP\arabic*)}},leftmargin=12mm]
\item \label{LP1} $e\times e = e$,
\item \label{LP2} $(e\times f)\times e = e\times f$,
\item \label{LP3} $g\times(f\times e) = (g\times f)\times(f\times e)$,
\item \label{LP4} $(g\times f)\times e = (g\times f)\times(g\times e)$,
\end{enumerate}
\begin{enumerate}[label=\textup{\textsf{(RP\arabic*)}},leftmargin=12mm]
\item \label{RP1} $e\star e = e$,
\item \label{RP2} $e\star(f\star e) = f\star e$,
\item \label{RP3} $(e\star f)\star g = (e\star f)\star(f\star g)$,
\item \label{RP4} $e\star(f\star g) = (e\star g)\star(f\star g)$,
\end{enumerate}
\emc
\begin{enumerate}[label=\textup{\textsf{(PA\arabic*)}},leftmargin=12mm]
\item \label{PA1} $(f\times e)\star f = f\times e$ and $f\times(e\star f) = e\star f$,
\item \label{PA2} $(e\star(f\times g))\star g = (e\star f)\star g$ and $g\times((g\star f)\times e) = g\times(f\times e)$.
\end{enumerate}
Our axioms \ref{P1}--\ref{P3} correspond to Jones' \ref{LP1}--\ref{LP3} and \ref{RP1}--\ref{RP3}, and our \ref{P4}--\ref{P5} are Jones' \ref{PA1}--\ref{PA2}.  
For example, directly translating \ref{RP3} to the unary setting yields $e\th_f\th_g = e\th_f\th_{f\th_g}$, which is equivalent to the equality of maps, $\th_f\th_g = \th_f\th_{f\th_g}$ (requiring one less parameter), the first part of our \ref{P3}.  

Jones' \ref{LP4} and \ref{RP4} are absent from our list, as they are consequences of the other axioms; see Lemma \ref{lem:P} and Remark \ref{rem:P7}.  The reason Jones included those axioms is that \ref{LP1}--\ref{LP4} are axioms for `left projection algebras', which are used to model projections of `DC-semigroups', and \ref{LP4} cannot be deduced from \ref{LP1}--\ref{LP3}, as can be shown for example with Mace4 \cite{Mace4}.  Similar comments apply to \ref{RP1}--\ref{RP4}.

One reason we have opted for the unary approach is due to associativity of function composition, meaning that bracketing becomes unnecessary for complex projection algebra terms.  Another is to emphasise the connection with certain natural mappings $\vt_a$, $\Th_a$, $\vd_a$ and $\De_a$ on associated categories, which will be introduced below.  For a more detailed comparison of the binary and unary approaches, we refer to \cite[Remark 4.2]{EPA2024}.  
\end{rem}

Here then is the connection between DRC-semigroups and projection algebras.  The following is well defined because of Lemma \ref{lem:PS}.

\begin{defn}[cf.~{\cite[Proposition 7.2]{Jones2021}}]
\label{defn:PS}
The projection algebra of a DRC-semigroup $S$ has
\bit
\item underlying set $\bP(S) = \im(D) = \im(R) = \set{p\in S}{p^2=p=D(p)=R(p)}$, and
\item $\th$ and $\de$ operations defined by
\eit
\[
q\th_p = R(qp) \AND q\de_p = D(pq) \qquad\text{for $p,q\in\bP(S)$.}
\]
\end{defn}

For the rest of this section, we fix a projection algebra $P$.  In calculations to follow, we use~$=_1$ to indicate equality by one or more applications of \ref{P1}, and similarly for the other axioms.  

\begin{lemma}\label{lem:leqleqF}
For any $p,q\in P$ we have
\[
p = p\th_q \implies p = q\th_p \AND p = p\de_q \implies p = q\de_p.
\]
\end{lemma}

\pf
For the first implication (the second is dual), suppose $p = p\th_q$.  Then 
\[
q\th_p =_2 p\th_{q\th_p} = p\th_q\th_{q\th_p} =_3 p\th_q\th_p = p\th_p =_1 p.  \qedhere
\]
\epf

The next lemma gathers some fundamental consequences of \ref{P1}--\ref{P5}, which will be used so frequently as to warrant naming them \ref{P6}--\ref{P10}.

\newpage

\begin{lemma}\label{lem:P}
For any $p,q\in P$ we have
\begin{enumerate}[label=\textup{\textsf{(P\arabic*)}},leftmargin=11mm]\bmc2 \addtocounter{enumi}5
\item \label{P6} \quad $\th_p\th_p = \th_p$,
\item \label{P7} \quad $\th_{q\th_p} = \th_p\th_{q\th_p} = \th_{q\th_p}\th_p$,
\item \label{P8} \quad $q\th_{q\th_p} = q\th_p$,
\item \label{P9} \quad $p\de_q\th_p = q\th_p$,
\item \label{P10} \quad $\de_q\th_{q\th_p} = \de_q\th_p$,
\item[] $\de_p\de_p=\de_p$,
\item[] $\de_{q\de_p} = \de_p\de_{q\de_p} = \de_{q\de_p}\de_p$,
\item[] $q\de_{q\de_p} = q\de_p$,
\item[] $p\th_q\de_p = q\de_p$,
\item[] $\th_q\de_{q\de_p} = \th_q\de_p$.
\emc
\end{enumerate}
\end{lemma}

\pf
\firstpfitem{\ref{P6}} We have $\th_p =_4 \th_p\de_p =_4 \th_p\de_p\th_p =_4 \th_p\th_p$.

\pfitem{\ref{P7}}  Using the now-established \ref{P6}, we have 
\[
\th_p\th_{q\th_p} =_5 \th_{q\th_p\de_p}\th_{q\th_p} =_4 \th_{q\th_p}\th_{q\th_p} =_6 \th_{q\th_p}
\AND
\th_{q\th_p}\th_p =_3 \th_{q\th_p}\th_{q\th_p\th_p} =_6 \th_{q\th_p}\th_{q\th_p} =_6 \th_{q\th_p}.
\]
\firstpfitem{\ref{P8}}  We have $q\th_{q\th_p} =_1 q\th_q\th_{q\th_p} =_3 q\th_q\th_p =_1 q\th_p$.

\pfitem{\ref{P9}}  Since $p\de_q =_4 (p\de_q)\th_q$, it follows by Lemma \ref{lem:leqleqF} that $p\de_q = q\th_{p\de_q}$, and then
\[
p\de_q\th_p = q\th_{p\de_q}\th_p =_5 q\th_q\th_p =_1 q\th_p.
\]
\firstpfitem{\ref{P10}}  We have $\de_q\th_{q\th_p} =_4 \de_q\th_q\th_{q\th_p} =_3 \de_q\th_q\th_p =_4 \de_q\th_p$.
\epf

\begin{rem}\label{rem:P7}
The $\th_{q\th_p} = \th_p\th_{q\th_p}$ and $\de_{q\de_p} = \de_p\de_{q\de_p}$ parts of \ref{P7} are Jones' axioms~\ref{RP4} and~\ref{LP4}; cf.~Remark \ref{rem:unary}.
\end{rem}

\begin{lemma}\label{lem:leqP}
For $p,q\in P$, the following are equivalent:
\ben
\item \label{leqP1} $p = p\th_q$,
\item \label{leqP2} $p = r\th_q$ for some $r\in P$, i.e.~$p\in\im(\th_q)$,
\item \label{leqP3} $p = p\de_q$,
\item \label{leqP4} $p = r\de_q$ for some $r\in P$, i.e.~$p\in\im(\de_q)$,
\item \label{leqP5} $p = p\th_q = p\de_q = q\th_p = q\de_p$.
\een
\end{lemma}

\pf
The following implications are clear:
\[
\ref{leqP5}\implies\ref{leqP1}\implies\ref{leqP2} \AND \ref{leqP5}\implies\ref{leqP3}\implies\ref{leqP4}.
\]
We will show that $\ref{leqP2}\implies\ref{leqP1}\implies\ref{leqP5}$, and then $\ref{leqP4}\implies\ref{leqP3}\implies\ref{leqP5}$ will follow by duality.

\pfitem{$\ref{leqP2}\implies\ref{leqP1}$}  If $p = r\th_q$ for some $r\in P$, then $p = r\th_q =_6 r\th_q\th_q = p\th_q$.

\pfitem{$\ref{leqP1}\implies\ref{leqP5}$}  If $p = p\th_q$, then $p = p\th_q =_4 p\th_q\de_q = p\de_q$.  Lemma \ref{lem:leqleqF} then gives $p = q\th_p = q\de_p$.
\epf

\begin{defn}
For $p,q\in P$ we say that $p\leq q$ if any (and hence all) of conditions \ref{leqP1}--\ref{leqP5} of Lemma \ref{lem:leqP} hold.  That is,
\begin{equation}\label{eq:leqP}
p\leq q \IFf p \in\im(\th_q) \IFf p \in\im(\de_q) \IFf p = p\th_q = p\de_q = q\th_p = q\de_p.
\end{equation}
\end{defn}

\begin{prop}[cf.~{\cite[Lemma 5.2]{Jones2021}}]
$\leq$ is a partial order.
\end{prop}

\pf
Reflexivity follows from \ref{P1}.  For anti-symmetry, suppose $p\leq q$ and $q\leq p$, so that $p=p\th_q$ and $q=q\th_p$.  Then $p = p\th_q =_2 q\th_{p\th_q} = q\th_p = q$.  For transitivity, suppose $p\leq q$ and $q\leq r$, so that $p=p\th_q$ and $q=q\th_r$.  Then $p = p\th_q =_6 p\th_q\th_q = p\th_q\th_{q\th_r} =_3 p\th_q\th_r = p\th_r$, which gives $p\leq r$.
\epf

It follows from \eqref{eq:leqP} that
\[
\im(\th_p) = \im(\de_p) = p^\da = \set{q\in P}{q\leq p} \qquad\text{for all $p\in P$.}
\]
The next basic result will be used several times.

\begin{lemma}\label{lem:thpthq}
If $p\leq q$, then $\th_p=\th_p\th_q=\th_q\th_p$ and $\de_p=\de_p\de_q=\de_q\de_p$.
\end{lemma}

\pf
If $p\leq q$, then $p=p\th_q$, and so $\th_p = \th_{p\th_q} =_7 \th_q\th_{p\th_q} = \th_q\th_p$.  Essentially the same calculation gives $\th_p=\th_p\th_q$.
\epf

In addition to the partial order $\leq$, a crucial role in all that follows will be played by the following relation.

\begin{defn}\label{defn:F}
We define the relation $\F$ on a projection algebra $P$ by
\[
p\F q \IFF p = q\de_p \ANd q = p\th_q \qquad\text{for $p,q\in P$.}
\]
Note that $\F$ is reflexive by \ref{P1}, but need not be symmetric or transitive.
\end{defn}

In the next sequence of results we gather some basic facts about the relation $\F$, and its interaction with the partial order $\leq$.

\begin{lemma}\label{lem:p'q'}
Let $p,q\in P$, and put $p' = q\de_p$ and $q' = p\th_q$.  Then:
\ben
\item \label{p'q'1} $p' \leq p$, $q'\leq q$ and $p' \F q'$,
\item \label{p'q'2} $\th_p\th_q = \th_{p'}\th_{q'}$ and $\de_q\de_p = \de_{q'}\de_{p'}$,
\item \label{p'q'3} if $p\leq r \F q$ for some $r\in P$, then $p'=p$,
\item \label{p'q'4} if $p \F s \geq q$ for some $s\in P$, then $q'=q$.
\een
\end{lemma}

\pf
\firstpfitem{\ref{p'q'1}}  We obtain $p'\leq p$ and $q'\leq q$ by definition.  We also have
\[
q'\de_{p'} = p\th_q \de_{q\de_p} =_{10} p\th_q\de_p =_9 q\de_p = p'.
\]
We obtain $p'\th_{q'} = q'$ by symmetry, and hence $p'\F q'$.

\pfitem{\ref{p'q'2}}  We have $\th_p\th_q =_3 \th_p\th_{p\th_q} =_5 \th_{p\th_q\de_p} \th_{p\th_q} =_9 \th_{q\de_p}\th_{p\th_q} = \th_{p'}\th_{q'}$.  

\pfitem{\ref{p'q'3}}  The assumptions give $p=r\de_p$ and $r=q\de_r$ (among other things).  We also have $\de_p=\de_r\de_p$ by Lemma \ref{lem:thpthq}.  But then $p' = q\de_p = q\de_r\de_p = r\de_p = p$.
\epf

\begin{cor}\label{cor:pqrs}
If $p\F q$, then for any $r\leq p$ and $s\leq q$ we have $r\F r\th_q$ and $s\de_p \F s$.
\end{cor}

\pf
For the first (the second is dual), part \ref{p'q'1} of Lemma \ref{lem:p'q'} gives $r' \F q'$, where $r' = q\de_r$ and $q' = r\th_q$.  Part \ref{p'q'3} of the same lemma gives $r'=r$.
\epf

\begin{lemma}\label{lem:pkql}
If $p_1,\ldots,p_k,q_1,\ldots,q_l\in P$ are such that $p_1\F\cdots\F p_k$ and $q_1\F\cdots\F q_k$, then
\[
\th_{p_1}\cdots\th_{p_k} = \th_{q_1}\cdots\th_{q_l} \implies p_k=q_l
\AND
\de_{p_k}\cdots\de_{p_1} = \de_{q_l}\cdots\de_{q_1} \implies p_1=q_1.
\]
\end{lemma}

\pf
For the first (the second is dual), we have
\[
p_k = p_1\th_{p_1}\cdots\th_{p_k} = p_1\th_{q_1}\cdots\th_{q_l} \leq q_l \ANDSIM q_l\leq p_k.  \qedhere
\]
\epf

The $k=l=1$ case of the previous result has the following simple consequence.

\begin{cor}\label{cor:thp=thq}
For $p,q\in P$, we have $\th_p=\th_q \iff \de_p=\de_q \iff p=q$.  \qed
\end{cor}

It will be convenient to record a result concerning products of projections in a DRC-semigroup.

\begin{lemma}\label{lem:Dpp}
If $p_1,\ldots,p_k\in\bP(S)$ for a DRC-semigroup $S$, then
\ben
\item \label{Dpp1} $D(p_1\cdots p_k) = p_k\de_{p_{k-1}}\cdots\de_{p_1}$ and $R(p_1\cdots p_k) = p_1\th_{p_2}\cdots\th_{p_k}$,
\item \label{Dpp2} $D(p_1\cdots p_k) = p_1$ and $R(p_1\cdots p_k) = p_k$ if $p_1\F\cdots\F p_k$
\een
\end{lemma}

\pf
\firstpfitem{\ref{Dpp1}}  The $k=1$ case says $D(p_1)=p_1$, which holds by \ref{DRC5}.  For $k\geq2$ we have
\begin{align*}
D(p_1p_2\cdots p_k) &= D(p_1D(p_2\cdots p_k)) &&\text{by \ref{DRC2}}\\
&= D(p_2\cdots p_k) \de_{p_1} &&\text{by definition}\\
&=  p_k\de_{p_{k-1}}\cdots\de_{p_2} \de_{p_1} &&\text{by induction.}
\end{align*}
\firstpfitem{\ref{Dpp2}}  This follows from \ref{Dpp1} and the definition of the $\F$ relation.
\epf

\subsection{The category of projection algebras}\label{subsect:PA}

\begin{defn}\label{defn:PA}
We write $\PA$ for the category of projection algebras.  A \emph{projection algebra morphism} ${(P,\th,\de)\to(P',\th',\de')}$ in $\PA$ is a map $\phi:P\to P'$ such that
\[
(q\th_p)\phi = (q\phi)\th'_{p\phi} \AND (q\de_p)\phi = (q\phi)\de'_{p\phi} \qquad\text{for all $p,q\in P$.}
\]
\end{defn}

In Definition \ref{defn:PS} (cf.~Lemma \ref{lem:PS}) we saw that a DRC-semigroup $S$ gives rise to a projection algebra $\bP(S)$.  In this way, $\bP$ can be thought of as an object map $\DRC\to\PA$.  Since a DRC-morphism $\phi:S\to S'$ maps projections to projections (by the law $D(a)\phi=D'(a\phi)$), it follows that we can define
\[
\bP(\phi) = \phi|_{\bP(S)}:\bP(S)\to\bP(S')
\]
to be the (set-theoretic) restriction of $\phi$ to $\bP(S)$.

\begin{prop}\label{prop:Pfunc}
$\bP$ is a functor $\DRC\to\PA$.
\end{prop}

\pf
It only remains to show that $\bP(\phi):\bP(S)\to\bP(S')$ is a projection algebra morphism for any DRC-morphism $\phi:S\to S'$.  But for any $p,q\in \bP(S)$ we have
\[
(q\th_p)\phi = R(qp)\phi = R'((qp)\phi) = R'((q\phi)(p\phi)) = (q\phi)\th'_{p\phi}
\ANDSIM
(q\de_p)\phi = (q\phi)\de'_{p\phi}.
\qedhere
\]
\epf

\subsection{The chain category of a projection algebra}\label{subsect:CP}

\begin{defn}\label{defn:PP}
Let $P$ be a projection algebra.  A \emph{($P$-)path} is a non-empty tuple $\p = (p_1,\ldots,p_k)$ where $p_1,\ldots,p_k\in P$ and $p_1\F\cdots\F p_k$.  We write $\bd(\p)=p_1$ and $\br(\p)=p_k$.  The set $\P = \P(P)$ of all such paths is the \emph{path category} of $P$, under the composition defined by
\[
(p_1,\ldots,p_k) \circ (p_k,\ldots,p_l) = (p_1,\ldots,p_{k-1},p_k,p_{k+1},\ldots,p_l).
\]
For $p\in P$, we identify the path $(p)$ with $p$ itself, and in this way we have $v\P = P$.
\end{defn}

We now wish to give the path category $\P = \P(P)$ the structure of a biordered category.  To this end, consider a path $\p = (p_1,\ldots,p_k) \in \P$, and let $q\leq \bd(\p) = p_1$ and $r \leq \br(\p) = p_k$.  We then define the restrictions
\begin{align}
\label{eq:q|p} {}_q\corest \p &= (q_1,\ldots,q_k) , &\text{where}&& q_i &= q\th_{p_1}\cdots\th_{p_i} = q\th_{p_2}\cdots\th_{p_i} &&\text{for each $1\leq i\leq k$,}
\intertext{and}
\label{eq:p|r} \p\rest_r &= (r_1,\ldots,r_k), &\text{where}&& r_i &= r\de_{p_k}\cdots\de_{p_i} = r\de_{p_{k-1}}\cdots\de_{p_i}  &&\text{for each $1\leq i\leq k$.}
\end{align}
Note here that $q\th_{p_1} = q$ and $r\de_{p_k} = r$, since $q\leq p_1$ and $r\leq p_k$, so in particular we have $q_1=q$ and $r_k=r$.
We now prove a sequence of results that show these restrictions satisfy conditions~\ref{O1}--\ref{O4} from Lemma \ref{lem:BC}.

\begin{lemma}\label{lem:leql1}
If $\p\in\P$, and if $q\leq \bd(\p)$, then ${}_q\corest\p\in\P$, and we have $\bd({}_q\corest\p) = q$ and $\br({}_q\corest\p) \leq \br(\p)$.
\end{lemma}

\pf
Let $\p=(p_1,\ldots,p_k)$, and write ${}_q\corest\p = (q_1,\ldots,q_k)$ as in \eqref{eq:q|p}.  Since $\p\in\P$ we have
\[
p_i \F p_{i+1} , \IE p_i\th_{p_{i+1}} = p_{i+1} \text{ \ and \ } p_{i+1}\de_{p_i} = p_i \qquad\text{for all $1\leq i<k$.}
\]
To show that ${}_q\corest\p\in\P$, we must show that $q_i\F q_{i+1}$ for all $1\leq i<k$, i.e.~that $q_i\th_{q_{i+1}} = q_{i+1}$ and $q_{i+1}\de_{q_i} = q_i$.  By definition, we note that $q_{i+1} = q_i\th_{p_{i+1}}$.  We then have
\[
q_i\th_{q_{i+1}} = q_i\th_{q_i\th_{p_{i+1}}} =_8 q_i\th_{p_{i+1}} = q_{i+1}.
\]
Since $q_i\leq p_i$, we have  $p_i\de_{q_i} = q_i$, and also $\de_{q_i} = \de_{p_i}\de_{q_i}$ by Lemma \ref{lem:thpthq}, so 
\[
q_{i+1}\de_{q_i} = q_i \th_{p_{i+1}}\de_{q_i} =_9 p_{i+1}\de_{q_i} = p_{i+1}\de_{p_i}\de_{q_i} = p_i\de_{q_i} = q_i.
\]

By definition, we have $\bd({}_q\corest\p) = q_1 = q$, and $\br({}_q\corest\p) = q_k = q\th_{p_1}\cdots\th_{p_k} \leq p_k$.
\epf

\begin{lemma}\label{lem:leql2}
If $\p\in\P$, and if $q=\bd(\p)$, then ${}_q\corest\p = \p$.
\end{lemma}

\pf
Let $\p=(p_1,\ldots,p_k)$, and write ${}_q\corest\p = (q_1,\ldots,q_k)$ as in \eqref{eq:q|p}.  We need to show that $q_i=p_i$ for all $1\leq i\leq k$.  This is true for $i=1$, as $q_1 = q = \bd(\p) = p_1$.  If $i\geq2$, then by induction, and since $p_{i-1}\F p_i$, we have $q_i = q_{i-1}\th_{p_i} = p_{i-1}\th_{p_i} = p_i$.
\epf

\begin{lemma}\label{lem:leql3}
If $\p\in\P$, and if $r\leq q\leq \bd(\p)$, then ${}_r\corest{}_q\corest\p = {}_r\corest\p$.
\end{lemma}

\pf
Write
\[
\p=(p_1,\ldots,p_k) \COMMA {}_q\corest\p=(q_1,\ldots,q_k) \COMMA {}_r\corest\p=(r_1,\ldots,r_k) \AND {}_r\corest{}_q\corest\p=(s_1,\ldots,s_k).
\]
So
\begin{equation}\label{eq:pqrsi}
q_i=q\th_{p_1}\cdots\th_{p_i} \COMMA r_i=r\th_{p_1}\cdots\th_{p_i} \AND s_i=r\th_{q_1}\cdots\th_{q_i} \qquad\text{for all $1\leq i\leq k$,}
\end{equation}
and we must show that $r_i=s_i$ for all $i$.  For $i=1$ we have $r_1=r=s_1$.  For $i\geq2$,
\begin{align*}
s_i &= s_{i-1}\th_{q_i}   &&\text{by \eqref{eq:pqrsi}}\\
&= s_{i-1}\th_{q_{i-1}\th_{p_i}}  &&\text{by \eqref{eq:pqrsi}}\\
&= (s_{i-1}\th_{q_{i-1}})\th_{q_{i-1}\th_{p_i}}  &&\text{by \eqref{eq:pqrsi} and \ref{P6}}\\
&= s_{i-1}\th_{q_{i-1}}\th_{p_i}  &&\text{by \ref{P3}}\\
&= s_{i-1}\th_{p_i}  &&\text{by \eqref{eq:pqrsi} and \ref{P6}}\\
&= r_{i-1}\th_{p_i}  &&\text{by induction}\\
&= r_i  &&\text{by \eqref{eq:pqrsi}.}  \qedhere
\end{align*}
\epf

\begin{lemma}\label{lem:leql4}
If $\p,\q\in\P$ with $\br(\p) = \bd(\q)$, and if $r\leq\bd(\p)$, then with $s=\br({}_r\corest\p)$ we have
\[
{}_r\corest(\p\circ\q) = {}_r\corest\p \circ {}_{s}\corest\q.
\]
\end{lemma}

\pf
Write $\p=(p_1,\ldots,p_k)$ and $\q=(q_1,\ldots,q_l)$, noting that $p_k=q_1$.  Then
\begin{align*}
{}_r\corest\p&=(r,r\th_{p_2},r\th_{p_2}\th_{p_3},\ldots,r\th_{p_2}\cdots\th_{p_k}), \qquad\text{so}\qquad s=r\th_{p_2}\cdots\th_{p_k}.
\intertext{It follows that}
{}_r\corest(\p\circ\q) &= {}_r\corest(p_1,\ldots,p_k,q_2,\ldots,q_l) \\
&= (r,r\th_{p_2},r\th_{p_2}\th_{p_3},\ldots,\underbrace{r\th_{p_2}\cdots\th_{p_k}}_{=s},s\th_{q_2},s\th_{q_2}\th_{q_3},\ldots,s\th_{q_2}\cdots\th_{q_l}) \\
&= (r,r\th_{p_2},r\th_{p_2}\th_{p_3},\ldots,r\th_{p_2}\cdots\th_{p_k})\circ(s,s\th_{q_2},s\th_{q_2}\th_{q_3},\ldots,s\th_{q_2}\cdots\th_{q_l})  = {}_r\corest\p \circ {}_{s}\corest\q. \qedhere
\end{align*}
\epf

Lemmas \ref{lem:leql1}--\ref{lem:leql4} and their duals show that the restrictions given in \eqref{eq:q|p} and \eqref{eq:p|r} satisfy the conditions of Lemma \ref{lem:BC}.  It follows that $\P = \P(P)$ is a biordered category, under the orders~$\leql$ and~$\leqr$ defined, for $\p,\q\in\P$, by
\[
\p\leql\q \iff \p = {}_{\bd(\p)}\corest\q
\AND
\p\leqr\q \iff \p = \q\rest_{\br(\p)}.
\]

\begin{defn}\label{defn:CP}
Let $P$ be a projection algebra, and let $\P = \P(P)$ be the path category.  We let $\approx$ be the congruence on $\P$ generated by the relations $(p,p) \approx p$ for all $p\in P$.  Since
\[
{}_q\corest(p,p) = (q,q) \approx q = {}_q\corest p \ANDSIM (p,p)\rest_q \approx p\rest_q \qquad\text{for $q\leq p$,}
\]
it follows that $\approx$ is a biordered congruence.  Consequently, the quotient
\[
\CC = \CC(P) = \P/{\approx}
\]
is a biordered category under the induced $\leql$ and $\leqr$ orders.  We call $\CC = \CC(P)$ the \emph{chain category of $P$}.  An element of $\CC$ is called a \emph{($P$-)chain}, and is an $\approx$-class of a path.  For a path $\p=(p_1,\ldots,p_k)\in\P$, we write $[\p] = [p_1,\ldots,p_k]$ for the corresponding chain, so that $\CC = \set{[\p]}{\p\in\P}$.  We again identify a projection $p\in P$ with the chain $[p]\in\CC$, so that~$v\CC = P$.  Restrictions in $\CC$ are given by
\[
{}_q\corest[\p] = [{}_q\corest\p] \AND [\p]\rest_r = [\p\rest_r] \qquad\text{for $\p\in\P$, and $q\leq\bd(\p)$ and $r\leq\br(\p)$.}
\]
\end{defn}

\section{Chained projection categories}\label{sect:CPC}

In Proposition \ref{prop:CS} we saw that a DRC-semigroup $S$ gives rise to a biordered category $\C(S)$ whose object set is the projection algebra $\bP(S)$.  It turns out that the category $\C(S)$ carries rather a lot more information, and has the structure of what we will call a \emph{chained projection category}.  In this section we introduce this class of categories, and the category $\CPC$ formed by them.  We do so sequentially, by first defining weak projection categories (Section \ref{subsect:WPC}), projection categories (Section \ref{subsect:PC}), weak chained projection categories (Section \ref{subsect:WCPC}), and finally chained projection categories (Section \ref{subsect:CPC}).  At each step of the way we will show that a DRC-semigroup gives rise to a structure of the relevant kind, ultimately leading to a functor $\bC:\DRC\to\CPC$; see Theorem~\ref{thm:Cfunctor}.  We will eventually see in Section \ref{sect:iso} that $\bC$ is an isomorphism.

\subsection{Weak projection categories}\label{subsect:WPC}

\begin{defn}
A \emph{weak projection category} is a pair $(P,\C)$, consisting of a biordered category $\C=(\C,\leql,\leqr)$ and a projection algebra $P = v\C$, for which the restriction of both orders~$\leql$ and~$\leqr$ to $P$ is the order $\leq$ from \eqref{eq:leqP}.
\end{defn}

Let $(P,\C)$ be a weak projection category, and let $a\in\C$.  As in \eqref{eq:vtvd}, we have the two maps
\[
\vt_a:\bd(a)^\da \to \br(a)^\da \AND \vd_a:\br(a)^\da\to\bd(a)^\da,
\]
given by $p\vt_a = \br({}_p\corest a)$ and $q\vd_a = \bd(a\rest_q)$, for $p\leq\bd(a)$ and $q\leq\br(a)$.
Since $\im(\th_{\bd(a)}) = \bd(a)^\da$ and $\im(\de_{\br(a)}) = \br(a)^\da$, we can also define the maps
\begin{equation}\label{eq:ThDe}
\Th_a = \th_{\bd(a)}\vt_a:P \to \br(a)^\da \AND \De_a = \de_{\br(a)}\vd_a:P\to\bd(a)^\da.
\end{equation}

\begin{lemma}\label{lem:Thpa}
If $(P,\C)$ is a weak projection category, and if $a\in\C$, then
\ben
\item \label{Thpa1} $\vt_a\th_{\br(a)} = \vt_a$ and $\vd_a\de_{\bd(a)} = \vd_a$,
\item \label{Thpa2} $\Th_a\th_{\br(a)} = \Th_a$ and $\De_a\de_{\bd(a)} = \De_a$,
\item \label{Thpa3} $\Th_{{}_p\corest a} = \th_p\Th_a$ and $\De_{a\rest_q} = \de_q\De_a$ for any $p\leq\bd(a)$ and $q\leq\br(a)$,
\item \label{Thpa4} $\vt_a = \Th_a|_{\bd(a)^\da}$ and $\vd_a = \De_a|_{\br(a)^\da}$.
\een
\end{lemma}

\pf
\firstpfitem{\ref{Thpa1}}  For any $t\in\dom(\vt_a) = \bd(a)^\da$, we have $t\vt_a\leq\br(a)$, and so $t\vt_a = (t\vt_a)\th_{\br(a)}$.

\pfitem{\ref{Thpa2}}  This follows from \ref{Thpa1}, as $\Th_a = \th_{\bd(a)}\vt_a$ and $\De_a = \de_{\br(a)}\vd_a$.

\pfitem{\ref{Thpa3}}  For any $t\in P$ we have
\begin{align*}
t\Th_{{}_p\corest a} &= t\th_p\vt_{{}_p\corest a} &&\text{by definition}\\
&= t\th_p\vt_a &&\text{by \eqref{eq:vtvd2}}\\
&= t\th_p\th_{\bd(a)}\vt_a &&\text{by Lemma \ref{lem:thpthq}, as $p\leq\bd(a)$}\\
&= t\th_p\Th_a &&\text{by definition.}
\end{align*}
\firstpfitem{\ref{Thpa4}}  The maps $\vt_a$ and $\Th_a|_{\bd(a)^\da}$ both have domain $\bd(a)^\da$, and for $p\in\bd(a)^\da$ we have $p = p\th_{\bd(a)}$ by definition, and so $p\Th_a = p\th_{\bd(a)}\vt_a = p\vt_a$.
\epf

\begin{lemma}\label{lem:PSCS}
If $S$ is a DRC-semigroup, then $(\bP(S),\C(S))$ is a weak projection category.
\end{lemma}

\pf
Write $P=\bP(S)$ and $\C=\C(S)$.
By Lemma \ref{lem:PS} and Proposition \ref{prop:CS}, $P$ is a projection algebra and $\C$ a biordered category, and by construction we have $P=v\C$.  It remains to check that:
\[
p\leql q \Iff p\leq q \Iff p\leqr q \qquad\text{for all $p,q\in P$.}
\]
As ever, it is enough to establish the equivalence involving $\leql$ and $\leq$.  For this, first suppose $p\leql q$, so that $p = {}_r\corest q = rq$ for some $r\leq q$.  It follows that $p=pq$, so that $p\leq q$ (cf.~\eqref{eq:leqPS}).
Conversely, if $p\leq q\ (=\bd(q))$, then the restriction ${}_p\corest q$ exists, and ${}_p\corest q = pq = p$; it follows that~$p\leql q$.
\epf

In what follows, we will need to understand the $\vt/\Th$ and $\vd/\De$ maps associated to the weak projection category $(P,\C) = (\bP(S),\C(S))$ arising from a DRC-semigroup $S$.  For this, consider a morphism $a\in\C\ (=S)$.  For $p\leq\bd(a) = D(a)$ we have
\[
p\vt_a = \br({}_p\corest a) = R(pa) . 
\]
It follows that for arbitrary $q\in P$ we have
\[
q\Th_a = q\th_{\bd(a)}\vt_a = R(qD(a))\vt_a = R(R(qD(a))a) =_2 R(qD(a)a) =_1 R(qa) .
\]
Similar calculations apply to the $\vd$ and $\De$ maps, and in summary we have:
\begin{align}
\nonumber p\vt_a &= R(pa) \text{ \ for $p\leq\bd(a)$,} & p\Th_a = R(pa) \text{ \ for $p\in P$,}\\
\label{eq:Rpa} p\vd_a &= D(ap) \text{ \ for $p\leq\br(a)$,} & p\De_a = D(ap) \text{ \ for $p\in P$.}
\end{align}

\subsection{Projection categories}\label{subsect:PC}

We saw in Lemma \ref{lem:Thpa}\ref{Thpa3} that the $\Th$ and $\De$ maps in a weak projection category $(P,\C)$ behave well with regard to left and right restrictions, respectively, in the sense that $\Th_{{}_p\corest a} = \th_p\Th_a$ and $\De_{a\rest_q} = \de_q\De_a$ for any $p\leq\bd(a)$ and $q\leq\br(a)$.  We are particularly interested in the situation in which the analogous statements hold for the other restrictions.

\begin{defn}\label{defn:PC}
A \emph{projection category} is a weak projection category $(P,\C)$ satisfying:
\begin{enumerate}[label=\textup{\textsf{(C\arabic*)}},leftmargin=10mm]
\item \label{C1} For every $a\in\C$, and every $p\leq\bd(a)$ and $q\leq\br(a)$, we have
\[
\Th_{a\rest_q} = \Th_a\th_q \AND \De_{{}_p\corest a} = \De_a\de_p.
\]
\end{enumerate}
We write $\PC$ for the category of projection categories.  A morphism $(P,\C)\to(P',\C')$ in $\PC$ is a biordered morphism $\phi:\C\to\C'$ whose object map $v\phi:P\to P'$ is a projection algebra morphism.
\end{defn}

\begin{prop}\label{prop:PSCS}
The assignment $S\mt(\bP(S),\C(S))$ is the object part of a functor $\DRC\to\PC$.
\end{prop}

\pf
We first check that $(P,\C) = (\bP(S),\C(S))$ is a projection category for any DRC-semigroup~$S$.  By Lemma \ref{lem:PSCS}, it remains to check that \ref{C1} holds.
By symmetry, it suffices to show that $\Th_{a\rest_q} = \Th_a\th_q$ for all $a\in\C$, and all $q\leq\br(a)$.  For this we fix $t\in P$, and use \eqref{eq:Rpa} to calculate
\[
t\Th_{a\rest_q} = t\Th_{aq} = R(taq) =_2 R(R(ta)q) = R(ta)\th_q = t\Th_a\th_q.
\]

Now we know that $\bC$ maps objects of $\DRC$ to objects of $\PC$.
It remains to check that any DRC-morphism $\phi:S\to S'$ is a projection category morphism from $(P,\C)=(\bP(S),\C(S))$ to $(P',\C') = (\bP(S'),\C(S'))$, i.e.~that:
\bit
\item $\phi$ is a biordered category morphism $\C\to\C'$, and
\item $v\phi$ is a projection algebra morphism $P\to P'$.
\eit
But these follow from Propositions \ref{prop:Cfunc} and \ref{prop:Pfunc}, respectively.
\epf

\subsection{Weak chained projection categories}\label{subsect:WCPC}

We now bring in an extra layer of structure.

\begin{defn}
Given a projection category $(P,\C)$, an \emph{evaluation map} is a biordered $v$-functor $\ve:\CC\to\C$, where $\CC=\CC(P)$ is the chain category of $P$, and where \emph{$v$-functor} means that $\ve(p) = p$ for all $p\in P=v\C$.  (We write evaluation maps to the left of their arguments.)  It quickly follows that $\bd(\ve(\c)) = \bd(\c)$ and $\br(\ve(\c)) = \br(\c)$ for all $\c\in\CC$.
\end{defn}

Within the image of $\ve$, we will mainly be interested in elements of the form $\ve[p,q]$, for $(p,q)\in{\F}$.  The next lemma gathers some simple properties of these elements, all of which follow quickly from the definitions.

\begin{lemma}\label{lem:vepq}
Let $(P,\C)$ be a projection category, and $\ve:\CC\to\C$ an evaluation map.
\ben
\item \label{vepq1} For all $p\in P$ we have $\ve[p,p] = p$.
\item \label{vepq2} For all $(p,q)\in{\F}$ we have $\bd(\ve[p,q]) = p$ and $\br(\ve[p,q]) = q$,
\item \label{vepq3} For all $(p,q)\in{\F}$, and for all $r\leq p$ and $s\leq q$, we have
\[
{}_r\corest\ve[p,q] = \ve[r,r\th_q] \AND \ve[p,q]\rest_s = \ve[s\de_p,s].  \epfreseq
\]
\een
\end{lemma}

\begin{defn}\label{defn:WCPC}
A \emph{weak chained projection category} is a triple $(P,\C,\ve)$, where $(P,\C)$ is a projection category, and $\ve:\C\to\C$ is an \emph{evaluation map}.  We write $\WCPC$ for the category of weak chained projection categories.  A morphism $(P,\C,\ve)\to(P',\C',\ve')$ in $\WCPC$ is called a \emph{chained projection functor}, and is a projection category morphism $\phi:(P,\C)\to(P',\C')$ that respects the evaluation maps, in the sense that 
\[
(\ve[p_1,\ldots,p_k])\phi = \ve'[p_1\phi,\ldots,p_k\phi] \qquad\text{for $p_1,\ldots,p_k\in P$ with $p_1\F\cdots\F p_k$.}
\]
\end{defn}

Consider a DRC-semigroup $S$, and write $(P,\C)=(\bP(S),\C(S))$.  Since projections are idempotents, there is a well-defined map
\[
\ve = \ve(S):\CC=\CC(P) \to \C\ (=S) \GIVENBY \ve[p_1,\ldots,p_k] = p_1\cdots p_k,
\]
where the latter product is taken in $S$.

\begin{lemma}\label{lem:veS}
$\ve=\ve(S)$ is an evaluation map.
\end{lemma}

\pf
We first claim that for $\c = [p_1,\ldots,p_k]\in\P$ we have
\[
\bd(\ve(\c)) = p_1 = \bd(\c) \AND \br(\ve(\c)) = p_k = \br(\c),
\]
i.e.~that
\begin{equation}\label{eq:Dp1k}
D(p_1\cdots p_k) = p_1 \AND R(p_1\cdots p_k) = p_k \qquad\text{when $p_1\F\cdots\F p_k$.}
\end{equation}
By symmetry it suffices to prove the first.  The $k=1$ case being clear, we assume that $k\geq2$.  We then use \ref{DRC2}, induction and $p_1\F p_2$ to calculate
\[
D(p_1p_2\cdots p_k) = D(p_1D(p_2\cdots p_k)) = D(p_1p_2) = p_2\de_{p_1} = p_1.
\]

Next we check that $\ve$ is a $v$-functor.  Certainly $\ve(p) = p$ for all $p\in P = v\C$.  Now consider composable chains $\c=[p_1,\ldots,p_k]$ and $\d=[p_k,\ldots,p_l]$ in $\CC$.  We then have
\[
\ve(\c\circ\d) = \ve[p_1,\ldots,p_k,\ldots,p_l] = p_1\cdots p_k\cdots p_l = p_1\cdots p_k\cdot p_k\cdots p_l = \ve(\c)\cdot\ve(\d) = \ve(\c)\circ\ve(\d).
\]
Note that the final product is a composition in $\C$ since $\br(p_1\cdots p_k) = p_k = \bd(p_k\cdots p_l)$ by \eqref{eq:Dp1k}.

Finally, we need to check that $\ve$ is biordered, i.e.~that
\[
\ve({}_q\corest\c) = {}_q\corest\ve(\c) \ANd \ve(\c\rest_r) = \ve(\c)\rest_r \qquad\text{for all $\c\in\CC$, and for all $q\leq\bd(\c)$ and $r\leq\br(\c)$.}
\]
For the first (the second is dual), write $\c = [p_1,\ldots,p_k]$, so that ${}_q\corest\c = [q_1,\ldots,q_k]$, where each $q_i = q\th_{p_2}\cdots\th_{p_i}$.  Noting that $\ve({}_q\corest\c) = q_1\cdots q_k$ and $ {}_q\corest\ve(\c) = q\cdot p_1\cdots p_k$, we need to show that
\[
q\cdot p_1\cdots p_k = q_1\cdots q_k.  
\]
When $k=1$ we have $q\cdot p_1 = q = q_1$ (as $q\leq p_1$), so now suppose $k\geq2$.  We then have $q\cdot p_1\cdots p_{k-1}p_k = q_1\cdots q_{k-1}\cdot p_k$ by induction, so it remains to show that $q_{k-1}p_k = q_{k-1}q_k$.  For this we use Lemma \ref{lem:Dpa} to calculate
\[
q_{k-1}p_k =_1 q_{k-1}p_kR(q_{k-1}p_k) = q_{k-1}R(q_{k-1}p_k) = q_{k-1}\cdot q_{k-1}\th_{p_k} = q_{k-1}q_k.  \qedhere
\]
\epf

\begin{prop}\label{prop:PSCSveS}
The assignment $S\mt(\bP(S),\C(S),\ve(S))$ is the object part of a functor $\DRC\to\WCPC$.  
\end{prop}

\pf
It follows from Proposition \ref{prop:PSCS} and Lemma \ref{lem:veS} that $(\bP(S),\C(S),\ve(S))$ is a weak chained projection category for any DRC-semigroup $S$.  It remains to check that any DRC-morphism $\phi:S\to S'$ is a chained projection functor $(P,\C,\ve)\to(P',\C',\ve')$.  We already know from Proposition~\ref{prop:PSCS} that $\phi$ is a projection category morphism $(P,\C)\to(P',\C')$, so it remains to check that $\phi$ preserves the evaluation maps.  But for $p_1,\ldots,p_k\in P$ with $p_1\F\cdots\F p_k$ we have
\[
(\ve[p_1,\ldots,p_k])\phi = (p_1\cdots p_k)\phi = (p_1\phi)\cdots(p_k\phi) = \ve'[p_1\phi,\ldots,p_k\phi].  \qedhere
\]
\epf

\subsection{Chained projection categories}\label{subsect:CPC}

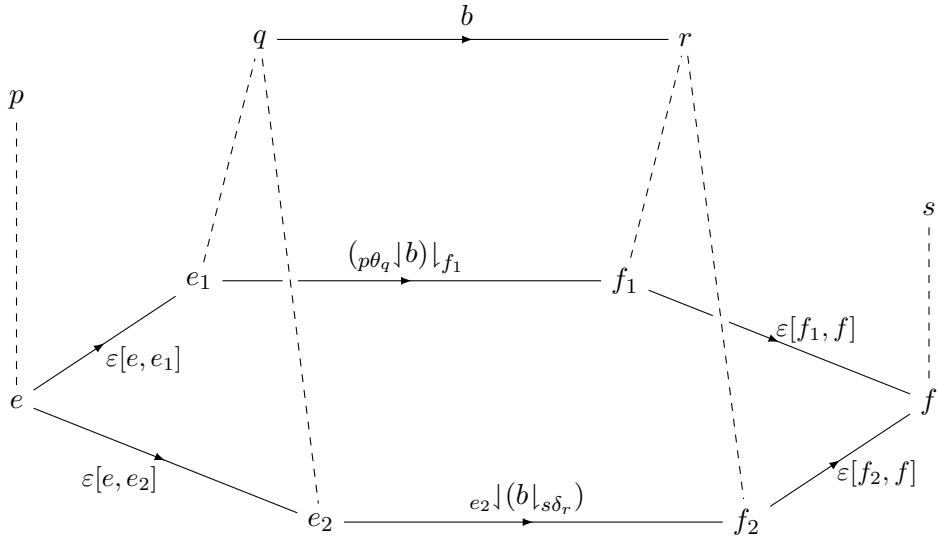
\begin{figure}[h]
\begin{center}
\begin{tikzpicture}[scale=0.8]
\tikzstyle{vertex}=[circle,draw=black, fill=white, inner sep = 0.07cm]
\node (e) at (0,0){$e$};
\node (e1) at (3,2){$e_1$};
\node (e2) at (5,-2){$e_2$};
\node (f) at (15,0){$f$};
\node (f1) at (10,2){$f_1$};
\node (f2) at (12,-2){$f_2$};
\node (p) at (0,5){$p$};
\node (q) at (4,6){$q$};
\node (r) at (11,6){$r$};
\node (s) at (15,3.2){$s$};
\draw[->-=0.5] (q)--(r);
\draw[->-=0.5] (e)--(e1);
\draw[->-=0.5] (e1)--(f1);
\draw[->-=0.5] (f1)--(f);
\draw[->-=0.5] (e)--(e2);
\draw[->-=0.5] (e2)--(f2);
\draw[->-=0.5] (f2)--(f);

\draw[white,line width=2mm] (q)--(e2);
\draw[white,line width=2mm] (r)--(f2);
\draw[dashed] (e1)--(q)--(e2) (f1)--(r)--(f2);
\draw[dashed] (e)--(p) (f)--(s);
\node () at (7.4,6.4) {$b$};
\node () at (6.4,2.4) {$({}_{p\th_q}\corest b)\rest_{f_1}$};
\node () at (8.4,-1.6) {${}_{e_2}\corest(b\rest_{s\de_r})$};
\node () at (2.1,0.7) {{\small $\ve[e,e_1]$}};
\node () at (1.7,-1.3) {{\small $\ve[e,e_2]$}};
\node () at (13.15,1.2) {{\small $\ve[f_1,f]$}};
\node () at (14.15,-1.2) {{\small $\ve[f_2,f]$}};
\end{tikzpicture}
\caption{The projections and morphisms associated to $b\in\C(q,r)$ and $p,s\in P$; see Definitions~\ref{defn:ef},~\ref{defn:lamrho} and~\ref{defn:CPC}, and Lemma \ref{lem:ef}.  Dashed lines indicate~$\leq$ relationships.  Axiom \ref{C2} says that the hexagon at the bottom of the diagram commutes.}
\label{fig:lamrho}
\end{center}
\end{figure}

We have almost achieved the main objective of this section, defining chained projection categories.  These will be the weak chained projection categories $(P,\C,\ve)$ satisfying a natural coherence condition stating that certain diagrams in $\C$ commute.  A diagrammatic representation of this condition can be seen in Figure \ref{fig:lamrho}.  To see that the relevant morphisms are well defined requires a preliminary definition and lemma.

In what follows, for a small category $\C$, and objects $p,q\in v\C$, it will be convenient to write
\[
\C(p,q) = \set{a\in\C}{\bd(a) = p ,\ \br(a) = q}
\]
for the set of all morphisms $p\to q$.

\begin{defn}\label{defn:ef}
Let $(P,\C)$ be a projection category, let $b\in\C(q,r)$, and let $p,s\in P$.  We define the projections
\begin{align}
\nonumber e(p,b,s) &= s\De_b\de_p , & e_1(p,b,s) &= s\De_b\de_{p\th_q} , & f_1(p,b,s) &= s\de_{p\Th_b} , & f(p,b,s) &= p\Th_b\th_s , \\
\label{eq:ef} & & e_2(p,b,s) &= p\th_{s\De_b} , & f_2(p,b,s) &= p\Th_b\th_{s\de_r}.
\end{align}
When no confusion arises, we will abbreviate these to $e = e(p,b,s)$, $e_1 = e_1(p,b,s)$, and so on.
\end{defn}

\begin{lemma}\label{lem:ef}
Let $(P,\C)$ be a projection category, and let $b\in\C(q,r)$ and $p,s\in P$.  Then the projections in~\eqref{eq:ef} satisfy the following:
\ben\bmc2
\item \label{ef1} $e\leq p$ and $e_1,e_2\leq q$,
\item \label{ef2} $e\F e_1,e_2$,
\item \label{ef3} $({}_{p\th_q}\corest b)\rest_{f_1}$ exists, and has domain $e_1$,
\item \label{ef4} $f\leq s$ and $f_1,f_2\leq r$,
\item \label{ef5} $f_1,f_2\F f$,
\item \label{ef6} ${}_{e_2}\corest(b\rest_{s\de_r})$ exists, and has range $f_2$.
\emc\een
\end{lemma}

\pf
\firstpfitem{\ref{ef1}}  We have
\[
e = s\De_b\de_p \leq p \COMMA e_1 = s\De_b\de_{p\th_q} \leq p\th_q \leq q \AND e_2 = p\th_{s\De_b} \leq s\De_b \leq \bd(b) = q.
\]
\firstpfitem{\ref{ef2}}  For $e\F e_2$, we have
\[
e\th_{e_2} = s\De_b\de_p\th_{p\th_{s\De_b}} =_{10} s\De_b\de_p\th_{s\De_b} =_9 p\th_{s\De_b} = e_2,
\]
and
\[
e_2\de_e = p\th_{s\De_b} \de_{s\De_b\de_p} =_{10} p\th_{s\De_b} \de_p =_9 s\De_b\de_p = e.
\]
For $e\F e_1$, let $p' = q\de_p$ and $q' = p\th_q$, so that $p'\F q'$ by Lemma \ref{lem:p'q'}\ref{p'q'1}.  
Since $e_1 = s\De_b\de_{q'} \leq q'$, it then follows from Corollary \ref{cor:pqrs} that $e_1\de_{p'} \F e_1$, so we can complete the proof of this part by showing that $e_1\de_{p'} = e$.  For this we use Lemmas \ref{lem:p'q'}\ref{p'q'2} and \ref{lem:Thpa}\ref{Thpa2} to calculate
\[
e_1\de_{p'} = s\De_b\de_{q'}\de_{p'} = s\De_b\de_q\de_p = s\De_b\de_p = e.
\]
\firstpfitem{\ref{ef3}}  As in the previous part we write $q' = p\th_q$, noting that $q' \leq q = \bd(b)$.  This means that ${}_{q'}\corest b$ exists, and we denote this by $b' = {}_{q'}\corest b$.  Note that
\[
\br(b') = q' \vt_b = p\th_q\vt_b = p\Th_b.
\]
Since $f_1 = s\de_{p\Th_b} \leq p\Th_b = \br(b')$, it follows that $b'\rest_{f_1}$ exists, and it remains to show that this has domain $e_1$.  For this we have
\begin{align*}
\bd(b'\rest_{f_1}) = f_1\vd_{b'} = s\de_{p\Th_b}\vd_{b'} &= s\De_{b'} &&\text{as $p\Th_b = \br(b')$}\\
&= s\De_b\de_{q'} = e_1 &&\text{by \ref{C1}, as $b' = {}_{q'}\corest b$.} 
\end{align*}
\firstpfitem{\ref{ef4}--\ref{ef6}}  These are dual to \ref{ef1}--\ref{ef2}.
\epf

\begin{defn}\label{defn:lamrho}
Let $(P,\C)$ be a projection category, and let $b\in\C(q,r)$ and $p,s\in P$.  Then with the projections in~\eqref{eq:ef}, it follows from Lemma \ref{lem:ef} that $\C$ contains the following well-defined morphisms:
\begin{equation}\label{eq:lamrho}
\lam(p,b,s) = \ve[e,e_1] \circ ({}_{p\th_q}\corest b)\rest_{f_1} \circ \ve[f_1,f]
\AND
\rho(p,b,s) = \ve[e,e_2] \circ {}_{e_2}\corest(b\rest_{s\de_r}) \circ \ve[f_2,f].
\end{equation}
These morphisms are shown in Figure \ref{fig:lamrho}.
\end{defn}

Here then is the main definition of this section:

\begin{defn}\label{defn:CPC}
A \emph{chained projection category} is a weak chained projection category $(P,\C,\ve)$ satisfying the following coherence condition:
\begin{enumerate}[label=\textup{\textsf{(C\arabic*)}},leftmargin=10mm] \addtocounter{enumi}1
\item \label{C2} For every $b\in\C$, and for every $p,s\in P$, we have $\lam(p,b,s) = \rho(p,b,s)$, where these morphisms are as in \eqref{eq:lamrho}.
\end{enumerate}
We denote by $\CPC$ the full subcategory of $\WCPC$ consisting of all chained projection categories, and all chained projection functors between them.
\end{defn}

For a DRC-semigroup $S$ we write $\bC(S) = (\bP(S),\C(S),\ve(S))$ for the weak chained projection category from Proposition \ref{prop:PSCSveS}.  

\begin{thm}\label{thm:Cfunctor}
$\bC$ is a functor $\DRC\to\CPC$.
\end{thm}

\pf
During the proof, we use \eqref{eq:Rpa} freely.

Given Proposition \ref{prop:PSCSveS}, we just need to check that $(P,\C,\ve) = \bC(S)$ satisfies~\ref{C2} for any DRC-semigroup $S$.  So fix $b\in\C(q,r)$ and $p,s\in P$, let $e,e_1,e_2,f,f_1,f_2\in P$ be as in \eqref{eq:ef}, and write
\[
b' = (_{p\th_q}\corest b)\rest_{f_1} \AND b'' = {}_{e_2}\corest(b\rest_{s\de_r}).
\]
By Lemma \ref{lem:ef}\ref{ef3} we have $D(b') = e_1$, and of course $R(b') = f_1$, so it follows from \ref{DRC1} that $b' = e_1b'f_1$, and so
\[
\lam(p,b,s) = ee_1 \cdot b' \cdot f_1f = eb'f \ANDSIM \rho(p,b,s) = eb''f.
\]
It therefore remains to show that $eb'f = eb''f$, and we claim that
\begin{equation}\label{eq:ebf}
eb'f = ebf \AND eb''f = ebf.
\end{equation}
It suffices by symmetry to prove the first.  For this we first note that $b' = ({}_{p\th_q}\corest b)\rest_{f_1} = R(pq)bf_1$, and we have
\begin{align}
\nonumber eb' &= ep\cdot R(pq)bf_1 &&\text{since $e\leq p$ by Lemma \ref{lem:ef}\ref{ef1}; cf.~\eqref{eq:leqPS}}\\
\nonumber &= ep\cdot qR(pq)bf_1 &&\text{by Lemma \ref{lem:Dpa}}\\
\nonumber &= epqbf_1 &&\text{by \ref{DRC1}}\\
\label{eq:eb'} &= epbf_1 &&\text{by \ref{DRC1}, since $q=D(b)$.}
\intertext{Next we note that}
\nonumber f_1f &= f_1sf &&\text{since $f\leq s$ by Lemma \ref{lem:ef}\ref{ef4}; cf.~\eqref{eq:leqPS}}\\
\nonumber &= D(R(pb)s)\cdot sf &&\text{by the definition of $f_1$}\\
\nonumber &= D(R(pb)s)R(pb)\cdot sf &&\text{by Lemma \ref{lem:Dpa}}\\
\nonumber &= R(pb)sf &&\text{by \ref{DRC1}}\\
\label{eq:f1f} &= R(pb)f &&\text{using $f\leq s$ again.}
\intertext{Putting everything together, it follows that}
\nonumber eb'f &= epbf_1f &&\text{by \eqref{eq:eb'}}\\
\nonumber &= epbR(pb)f &&\text{by \eqref{eq:f1f}}\\
\nonumber &= epbf &&\text{by \ref{DRC1}}\\
\nonumber &= ebf &&\text{as $e\leq p$.}
\end{align}
This completes the proof of \eqref{eq:ebf}, and hence of the theorem.
\epf

\section{From chained projection categories to DRC-semigroups}\label{sect:CtoS}

In Section \ref{sect:CPC} we defined the category $\CPC$ of chained projection categories, and constructed a functor $\bC:\DRC\to\CPC$.  We will see in Section \ref{sect:iso} that $\bC$ is an isomorphism.  The proof involves constructing an inverse functor $\bS:\CPC\to\DRC$ in the opposite direction, and that is the purpose of this section.  The main part of the work is undertaken in Section \ref{subsect:S}, where we show how to construct a DRC semigroup from a chained projection category; see Definition \ref{defn:pr} and Theorem \ref{thm:S}.  We then construct the functor $\bS$ in Section \ref{subsect:Sfunctor}; see Theorem \ref{thm:Sfunctor}.

\subsection{The DRC-semigroup associated to a chained projection category}\label{subsect:S}

\begin{figure}[h]
\begin{center}
\begin{tikzpicture}[scale=.8]
\tikzstyle{vertex}=[circle,draw=black, fill=white, inner sep = 0.07cm]
\node (p) at (5,2.5){$p$};
\node (q) at (10,3.5){$q$};
\node (p0) at (0,2.5){};
\node (q0) at (15,3.5){};
\node (p') at (5,0){$p'$};
\node (q') at (10,0){$q'$};
\node (p0') at (0,0){};
\node (q0') at (15,0){};
\draw[dashed] (p)--(p') 
(q)--(q') 
(p0.south east)--(p0'.north east)
(q0.south west)--(q0'.north west)
;
\draw[->-=0.5] (p0)--(p);
\draw[->-=0.5] (q)--(q0);
\draw[->-=0.5] (p0')--(p');
\draw[->-=0.5] (q')--(q0');
\draw[->-=0.5] (p')--(q');
\draw[->-=0.5] (p0'.east) to [bend right = 15] (q0'.west);
\node () at (2.4,3.2-.5+.1) {$a$};
\node () at (12.5,4.25-.5+.1) {$b$};
\node () at (2.4,.3+.1) {$a\rest_{p'}$};
\node () at (12.5,.25+.1) {${}_{q'}\corest b$};
\node () at (7.4,.3+.1) {$\ve[p',q']$};
\node () at (7.35,-1.35-.1) {$a\pr b$};
\end{tikzpicture}
\caption{Formation of the product $a\pr b$; see Definition \ref{defn:pr}.  Dashed lines indicate~$\leq$ relationships.}
\label{fig:pr}
\end{center}
\end{figure}
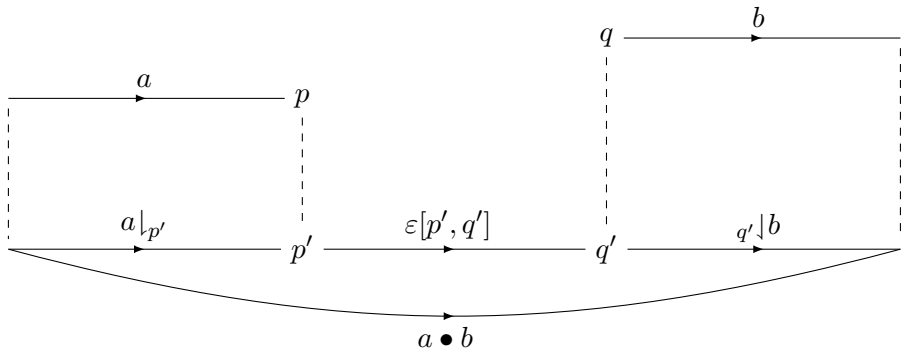

\begin{defn}\label{defn:pr}
Given a chained projection category $(P,\C,\ve)$, we define $\bS(P,\C,\ve) = (S,\pr,D,R)$ to be the $(2,1,1)$-algebra with:
\bit
\item underlying set $S = \C$,
\item unary operations $D$ and $R$ given by $D(a) = \bd(a)$ and $R(a) = \br(a)$ for $a\in \C$, and
\item binary operation $\pr$ defined for $a,b\in\C$, with $p=\br(a)$ and $q=\bd(b)$, by
\[
a\pr b = a\rest_{p'} \circ \ve[p',q'] \circ {}_{q'}\corest b, \WHERE p' = q\de_p \ANd q' = p\th_q.
\]
These elements and compositions are well defined because of Lemma \ref{lem:p'q'}\ref{p'q'1}.  See Figure \ref{fig:pr} for a diagrammatic representation of the product $a\pr b$.
\eit
\end{defn}

For the remainder of this section we fix a chained projection category $(P,\C,\ve)$.  In Theorem~\ref{thm:S} we will show that $\bS(P,\C,\ve)$ is indeed a DRC-semigroup.  Our first task in this direction is to show that~$\pr$ is associative, the proof of which requires the following lemma.

\begin{lemma}\label{lem:draprb}
For any $a,b\in\C$ we have
\[
\bd(a\pr b) = \bd(b)\De_a \AND \br(a\pr b) = \br(a)\Th_b.
\]
\end{lemma}

\pf
For the first (the second is dual), we keep the notation of Definition \ref{defn:pr}, and we have $\bd(a\pr b) = \bd(a\rest_{p'}) = p'\vd_a = q\de_p\vd_a = q\De_a = \bd(b) \De_a$.
\epf

\begin{lemma}\label{lem:abc}
For any $a,b,c\in\C$ we have $(a\pr b)\pr c = a\pr (b\pr c)$.
\end{lemma}

\pf
During the proof, we write 
\[
p = \br(a) \COMMA q = \bd(b) \COMMA r = \br(b) \AND s = \bd(c).
\]
We also denote the projections from Definition \ref{defn:ef} by $e=e(p,b,s)$, $e_1=e_1(p,b,s)$, and so on.  Given \ref{C2}, we can prove the lemma by showing that
\[
(a\pr b)\pr c = a\rest_e \circ \lam \circ {}_f\corest c
\AND
a\pr (b\pr c) = a\rest_e \circ \rho \circ {}_f\corest c,
\]
where $\lam=\lam(p,b,s)$ and $\rho=\rho(p,b,s)$ are as in Definition \ref{defn:lamrho}.  We just do the first, as the second follows by duality.

Writing $p'=q\de_p$ and $q'=p\th_q$, we first note that
\[
a\pr b = a' \circ g \circ b', \WHERE a' = a\rest_{p'}\COMMa g = \ve[p',q'] \ANd b' = {}_{q'}\corest b,
\]
and we also let $t = \br(a\pr b) = p\Th_b$ (cf.~Lemma \ref{lem:draprb}).  Then
\[
(a\pr b)\pr c = (a\pr b)\rest_{t'} \circ \ve[t',s'] \circ {}_{s'}\corest c, \WHERE t' = s\de_t \ \ \text{and} \ \ s' = t\th_s.
\]
We also use \ref{O4} to calculate
\[
(a\pr b)\rest_{t'} = (a' \circ g \circ b')\rest_{t'} = a'\rest_u \circ g\rest_v \circ b'\rest_{t'}, 
\]
where $v = \bd(b'\rest_{t'}) = t'\vd_{b'}$ and $u = \bd(g\rest_v) = v\vd_{g}$.  
Next we note that
\[
t' = s\de_t = s\de_{p\Th_b} = f_1 \AND s' = t\th_s = p\Th_b\th_s = f.
\]
Since $t = \br(a\pr b) = \br(a' \circ g \circ b') = \br(b')$, we also have
\[
v = t'\vd_{b'} = s\de_t\vd_{b'} = s\De_{b'} = s\De_{{}_{q'}\corest b} = s\De_b\de_{q'} = e_1,
\]
using \ref{C1} in the second-last step.  Combining this with Lemma \ref{lem:vepq}\ref{vepq3} it follows that
\[
g\rest_v = \ve[p',q']\rest_v = \ve[v\de_{p'},v].
\]
We then have
\begin{align*}
u = \bd(g\rest_v) = v\de_{p'} &= s\De_b\de_{q'}\de_{p'} &&\text{as shown above}\\
&= s\De_b\de_q\de_p &&\text{by Lemma \ref{lem:p'q'}\ref{p'q'2}}\\
&= s\De_b\de_p = e &&\text{by Lemma \ref{lem:Thpa}\ref{Thpa2}, noting that $q=\bd(b)$.}
\end{align*}
Putting everything together, we have
\begin{align*}
(a\pr b)\pr c = (a\pr b)\rest_{t'} \circ \ve[t',s'] \circ {}_{s'}\corest c &= a'\rest_u \circ g\rest_v \circ b'\rest_{t'} \circ \ve[t',s'] \circ {}_{s'}\corest c \\
&= a\rest_{p'}\rest_e \circ \ve[e,e_1] \circ ({}_{q'}\corest b)\rest_{f_1} \circ \ve[f_1,f] \circ {}_f\corest c
&= a\rest_e \circ \lam \circ {}_f\corest c,
\end{align*}
as required.
\epf

\begin{lemma}\label{lem:prcirc}
\ben
\item \label{prcirc1} If $a,b\in\C$ are such that $\br(a) = \bd(b)$, then $a\pr b = a\circ b$.
\item \label{prcirc2} If $p,q\in P$, then $p\pr q = \ve[p',q']$, where $p' = q\de_p$ and $q' = p\th_q$.
\item \label{prcirc3} If $a\in\C$, and if $p\leq\bd(a)$ and $q\leq\br(a)$, then $p\pr a = {}_p\corest a$ and $a\pr q = a\rest_q$.
\een
\end{lemma}

\pf
\firstpfitem{\ref{prcirc1}}  Keeping the notation of Definition \ref{defn:pr}, the assumption here is that $p=q$.  In this case we also have $p'=p = q = q'$ by \ref{P1}, and then
\[
a\pr b = a\rest_p\circ\ve[p,p]\circ{}_p\corest b = a\circ p \circ b = a\circ b.
\]
\firstpfitem{\ref{prcirc2}}  Noting that $p = \br(p)$ and $q = \bd(q)$, we have
\[
p\pr q = p\rest_{p'}\circ\ve[p',q']\circ{}_{q'}\corest q = p'\circ\ve[p',q']\circ q' = \ve[p',q'].
\]
\firstpfitem{\ref{prcirc3}}  For the first (the second is dual), we write $s = \bd(a)$, and we have
\[
p \pr a = p\rest_{p'} \circ \ve[p',s'] \circ {}_{s'}\corest a, \WHERE p' = s\de_p \ANd s' = p\th_s.
\]
But since $p\leq s$, we have $p' = s' = p$, so in fact
\[
p\pr a = p\rest_p \circ \ve[p,p] \circ {}_p\corest a = p\circ p\circ {}_p\corest a = {}_p\corest a.  \qedhere
\]
\epf

\begin{thm}\label{thm:S}
If $(P,\C,\ve)$ is a chained projection category, then $\bS(P,\C,\ve)$ is a DRC-semigroup.
\end{thm}

\pf
Given Lemma \ref{lem:abc}, it remains to verify \ref{DRC1}--\ref{DRC4}.  By symmetry, only the first part of each needs to be treated.

\pfitem{\ref{DRC1}}  We have $D(a)\pr a = \bd(a) \pr a = \bd(a)\circ a = a$, where we used Lemma \ref{lem:prcirc}\ref{prcirc1}, and the fact that $\br(\bd(a)) = \bd(a)$.

\pfitem{\ref{DRC2}}  Writing $p = \bd(b)$, we use Lemma \ref{lem:draprb} to calculate
\[
D(a\pr D(b)) = \bd(a\pr p) = \bd(p)\De_a = p\De_a = \bd(b)\De_a = \bd(a\pr b) = D(a\pr b).
\]
\firstpfitem{\ref{DRC3}}  Write $p = \bd(a)$ and $q = \bd(a\pr b)$.  We must show that $p\pr q = q = q\pr p$.  By Lemma~\ref{lem:prcirc}\ref{prcirc2}, we have
\[
p \pr q = \ve[p',q'], \WHERE p' = q\de_p \ANd q' = p\th_q.
\]
By Lemma \ref{lem:draprb}, and since $\im(\De_a)\sub p^\da$, we have $q = \bd(a\pr b) = \bd(b)\De_a \leq p$, so in fact $p'=q'=q$.  Thus,
\[
p \pr q = \ve[p',q'] = \ve[q,q] = q.
\]
Essentially the same argument gives $q\pr p = q$.

\pfitem{\ref{DRC4}}  This follows from $\br(\bd(a)) = \bd(a)$.
\epf

We will also need the following information concerning the projections of $\bS(P,\C,\ve)$.

\begin{prop}\label{prop:PSPCve}
Let $(P,\C,\ve)$ be a chained projection category, and let $S = \bS(P,\C,\ve)$.  Then
\ben
\item \label{PSPCve1} $\bP(S) = P$,
\item \label{PSPCve2} $D(p\pr q) = q\de_p$ and $R(p\pr q) = p\th_q$ for all $p,q\in P$.
\een
\end{prop}

\pf
\firstpfitem{\ref{PSPCve1}}  We have $\bP(S) = \set{D(a)}{a\in S} = \set{\bd(a)}{a\in\C} = v\C = P$.

\pfitem{\ref{PSPCve2}}  This follows from combining Lemma \ref{lem:prcirc}\ref{prcirc2} with Lemma \ref{lem:vepq}\ref{vepq2}.
\epf

\subsection{A functor}\label{subsect:Sfunctor}

We can think of the construction of the DRC-semigroup $\bS(P,\C,\ve)$ from the chained projection category $(P,\C,\ve)$ as an object map $\CPC\to\DRC$.  At the level of morphisms, any chained projection functor $\phi:(P,\C,\ve)\to(P',\C',\ve')$ is a map $\C\to\C'$.  Since the underlying sets of the semigroups $S=\bS(P,\C,\ve)$ and $S'=\bS(P',\C',\ve')$ are $\C$ and $\C'$, respectively, we can think of $\phi$ as a map $\phi=\bS(\phi):S\to S'$.  We now show that $\bS:\CPC\to\DRC$, interpreted in this way, is a functor.

\begin{thm}\label{thm:Sfunctor}
$\bS$ is a functor $\CPC\to\DRC$.
\end{thm}

\pf
It remains only to check that a morphism $\phi:(P,\C,\ve)\to(P',\C',\ve')$ in $\CPC$ is also a morphism $S\to S'$ in $\DRC$, where $S=\bS(P,\C,\ve)$ and $S'=\bS(P',\C',\ve')$.  For simplicity, we will write $\ol a = a\phi$ for $a\in S\ (=\C)$.  We also use dashes to distinguish between the various parameters and operations on $\C$ and $\C'$, and on $S$ and $S'$.

If $a\in S$, then writing $p = D(a) = \bd(a)$, we have
\[
D'(\ol a) = \bd'(\ol a) = \bd'(\ol{p\circ a}) = \bd'(\ol p\circ' \ol a) = \bd'(\ol p) = \ol p = \ol{D(a)} \ANDSIM R'(\ol a) = \ol{R(a)}.
\]
It therefore remains to show that $\ol{a\pr b}=\ol a\pr' \ol b$ for all $a,b\in S$.  So fix some such $a,b$, and let
\[
p = \br(a) \COMMA q = \bd(b) \COMMA p' = q\de_p \AND q' = p\th_q.
\]
As shown above, we have $\br'(\ol a) = \ol{\br(a)} = \ol p$, and similarly $\bd'(\ol b) = \ol q$.  It follows that
\[
a\pr b = a\rest_{p'} \circ \ve[p',q'] \circ {}_{q'}\corest b \AND 
\ol a \pr' \ol b = \ol a\rest_{\ol p'} \circ' \ve'[\ol p',\ol q'] \circ' {}_{\ol q'}\corest \ol b, 
\]
where $\ol p' = \ol q\de'_{\ol p}$ and $\ol q' = \ol p\th'_{\ol q}$.  Since $v\phi:P\to P'$ is a projection algebra morphism, we have
\[
\ol{p'} = \ol{q\de_p} = \ol q\de'_{\ol p} = \ol p' \ANDSIM \ol{q'} = \ol q'.
\]
Since $\phi$ respects evaluation maps, we have
\[
\ol{\ve[p',q']} = \ve'[\ol{p'},\ol{q'}] = \ve'[\ol p',\ol q'].
\]
Since $\phi$ is a biordered functor, we have
\[
\ol{a\rest_{p'}} = \ol a\rest_{\ol{p'}} = \ol a\rest_{\ol p'} \ANDSIM \ol{{}_{q'}\corest b} = {}_{\ol q'}\corest\ol b.
\]
Putting everything together, we have
\[
\ol{a\pr b} = \ol{a\rest_{p'} \circ \ve[p',q'] \circ {}_{q'}\corest b} = \ol{a\rest_{p'}} \circ' \ol{\ve[p',q']} \circ' \ol{{}_{q'}\corest b} = \ol a\rest_{\ol p'} \circ' \ve'[\ol p',\ol q'] \circ' {}_{\ol q'}\corest \ol b = \ol a\pr' \ol b.  \qedhere
\]
\epf

\section{The category isomorphism}\label{sect:iso}

In this section we will show that the functors
\[
\bC:\DRC\to\CPC \AND \bS:\CPC\to\DRC
\]
from Theorems \ref{thm:Cfunctor} and \ref{thm:Sfunctor} are mutually inverse isomorphisms, thereby proving the main result of the paper:

\begin{thm}\label{thm:iso}
The category $\DRC$ of DRC-semigroups (with DRC-morphisms) is isomorphic to the category $\CPC$ of chained projection categories (with chained projection functors).
\end{thm}

\pf
Since the functors $\bC$ and $\bS$ both act identically on morphisms, we just need to show that
\[
\bC(\bS(P,\C,\ve)) = (P,\C,\ve) \AND \bS(\bC(S)) = S
\]
for any chained projection category $(P,\C,\ve)$, and any DRC-semigroup $S$.  We do this in Propositions \ref{prop:CoS} and \ref{prop:SoC}.
\epf

\begin{prop}\label{prop:CoS}
For any chained projection category $(P,\C,\ve)$ we have $\bC(\bS(P,\C,\ve)) = (P,\C,\ve)$.
\end{prop}

\pf
Throughout the proof we will write
\[
S = \bS(P,\C,\ve) \AND (P',\C',\ve') = \bC(S) = (\bP(S),\C(S),\ve(S)).
\]
To prove the result we need to show that:
\ben
\item \label{CoS1} $P = P'$ as projection algebras,
\item \label{CoS2} $\C = \C'$ as biordered categories, and
\item \label{CoS3} $\ve = \ve'$ as maps.
\een
As usual we use dashes in what follows to denote the various parameters and operations on $\C'$.

\newpage  

\pfitem{\ref{CoS1}}  This follows immediately from Proposition \ref{prop:PSPCve}.

\pfitem{\ref{CoS2}}  By construction, the underlying set of $\C'$ is the same as that of $S$ and $\C$.  We also have $v\C' = \bP(S) = P = v\C$, and
\[
\bd'(a) = D(a) = \bd(a) \ANDSIM \br'(a) = \br(a) \qquad\text{for all $a\in\C'\ (=\C)$.}
\]
It follows that for any $a,b\in\C'$ with $\br'(a) = \bd'(b)$, we have $\br(a) = \bd(b)$, and Lemma \ref{lem:prcirc}\ref{prcirc1} then gives
\[
a\circ' b = a \pr b = a\circ b.
\]
It remains to check that the orders on $\C$ and $\C'$ coincide, and for this it suffices to show that
\[
{}_p\corestt a = {}_p\corest a \AND a\restt_q = a\rest_q \qquad\text{for all $a\in\C'\ (=\C)$, and all $p\leq\bd(a)$ and $q\leq\br(a)$.}
\]
Here we write $\corestt$ and $\restt$ for restrictions in $\C'$, and $\corest$ and $\rest$ for restrictions in $\C$.  For this we use Lemma \ref{lem:prcirc}\ref{prcirc3} to calculate
\[
{}_p\corestt a = p\pr a ={}_p\corest a \ANDSIM a\restt_q = a\rest_q.
\]
\firstpfitem{\ref{CoS3}}  We must show that $\ve'(\c) = \ve(\c)$ for any chain $\c\in\C(P)$.  So fix some such~${\c = [p_1,\ldots,p_k]}$.  First observe that for any $1\leq i<k$, Lemma \ref{lem:prcirc}\ref{prcirc2} gives 
\[
p_i\pr p_{i+1} = \ve[p_{i+1}\de_{p_i},p_i\th_{p_{i+1}}] = \ve[p_i,p_{i+1}], \qquad\text{as $p_i\F p_{i+1}$.}
\]
We then have
\begin{align*}
\ve'(\c) &= p_1\pr\cdots\pr p_k &&\text{by definition}\\
&= (p_1\pr p_2) \pr (p_2\pr p_3) \pr\cdots\pr (p_{k-1}\pr p_k) &&\text{as each $p_i$ is an idempotent of $S$}\\
&= \ve[p_1,p_2] \pr \ve[p_2,p_3] \pr\cdots\pr \ve[p_{k-1},p_k] &&\text{as just observed}\\
&= \ve[p_1,p_2] \circ \ve[p_2,p_3] \circ\cdots\circ \ve[p_{k-1},p_k] &&\text{by Lemma \ref{lem:prcirc}\ref{prcirc1}}\\
&= \ve[p_1,p_2,p_3,\ldots,p_k] = \ve(\c) &&\text{as $\ve$ is a functor.}
\end{align*}
This completes the proof.
\epf

\begin{prop}\label{prop:SoC}
For any DRC-semigroup $S$ we have $\bS(\bC(S)) = S$.
\end{prop}

\pf
This time we write
\[
(P,\C,\ve) = \bC(S) = (\bP(S),\C(S),\ve(S)) \AND S' = \bS(P,\C,\ve).
\]
Again the underlying sets of $S$ and $S'$ are the same, and so are the unary operations.  It therefore remains to show that
\[
a\pr b = ab \qquad\text{for all $a,b\in S'\ (=S)$.}
\]
For this, write $p = \br(a) = R(a)$ and $q = \bd(b) = D(b)$.  Then with $p' = q\de_p = D(pq)$ and $q' = p\th_q = R(pq)$, we have
\[
a\pr b = a\rest_{p'} \circ \ve[p',q'] \circ {}_{q'}\corest b = (ap') \cdot (p'q') \cdot (q'b) = ap'q'b.
\]
We then use Lemma \ref{lem:Dpa} and \ref{DRC1} to calculate
\begin{equation}\label{eq:DpqRpq}
p'q' = D(pq)\cdot R(pq) = D(pq)p \cdot qR(pq) = pq.
\end{equation}
Thus, continuing from above, we have
\[
a\pr b = ap'q'b = ap \cdot qb = ab,
\]
by \ref{DRC1}, as $p=R(a)$ and $q=D(b)$.
\epf

\section{Free and fundamental projection-generated DRC-semigroups}\label{sect:PGP}

Theorem \ref{thm:iso} shows that the categories $\DRC$ and $\CPC$ are isomorphic.  It also follows from the definition of the functors $\bC:\DRC\to\CPC$ and $\bS:\CPC\to\DRC$ that the algebras of projections of DRC-semigroups are precisely the same as the (structured) object sets of chained projection categories.  But the theorem does not answer the following:

\begin{que}\label{qu:P}
Given an abstract projection algebra $P$, as in Definition \ref{defn:P}, does there exist a DRC-semigroup $S$ with $\bP(S) = P$?  Or equivalently, does there exist a chained projection category~$(P,\C,\ve)$?
\end{que}

The semigroup version of this question was answered in the affirmative by Jones in his study of \emph{fundamental} DRC-semigroups \cite{Jones2021}, i.e.~the DRC-semigroups with no non-trivial projection-separating congruences.  In \cite[Section 9]{Jones2021}, Jones constructed a DRC-semigroup, denoted $C(P)$, from a projection algebra $P$, and showed it to be the \emph{maximum} fundamental DRC-semigroup with projection algebra $P$, meaning that:
\bit
\item $C(P)$ is fundamental, and has projection algebra $P = \bP(C(P))$, and
\item any fundamental DRC-semigroup with projection algebra $P$ embeds canonically into $C(P)$.
\eit

In a sense, one of the purposes of the current section is to provide another way to answer Question~\ref{qu:P}, in that we construct two different DRC-semigroups, denoted $\FP$ and $\MP$, both with projection algebra~$P$.  However, we have a deeper purpose than this (as did Jones), and our main objective here is to demonstrate the existence of \emph{free (projection-generated) DRC-semigroups}.  Specifically, we show that the assignment $P\mt\FP$ is the object part of a functor $\bF:\PA\to\DRC$, and that this is in fact a left adjoint to the forgetful functor $\bP:\DRC\to\PA$ from Proposition \ref{prop:Pfunc}, which maps a DRC-semigroup $S$ to its underlying projection algebra~$\bP(S)$.  It follows that $\FP$ is the free (projection-generated) DRC-semigroup with projection algebra $P$.

The semigroup $\FP$ is defined abstractly in terms of a presentation by generators and relations.  By contrast, the semigroup $\MP$ is a concrete semigroup consisting of pairs of self-maps of $P$, built from the underlying $\th_p$ and $\de_p$ operations of $P$.  The main initial purpose of $\MP$ is to provide a concrete homomorphic image of $\FP$, which will be useful in certain proofs concerning the latter.  But we will also see that $\MP$ is very special in its own right, in that it is the unique projection-generated fundamental DRC-semigroup with projection algebra $P$.

The definitions of $\FP$ and $\MP$ are given in Section \ref{subsect:FP}, and their key properties are established in Theorem \ref{thm:FPMP}.  Categorical freeness of $\FP$ is established in Section \ref{subsect:adjoint}; see Theorem \ref{thm:free}.  The unique fundamentality properties of $\MP$ are established in Section \ref{subsect:fund}; see Theorem \ref{thm:fund}.  

Before we begin, we briefly revise our notation for semigroup presentations.  Let $X$ be an alphabet, and write $X^+$ for the free semigroup over $X$, which is the semigroup of all non-empty words over~$X$, under concatenation.  Also let $R \sub X^+\times X^+$ be a set of pairs of words over $X$, and write $R^\sharp$ for the congruence on $X^+$ generated by $R$.  We then write
\[
\pres XR = X^+/R^\sharp,
\]
and call this the semigroup defined by the \emph{presentation} $\pres XR$.  Elements of $X$ are called \emph{letters} or \emph{generators}, and elements of $R$ are called \emph{relations}.  A relation $(u,v)\in R$ is often written as an equality,~$u=v$.

\newpage

\subsection[Definition of $\FP$ and $\MP$]{\boldmath Definition of $\FP$ and $\MP$}\label{subsect:FP}

Fix a projection algebra $P$ for the duration of this section.  Also fix an alphabet
\[
X_P = \set{x_p}{p\in P}
\]
in one-one correspondence with $P$, and let $R_P$ be the set of relations, quantified over all $p,q\in P$:
\begin{enumerate}[label=\textup{\textsf{(R\arabic*)}},leftmargin=9mm]\bmc3
\item \label{R1} $x_p^2 = x_p$,
\item \label{R2} $x_px_q = x_px_{p\th_q}$,
\item \label{R3} $x_px_q = x_{q\de_p}x_q$.
\emc
\end{enumerate}
Write ${\sim} = R_P^\sharp$ for the congruence on $X_P^+$ generated by the relations \ref{R1}--\ref{R3}, and denote by~$\ol w$ the $\sim$-class of a word $w\in X_P^+$.  We then take $\FP$ to be the semigroup defined by the presentation
\[
\FP = \pres{X_P}{R_P} = X_P^+/{\sim} = \set{\ol w}{w\in X_P^+}.
\]
We will soon see that $\FP$ is a DRC-semigroup with projection algebra (isomorphic to) $P$.

We also define a second semigroup associated to $P$.  First, we denote by $\T_P$ the \emph{full transformation semigroup} over $P$, i.e.~the semigroup of all maps $P\to P$ under composition.  As before, we also write $\T_P^\opp$ for the opposite semigroup.  Our second semigroup will be a subsemigroup of the direct product $\T_P\times\T_P^\opp$, in which the operation is given by $(\al,\al')(\be,\be') = (\al\be,\be'\al')$.  For each $p\in P$, we define the pair
\[
\wh p = (\th_p,\de_p) \in \T_P \times \T_P^\opp,
\]
and the semigroup
\[
\MP = \pres{\wh p}{p\in P} \leq \T_P\times\T_P^\opp.
\]
A typical element of $\MP$ has the form
\[
\wh p_1\cdots\wh p_k = (\th_{p_1}\cdots\th_{p_k},\de_{p_k}\cdots\de_{p_1}) \qquad\text{for $p_1,\ldots,p_k\in P$.}
\]
Define the (surjective) semigroup homomorphism
\[
\psi:X_P^+\to\MP:x_p\mt\wh p = (\th_p,\de_p).
\]

\begin{lemma}
We have $R_P \sub \ker(\psi)$.
\end{lemma}

\pf
We need to check that each relation from $R_P$ is preserved by $\psi$, in the sense that $u\psi=v\psi$ for each such relation $(u,v)\in R_P$.  To do so, let $p,q\in P$.  For \ref{R1} we have
\[
x_p^2\psi = \wh p\cdot\wh p = (\th_p\th_p,\de_p\de_p) = (\th_p,\de_p) = \wh p = x_p\psi,
\]
by \ref{P6}.  For \ref{R2} we use \ref{P3} and \ref{P5} to calculate:
\[
(x_px_q)\psi = \wh p\cdot\wh q = (\th_p\th_q,\de_q\de_p) = (\th_p\th_{p\th_q},\de_{p\th_q}\de_p) = \wh p \cdot\wh{p\th_q} = (x_px_{p\th_q})\psi.
\]
We deal with \ref{R3} in similar fashion.
\epf

It follows that $\psi$ induces a well-defined (surjective) semigroup homomorphism
\[
\Psi:\FP \to \MP \GIVENBY \ol w\Psi = w\psi \qquad\text{for $w\in X_P^+$.}
\]
Note in particular that $\ol x_p\Psi = \wh p = (\th_p,\de_p)$ for $p\in P$.  

In what follows, we write $\sim_1$ to denote $\sim$-equivalence of words from $X_P^+$ by one or more applications of relations from \ref{R1}, with similar meanings for $\sim_2$ and $\sim_3$.

\begin{lemma}\label{lem:xp'xq'}
If $p,q\in P$, then $x_px_q \sim x_{p'}x_{q'}$, where $p' = q\de_p$ and $q' = p\th_q$.
\end{lemma}

\pf
We have $x_px_q \sim_2 x_px_{p\th_q} \sim_3 x_{(p\th_q)\de_p} x_{p\th_q} = x_{q\de_p}x_{p\th_q}$, using \ref{P9} in the last step.
\epf

For a path $\p = (p_1,\ldots,p_k)$ from the path category $\P = \P(P)$, we define the word
\[
w_\p = x_{p_1}\cdots x_{p_k} \in X_P^+.
\]

\begin{lemma}\label{lem:wp}
Any word over $X_P$ is $\sim$-equivalent to $w_\p$ for some $\p\in\P$.
\end{lemma}

\pf
Consider an arbitrary word $w = x_{p_1}\cdots x_{p_k} \in X_P^+$.  We will show by induction on $k$ that
\[
w \sim x_{p_1'}\cdots x_{p_k'} \qquad\text{for some $p_i'\leq p_i$ with $p_1'\F\cdots\F p_k'$.}
\]
The $k=1$ case is trivial, so we now assume that $k\geq2$.  First, Lemma \ref{lem:xp'xq'} gives
\[
x_{p_{k-1}}x_{p_k} \sim x_{p_{k-1}''}x_{p_k''}, \WHERE p_{k-1}'' = p_k\de_{p_{k-1}} \ANd p_k'' = p_{k-1}\th_{p_k},
\]
and by Lemma \ref{lem:p'q'}\ref{p'q'1} we have $p_{k-1}''\leq p_{k-1}$, $p_k''\leq p_k$ and $p_{k-1}''\F p_k''$.  By induction, we have
\[
x_{p_1}\cdots x_{p_{k-2}} x_{p_{k-1}''} \sim x_{p_1'}\cdots x_{p_{k-2}'}x_{p_{k-1}'},
\]
for some $p_i'\leq p_i$ ($i=1,\ldots,k-2$) and $p_{k-1}'\leq p_{k-1}''$, with $p_1'\F\cdots\F p_{k-1}'$.  We also have
\[
x_{p_{k-1}'}x_{p_k''} \sim_2 x_{p_{k-1}'}x_{p_k'}, \WHERE p_k' = p_{k-1}'\th_{p_{k}''}.
\]
Since $p_{k-1}' \leq p_{k-1}'' \F p_k''$, it follows from Corollary \ref{cor:pqrs} that $p_{k-1}'\F p_k'$.  Putting everything together, we have
\[
w = x_{p_1}\cdots x_{p_k} \sim x_{p_1}\cdots x_{p_{k-2}}x_{p_{k-1}''}x_{p_k''} \sim x_{p_1'}\cdots x_{p_{k-2}'}x_{p_{k-1}'}x_{p_k''} \sim x_{p_1'}\cdots x_{p_{k-2}'}x_{p_{k-1}'}x_{p_k'},
\]
with all conditions met.  (Note that $p_i' \leq p_i''\leq p_i$ for $i=k-1,k$.)
\epf

For a path $\p=(p_1,\ldots,p_k)\in\P$, we define
\[
\wh\p = w_\p\psi = \wh p_1\cdots\wh p_k = (\th_{p_1}\cdots\th_{p_k},\de_{p_k}\cdots\de_{p_1}) \in \MP.
\]
It follows from Lemma \ref{lem:wp}, and surjectivity of $\Psi$, that
\[
\FP = \set{\ol w_\p}{\p\in\P} \AND \MP = \set{\wh\p}{\p\in\P}.
\]

\begin{lemma}\label{lem:whpq}
If $\p,\q\in\P$, then
\ben
\item \label{whpq1}  $\wh\p = \wh\q \implies \bd(\p)=\bd(\q)$ and $\br(\p)=\br(\q)$,
\item \label{whpq2}  $\ol w_\p = \ol w_\q \implies \bd(\p)=\bd(\q)$ and $\br(\p)=\br(\q)$.
\een
\end{lemma}

\pf
Since $\ol w_\p = \ol w_\q \implies w_\p\sim w_\q \implies\wh\p = w_\p\psi = w_\q\psi = \wh\q$, it suffices to prove \ref{whpq1}.  For this, write ${\p=(p_1,\ldots,p_k)}$ and $\q=(q_1,\ldots,q_l)$, so that
\[
\wh\p = (\th_{p_1}\cdots\th_{p_k},\de_{p_k}\cdots\de_{p_1})
\AND
\wh\q = (\th_{q_1}\cdots\th_{q_l},\de_{q_l}\cdots\de_{q_1}).
\]
It then follows from Lemma \ref{lem:pkql} that $\wh\p=\wh\q$ implies $p_1=q_1$ and $p_k=q_l$.
\epf

It follows that we have well-defined maps $D,R:\FP\to\FP$ and $D,R:\MP\to\MP$ given by
\begin{align}
\nonumber D(\ol w_\p) &= \ol x_{\bd(\p)} = \ol x_{p_1}, & D(\wh\p) &= \wh{\bd(\p)} = \wh p_1 = (\th_{p_1},\de_{p_1}), \\
\label{eq:DRwp} R(\ol w_\p) &= \ol x_{\br(\p)} = \ol x_{p_k},  & R(\wh\p) &= \wh{\br(\p)} = \wh p_k = (\th_{p_k},\de_{p_k}), && \ \ \  \text{for $\p=(p_1,\ldots,p_k)\in\P$.}
\end{align}
There should be no confusion in using $D$ and $R$ to denote these operations on both semigroups~$\FP$ and $\MP$, as they take on different kinds of arguments in the two semigroups.
Note also that
\begin{equation}\label{eq:DRPhi}
D(\ol w_\p)\Psi = (\ol x_{\bd(\p)})\Psi = \wh{\bd(\p)} = D(\wh \p) = D(\ol w_\p\Psi)
\ANDSIM
R(\ol w_\p)\Psi = R(\ol w_\p\Psi).
\end{equation}

\begin{lemma}\label{lem:wprest}
If $\p\in\P$, and if $s\leq\bd(\p)$ and $t\leq\br(\p)$, then
\[
x_sw_\p \sim w_{{}_s\corest\p} \AND w_\p x_t \sim w_{\p\rest_t}.
\]
\end{lemma}

\pf
For the first (the second is dual), write $\p = (p_1,\ldots,p_k)$, and fix $s\leq p_1$.  We then have
\[
{}_s\corest\p = (s_1,\ldots,s_k), \WHERE s_i = s\th_{p_1}\cdots\th_{p_i} \qquad\text{for each $i$,}
\]
and we must show that
\[
x_sx_{p_1}\cdots x_{p_k} \sim x_{s_1}\cdots x_{s_k}.
\]
If $k=1$, then $x_s x_{p_1} \sim_2 x_s x_{s\th_{p_1}} = x_sx_s \sim_1 x_s = x_{s_1}$ (using $s\leq p_1$ in the second step).  For $k\geq2$, we apply induction to calculate
\[
x_sx_{p_1}\cdots x_{p_{k-1}}x_{p_k} \sim x_{s_1}\cdots x_{s_{k-1}}x_{p_k} \sim_2 x_{s_1}\cdots x_{s_{k-1}}x_{s_{k-1}\th_{p_k}},
\]
and we are done, since $s_{k-1}\th_{p_k} = s_k$.
\epf

For paths $\p = (p_1,\ldots,p_k)$ and $\q = (q_1,\ldots,q_l)$ with $p_k\F q_1$, the concatenation of $\p$ with $\q$ is a path, which we denote by
\[
\p \op \q = (p_1,\ldots,p_k,q_1,\ldots,q_l).
\]

\begin{lemma}\label{lem:wpwq}
For paths $\p = (p_1,\ldots,p_k)$ and $\q = (q_1,\ldots,q_l)$, we have
\[
w_\p w_\q \sim w_{\p\rest_s \op {}_t\corest\q}, \WHERE \text{$s = q_1\de_{p_k}$ and $t = p_k\th_{q_1}$.}
\]
\end{lemma}

\pf
Using \ref{R1}, and Lemmas \ref{lem:xp'xq'} and \ref{lem:wprest}, we have
\[
w_\p w_\q \sim w_\p\cdot x_{p_k}x_{q_1} \cdot w_\q \sim w_\p\cdot x_sx_t \cdot w_\q \sim w_{\p\rest_s} \cdot w_{{}_t\corest\q}  = w_{\p\rest_s \op {}_t\corest\q}.  \qedhere
\]
\epf

\begin{prop}\label{prop:SP}
$\FP$ is a DRC-semigroup.
\end{prop}

\pf
We need to check the axioms \ref{DRC1}--\ref{DRC4}, and as usual we only have to prove one part of each.  To do so, fix paths $\p = (p_1,\ldots,p_k)$ and $\q = (q_1,\ldots,q_l)$.  

\pfitem{\ref{DRC1}}  We have $D(\ol w_\p) \ol w_\p = \ol x_{p_1} \cdot \ol x_{p_1}\cdots \ol x_{p_k} =_1 \ol x_{p_1}\cdots \ol x_{p_k} = \ol w_\p$.

\pfitem{\ref{DRC2}}  Let $s = q_1\de_{p_k}$ and $t = p_k\th_{q_1}$, and write $\p\rest_s = (s_1,\ldots,s_k)$ and ${}_t\corest\q = (t_1,\ldots,t_l)$.  Lemma~\ref{lem:wpwq} then gives
\[
\ol w_\p\ol w_\q = \ol w_{\p\rest_s \op {}_t\corest\q}
\AND
\ol w_\p D(\ol w_\q) = \ol w_\p \cdot \ol x_{q_1} = \ol w_\p \cdot \ol w_{(q_1)} = \ol w_{\p\rest_s \op {}_t\corest(q_1)}.
\]
Since $\p\rest_s \op {}_t\corest\q = (s_1,\ldots,s_k,t_1,\ldots,t_l)$ and $\p\rest_s \op {}_t\corest(q_1) = (s_1,\ldots,s_k,t)$, it then follows that
\[
D(\ol w_\p\ol w_\q) = \ol x_{s_1} = D(\ol w_\p D(\ol w_\q)).
\]
\firstpfitem{\ref{DRC3}}  Keeping the notation of the previous part, and keeping in mind $s_1\leq p_1$, we have
\[
D(\ol w_\p\ol w_\q)D(\ol w_\p) = \ol x_{s_1} \ol x_{p_1} =_2 \ol x_{s_1} \ol x_{s_1\th_{p_1}} = \ol x_{s_1}\ol x_{s_1} =_1 \ol x_{s_1} = D(\ol w_\p\ol w_\q).
\]
An analogous calculation gives $D(\ol w_\p)D(\ol w_\p\ol w_\q) = D(\ol w_\p\ol w_\q)$.

\pfitem{\ref{DRC4}}  We have $R(D(\ol w_\p)) = R(\ol x_{p_1}) = \ol x_{p_1} = D(\ol w_\p)$.
\epf

\begin{thm}\label{thm:FPMP}
If $P$ is a projection algebra, then $\FP$ and $\MP$ are projection-generated DRC-semigroups, with
\[
\bP(\FP) \cong \bP(\MP) \cong P.
\]
\end{thm}

\pf
We saw in Proposition \ref{prop:SP} that $\FP$ is a DRC-semigroup.  It follows that
\[
\bP(\FP) = \im(D) = \set{\ol x_p}{p\in P},
\]
so we have a surjective map
\[
\xi:P\to\bP(\FP):p\mt \ol x_p.
\]
Combining the definitions with Corollary \ref{cor:thp=thq} shows that $\xi$ is also injective, as for $p,q\in P$ we have
\begin{equation}\label{eq:PPbij}
p\xi = q\xi \implies \ol x_p = \ol x_q \implies \wh p = \ol x_p\Psi = \ol x_q\Psi = \wh q \implies (\th_p,\de_p) = (\th_q,\de_q) \implies p = q.
\end{equation}
So $\xi:P\to\bP(\FP)$ is a bijection, and we now wish to show that it is a projection algebra morphism, meaning that
\[
(p\th_q)\xi = (p\xi)\th'_{q\xi} \AND (p\de_q)\xi = (p\xi)\de'_{q\xi} \qquad\text{for all $p,q\in P$,}
\]
where here we write $\th'$ and $\de'$ for the projection algebra operations in $\bP(\FP)$.  For the first (the second is dual), we write $p'=q\de_p$ and $q'=p\th_q$, and use Lemma \ref{lem:xp'xq'} to calculate
\[
(p\xi)\th'_{q\xi} = (\ol x_p)\th'_{\ol x_q} = R(\ol x_p\ol x_q) = R(\ol x_{p'} \ol x_{q'}) = R(\ol w_{(p',q')}) = \ol x_{q'} = q'\xi = (p\th_q)\xi.
\]
By definition, $\FP = \pres{X_P}{R_P}$ is generated by the projections $\ol x_p$ ($p\in P$).  

This completes the proof of the assertions regarding $\FP$.  Those for $\MP$ follow quickly, in light of the fact that the surjective homomorphism $\Psi:\FP\to\MP$ preserves the $D$ and $R$ operations on the two semigroups (cf.~\eqref{eq:DRPhi}), and maps $\bP(\FP)$ bijectively (as $\ol x_p\Psi=\ol x_q\Psi\implies p=q$, as in~\eqref{eq:PPbij}).
\epf

\begin{rem}\label{rem:P=PFP}
By identifying each projection $p\in P$ with the $\sim$-class $\ol x_p\in\FP$, we can identify~$P$ itself with the projection algebra $\bP(\FP)$.  It will be convenient to do so in the next section.
\end{rem}

\subsection{Free projection-generated DRC-semigroups}\label{subsect:adjoint}

Recall that we have a functor $\bP:\DRC\to\PA$, which at the object level maps a DRC-semigroup~$S$ to its projection algebra
\[
\bP(S) = \im(D) = \im(R) = \set{p\in S}{p^2=p=D(p)=R(p)},
\]
with operations given by
\[
q\th_p = R(qp) \AND q\de_p = D(pq) \qquad\text{for $p,q\in \bP(S)$.}
\]
We can think of $\bP$ as a \emph{forgetful functor}, as when we construct $\bP(S)$ from $S$, we remember only a subset of the elements of $S$, and retain only partial information about their products.

It turns out that $\bP$ has a left adjoint, $\bF:\PA\to\DRC$, which involves the semigroups constructed in Section \ref{subsect:FP}.  At the object level, we define
\[
\bF(P) = \FP \qquad\text{for a projection algebra $P$.}
\]
To see how $\bF$ acts on a projection algebra morphism $\phi:P\to P'$, we begin by defining a semigroup morphism
\[
\varphi:X_P^+ \to \FPd = X_{P'}^+/R_{P'}^\sharp \BY x_p\varphi = \ol x_{p\phi} \qquad\text{for $p\in P$.}
\]
It is easy to see that $R_P \sub \ker(\varphi)$.  For example, if $p,q\in P$ then
\begin{align*}
(x_px_q)\varphi = \ol x_{p\phi} \ol x_{q\phi} &= \ol x_{p\phi} \ol x_{(p\phi)\th'_{q\phi}} &&\text{as $R_{P'}$ contains the relation $x_{p\phi}x_{q\phi} = x_{p\phi}x_{(p\phi)\th'_{q\phi}}$}\\
&= \ol x_{p\phi} \ol x_{(p\th_q)\phi} &&\text{as $\phi$ is a projection algebra morphism}\\
&= (x_px_{p\th_q})\varphi,
\end{align*}
which shows that $\varphi$ preserves \ref{R2}.  It follows that $\varphi$ induces a well-defined semigroup homomorphism
\begin{equation}\label{eq:Fphi}
\Phi :\FP = X_P^+/R_P^\sharp \to\FPd \GIVENBY \ol w\Phi = w\varphi \qquad\text{for $w\in X_P^+$.}
\end{equation}
We then define $\bF(\phi) = \Phi$.  It is essentially trivial to check that $\Phi$ preserves the $D$ and $R$ operations on~$\FP$ and~$\FPd$, in light of the definitions of these operations in \eqref{eq:DRwp}, meaning that $\bF(\phi)=\Phi$ is indeed a well-defined DRC-morphism $\bF(P)\to\bF(P')$.  

\begin{prop}\label{prop:free1}
$\bF$ is a functor $\PA\to\DRC$, and we have $\bP\bF = \id_{\PA}$.
\end{prop}

\pf
The first claim follows quickly from the above discussion.  For the second, we must show that:
\bit
\item $\bP(\bF(P)) = P$ for all projection algebras $P$, and 
\item $\bP(\bF(\phi)) = \phi$ for all projection algebra morphisms $\phi$.
\eit
This is routine, in light of the identification of $P$ with $\bP(\FP)$, as in Remark \ref{rem:P=PFP}; cf.~\eqref{eq:Fphi}.  
\epf

\begin{prop}\label{prop:free2}
For every projection algebra $P$, every DRC-semigroup $S$, and every projection algebra morphism $\phi:P\to\bP(S)$, there exists a unique DRC-morphism $\Phi:\bF(P)\to S$ such that $\phi = \bP(\Phi)$, i.e.~such that the following diagram of maps commutes (with vertical arrows being inclusions):
\[
\begin{tikzcd}[sep=small]
P \arrow{rr}{\ \ \phi} \arrow[swap,hookrightarrow]{dd} & ~ & \bP(S) \arrow[hookrightarrow]{dd} \\%
~&~&~\\
\bF(P) \arrow{rr}{\Phi\ \ }& ~ & S.
\end{tikzcd}
\]
\end{prop}

\pf
Fix a projection algebra morphism $\phi:P\to\bP(S)$.  For convenience we write
\[
\ul p = p\phi\in\bP(S) \qquad\text{for $p\in P$.}
\]
We write $D$ and $R$ for the unary operations on $\FP = \bF(P)$, and $\th$ and $\de$ for the projection algebra operations in $P=\bP(\FP)$, and use dashes to distinguish the corresponding operations in~$S$ and~$\bP(S)$.  

We begin by defining 
\[
\varphi:X_P^+\to S \BY x_p\varphi = \ul p \qquad\text{for $p\in P$.}
\]
To see that $R_P\sub\ker(\varphi)$, note for example that if $p,q\in P$ then 
\begin{align*}
(x_px_q)\varphi = \ul p \spc \ul q &= \ul p \spc \ul q \cdot R'(\ul p \spc \ul q) &&\text{by \ref{DRC1}}\\
&= \ul p \cdot R'(\ul p \spc \ul q) &&\text{by Lemma \ref{lem:Dpa}}\\
&= \ul p \cdot \ul p \th'_{\ul q} &&\text{by definition}\\
&= \ul p \cdot \ul{p\th_q} &&\text{as $\phi$ is a projection algebra morphism}\\
&= (x_px_{p\th_q})\varphi,
\end{align*}
showing that $\varphi$ preserves \ref{R2}.  As usual, it follows that $\varphi$ induces a well-defined semigroup morphism
\[
\Phi:\FP\to S \GIVENBY \ol w\Phi = w\varphi \qquad\text{for $w\in X_P^+$.}
\]
To see that $\Phi$ preserves $D$ ($R$ is dual), consider a typical element $\ol w_\p$ of $\FP$ (cf.~Lemma \ref{lem:wp}), where $\p = (p_1,\ldots,p_k) \in \P$.  We then have
\[
D(\ol w_\p)\Phi = \ol x_{p_1}\Phi = \ul p_1
\AND
D'(\ol w_\p\Phi) = D'(\ul p_1\cdots\ul p_k),
\]
so we must show that $D'(\ul p_1\cdots\ul p_k) = \ul p_1$.  But this follows from Lemma \ref{lem:Dpp}\ref{Dpp2}, together with the fact that $p_1\F\cdots\F p_k$ implies $\ul p_1\F\cdots\F\ul p_k$ (as $\phi$ is a projection algebra morphism).

So $\Phi$ is indeed a DRC-morphism $\FP\to S$.  Since $\Phi$ maps $p \equiv \ol x_p$ to $\ul p = p\phi$, it is clear that ${\bP(\Phi) = \Phi|_P = \phi}$.  This establishes the existence of $\Phi$.  Uniqueness follows from the fact that~${\FP = \la P\ra}$.
\epf

Propositions \ref{prop:free1} and \ref{prop:free2} verify the assumptions of \cite[Lemma 5.6]{EGPAR2024}, which immediately gives the following result.
For the definitions of the (standard) categorical terms in the statement see for example \cite{Awodey2010,MacLane1998}, and also \cite[Section 5]{EGPAR2024}.

\begin{thm}\label{thm:free}
The functor $\PA\to\DRC:P\mt\FP$ is a left adjoint to the forgetful functor $\DRC\to\PA:S\mt\bP(S)$, and $\PA$ is coreflective in $\DRC$.  \qed
\end{thm}

\begin{rem}
One might wonder if the assignment $P\mt\MP$ is also the object part of a functor~${\PA\to\DRC}$, but this turns out not to be the case, for exactly the same reason discussed in \cite[Remark 4.29]{EPA2024}.
\end{rem}

\begin{rem}
The study of free idempotent-generated semigroups over biordered sets is a topic of considerable interest in its own right \cite{GR2012,GR2012b,DGR2017, Dolinka2021,DDG2019,YDG2015,DG2014,BMM2009}, and the corresponding theory for regular $*$-semigroups has been recently initiated in \cite{EGPAR2024}.  We hope that the free projection-generated DRC-semigroups introduced here will likewise inspire further studies.
\end{rem}

\subsection{Fundamental projection-generated DRC-semigroups}\label{subsect:fund}

Following \cite{Jones2021}, a \emph{(DRC-)congruence} on a DRC-semigroup $S$ is an equivalence relation $\si$ on $S$ respecting all of the operations, in the sense that
\[
[a\mr\si b \text{ and } a'\mr\si b'] \implies aa' \mr\si bb' \AND a\mr\si b \implies [ D(a)\mr\si D(b) \text{ and } R(a) \mr\si R(b)].
\]
The congruence $\si$ is \emph{projection-separating} if $p\mr\si q\implies p=q$ for all projections $p,q\in\bP(S)$.  Among other things, \cite[Proposition 2.3]{Jones2021} shows that $S$ has a \emph{maximum projection-separating congruence}, denoted $\mu_S$; that is, $\mu_S$ is projection-separating, and any other projection-separating congruence is contained in $\mu_S$.  We say $S$ is \emph{(DRC-)fundamental} if $\mu_S$ is the trivial congruence $\io_S = \set{(a,a)}{a\in S}$.  In general, $S/\mu_S$ is fundamental, the maximum fundamental image of~$S$.  

\begin{rem}
Here we follow Jones' convention in dropping the `DRC-' prefix, and speak simply of `fundamental DRC-semigroups'.  But the reader should be aware that a `fundamental' (ordinary) semigroup is usually defined as a semigroup with no non-trivial congruence contained in Green's $\H$-relation.  For regular semigroups, this is equivalent to having no non-trivial idempotent-separating congruence \cite{Lallement1967}.  
\end{rem}

\begin{prop}\label{prop:muS}
If $S$ is a DRC-semigroup, then
\begin{align*}
\mu_S &= \set{(a,b)\in S\times S}{\Th_a=\Th_b\text{ and }\De_a=\De_b} \\
&= \set{(a,b)\in S\times S}{\vt_a=\vt_b\text{ and }\vd_a=\vd_b}
\end{align*}
\end{prop}

\pf
Writing $P = \bP(S)$, Proposition 2.3 of \cite{Jones2021} says that for $a,b\in S$ we have $(a,b)\in\mu_S$ if and only if
\[
D(ap)=D(bp) \ANd R(pa)=R(pb) \qquad\text{for all $p\in P^1$,}
\]
where $P^1$ is $P$ with an adjoined identity element.  Keeping \eqref{eq:Rpa} in mind, this is equivalent to all of the following holding:
\[
\Th_a=\Th_b \COMMA \De_a=\De_b \COMMA D(a)=D(b) \AND R(a)=R(b).
\]
Remembering also that $\bd=D$ and $\br=R$, we can therefore prove the first equality in the current proposition by showing that 
\begin{equation}\label{eq:DapDbp3}
[\Th_a=\Th_b \text{ and } \De_a=\De_b] \IMPLIES [\bd(a) = \bd(b) \text{ and } \br(a) = \br(b)].
\end{equation}
To do so, suppose $\Th_a=\Th_b$ and $\De_a=\De_b$.  Writing $q = \br(a)$, we have
\[
\bd(a) = \bd(a\rest_q) = q\vd_a = q\de_q\vd_a = q\De_a = q\De_b \leq \bd(b) \ANDSIM \bd(b) \leq \bd(a).
\]
The proof that $\br(a) = \br(b)$ is dual.

For the second equality in the proposition, we will show that
\begin{equation}\label{eq:ThaThb}
[\Th_a=\Th_b \text{ and } \De_a=\De_b] \IFF [\vt_a=\vt_b \text{ and } \vd_a=\vd_b] \qquad\text{for $a,b\in S$.}
\end{equation}
In light of the identities $\Th_a = \th_{\bd(a)}\vt_a$ and $\vt_a = \Th_a|_{\bd(a)^\da}$ (the second of which is Lemma~\ref{lem:Thpa}\ref{Thpa4}), and the analogous identities for $\De_a$ and~$\vd_a$, we can prove \eqref{eq:ThaThb} by showing that
\begin{align*}
[\Th_a=\Th_b \text{ and } \De_a=\De_b] &\IMPLIES [\bd(a) = \bd(b) \text{ and } \br(a) = \br(b)], \\
\text{and} \qquad [\vt_a=\vt_b \text{ and } \vd_a=\vd_b]\hspace{0.23cm} &\IMPLIES [\bd(a) = \bd(b) \text{ and } \br(a) = \br(b)].
\end{align*}
The first implication is \eqref{eq:DapDbp3}, which was proved above.
The second follows from the fact that $\dom(\vt_a) = \bd(a)^\da$ and $\dom(\vd_a) = \br(a)^\da$.  
\epf

Fundamental DRC-semigroups were classified in \cite{Jones2021} by means of a transformation representation.  We will not attempt to reprove this here, but instead we prove a result concerning \emph{projection-generated} fundamental DRC-semigroups (see Theorem~\ref{thm:fund} below), which was not given in \cite{Jones2021}.  For its proof, we will make use of the following result, which shows that the maps~$\Th_a$ and~$\De_a$ have very natural forms when $a$ is a product of projections.

\begin{lemma}\label{lem:Thppp}
If $S$ is a DRC-semigroup, and if $p_1,\ldots,p_k\in\bP(S)$, then
\[
\Th_{p_1\cdots p_k} = \th_{p_1}\cdots\th_{p_k} \AND \De_{p_1\cdots p_k} = \de_{p_k}\cdots\de_{p_1}.
\]
\end{lemma}

\pf
To prove the first statement (the second is dual), let $t\in\bP(S)$.  Then using \eqref{eq:Rpa} and Lemma~\ref{lem:Dpp}\ref{Dpp1}, we have $t\Th_{p_1\cdots p_k} = R(tp_1\cdots p_k) = t\th_{p_1}\cdots\th_{p_k}$.
\epf

We are now ready to prove our final result of this section:

\begin{thm}\label{thm:fund}
For any projection algebra $P$ there is a unique (up to isomorphism) projection-generated fundamental DRC-semigroup with projection algebra $P$, namely $\MP \cong \FP/\mu_{\FP}$.
\end{thm}

\pf
By Theorem \ref{thm:FPMP}, we can identify the projection algebras $\bP(\FP)$ and $\bP(\MP)$ with $P$ itself, and in this way the semigroups $\FP$ and $\MP$ are both generated by $P$.  To avoid confusion, we will denote the products in these semigroups by $\pr$ and $\star$, respectively, so that typical elements have the form
\[
p_1\pr\cdots\pr p_k \equiv \ol{x_{p_1}\cdots x_{p_k}} \in \FP \ANd p_1\star\cdots\star p_k \equiv \wh p_1\cdots\wh p_k = (\th_{p_1}\cdots\th_{p_k},\de_{p_k}\cdots\de_{p_1}) \in \MP,
\]
for some $p_1,\ldots,p_k\in P$.  In what follows, we use Proposition \ref{prop:muS} freely.

We have already noted that $\MP$ is a projection-generated DRC-semigroup with projection algebra $P$.  To see that it is fundamental, we need to show that $\mu_{\MP}$ is trivial.  To do so, suppose $(a,b)\in\mu_{\MP}$; we must show that $a=b$.  As above, we have
\[
a = p_1\star\cdots\star p_k = (\th_{p_1}\cdots\th_{p_k},\de_{p_k}\cdots\de_{p_1}) \AND b = q_1\star\cdots\star q_l = (\th_{q_1}\cdots\th_{q_l},\de_{q_l}\cdots\de_{q_1}),
\]
for some $p_1,\ldots,p_k,q_1,\ldots,q_l\in P$.  Since $(a,b)\in\mu_{\MP}$, we have $\Th_a=\Th_b$ and $\De_a=\De_b$.  It then follows from Lemma \ref{lem:Thppp} that
\[
\th_{p_1}\cdots\th_{p_k} = \Th_a = \Th_b = \th_{q_1}\cdots\th_{q_l} \ANDSIM \de_{p_k}\cdots\de_{p_1} = \de_{q_l}\cdots\de_{q_1}.
\]
But then $a =  (\th_{p_1}\cdots\th_{p_k},\de_{p_k}\cdots\de_{p_1}) = (\th_{q_1}\cdots\th_{q_l},\de_{q_l}\cdots\de_{q_1}) = b$.

Now let $S$ be an arbitrary projection-generated fundamental DRC-semigroup with projection algebra $\bP(S) = P$; we denote the product in $S$ simply by juxtaposition.  The proof will be complete if we can show that $S\cong\FP/\mu_{\FP}$.  Applying Proposition \ref{prop:free2} to the identity morphism $\phi:P\to P=\bP(S)$, we see that there is a DRC-morphism
\[
\Phi:\FP \to S:p\mt p.
\]
Since $S$ is projection-generated, $\Phi$ is surjective, and so $S\cong\FP/\ker(\Phi)$ by the fundamental homomorphism theorem.  It therefore remains to show that $\ker(\Phi) = \mu_{\FP}$.  Since $\Phi$ maps $P = \bP(\FP)$ identically, $\ker(\Phi)$ is projection-separating, and so $\ker(\Phi) \sub \mu_{\FP}$.  Conversely, let $(a,b)\in\mu_{\FP}$, so that $\Th_a=\Th_b$ and $\De_a=\De_b$.  Since $\FP$ is projection-generated, we have
\[
a = p_1\pr\cdots\pr p_k  \AND b = q_1\pr\cdots\pr q_l \qquad\text{for some $p_1,\ldots,p_k,q_1,\ldots,q_l\in P$.}
\]
But then applying Lemma \ref{lem:Thppp} in both $\FP$ and $S$, we have
\[
\Th_a = \Th_{p_1\pr\cdots\pr p_k} = \th_{p_1}\cdots\th_{p_k} = \Th_{p_1\cdots p_k} = \Th_{a\Phi} \ANDSIM \Th_b = \Th_{b\Phi}.
\]
It follows that $\Th_{a\Phi} = \Th_a = \Th_b = \Th_{b\Phi}$, and simlarly $\De_{a\Phi} = \De_{b\Phi}$, so that $(a\Phi,b\Phi) \in \mu_S$.  Since $S$ is fundamental, it follows that $a\Phi = b\Phi$, i.e.~that $(a,b)\in\ker(\Phi)$.
\epf

\section{\boldmath Regular $*$- and $*$-regular semigroups}\label{sect:*}

As mentioned in Section \ref{sect:intro}, the current paper extends to DRC-semigroups the groupoid-based approach to \emph{regular $*$-semigroups} from \cite{EPA2024}.  In the case of regular $*$-semigroups (whose definition will be recalled below), many of the current categorical constructions simplify, sometimes drammatically.  In this section we outline many of these simplifications, and indicate why the more elaborate approach of the current paper is necessary to encompass the full generality of DRC-semigroups.  In doing so, we will contrast the situation with Drazin's broader class of \emph{$*$-regular semigroups} \cite{Drazin1979}, which induce natural DRC-structures.  These are already complex enough to require the more complicated setup, and a ready source of (counter)examples will be provided by the real matrix monoids $\MnR$, which are known to be $*$-regular \cite{Penrose1955}.  (Due to quirks of terminology, $*$-regular and regular $*$-semigroups are distinct classes, though the former contains the latter.)  We will also observe in Section \ref{sect:other} that many of the simplifications hold in the \emph{DRC-restriction semigroups} considered by Die and Wang in \cite{DW2024}.

Before we begin, we recall some basic semigroup theoretical background; for more information see for example \cite{Howie1995,CPbook}.  Let $S$ be a semigroup, and $S^1$ its monoid completion.  So $S^1=S$ if $S$ is a monoid, or else $S^1=S\cup\{1\}$, where $1$ is a symbol not belonging to $S$, acting as an adjoined identity.  Green's~$\L$,~$\R$ and~$\J$ equivalences are defined, for $a,b\in S$, by
\[
a\L b \iff S^1a=S^1b \COMMA a\R b \iff aS^1=bS^1 \AND a\J b \iff S^1aS^1 = S^1bS^1.
\]
From these are defined the further equivalences ${\H}={\L}\cap{\R}$ and ${\D}={\L}\vee{\R}$, where the latter is the join in the lattice of all equivalences on $S$.  We have ${\D}={\L}\circ{\R} = {\R}\circ{\L}$, and if $S$ is finite then ${\D}={\J}$.

An element $a$ of a semigroup $S$ is \emph{regular} if $a=axa$ for some $x\in S$.  This is equivalent to having $a=aya$ and $y=yay$ for some $y\in S$ (if $a=axa$, then take $y=xax$), and such an element $y$ is called a \emph{(semigroup) inverse} of $a$.  Note then that $a\R ay$ and $a\L ya$, with~$ay$ and~$ya$ idempotents.  In fact, an element of a semigroup is regular if and only if it is $\D$-related to an idempotent.

\subsection{Matrix monoids}\label{subsect:MnR}

Fix an integer $n\geq1$, and let $\MnR$ be the monoid of $n\times n$ real matrices, under ordinary matrix multiplication.  (We consider real matrices merely for convenience, though what we say can be adapted to the complex field.)  Denoting the row and column spaces of $a\in\MnR$ by $\Row(a)$ and $\Col(a)$, and writing $\rank(a) = \dim(\Row(a)) = \dim(\Col(a))$ for the rank of $a$, Green's relations on $\MnR$ are given by
\begin{align}
\nonumber a\L b &\iff \Row(a)=\Row(b) ,\\
\label{eq:Green_MnR} a\R b &\iff \Col(a)=\Col(b) \AND a\J b\iff a\D b \iff \rank(a)=\rank(b).
\end{align}
(See for example \cite[Lemma 2.1]{Okninski1998}.)  Since idempotents of arbitrary rank $0\leq r\leq n$ clearly exist,~$\MnR$ is regular.

In fact, any matrix $a\in\MnR$ can be assigned a \emph{special} inverse \cite{Penrose1955}.  Using ${}^\tr$ to denote transpose, there exists a unique matrix $x\in\MnR$ satisfying 
\begin{equation}\label{eq:AXA}
a=axa \COMMA x=xax \COMMA (ax)^\tr=ax \AND (xa)^\tr=xa.
\end{equation}
This $x$ is known as the \emph{Moore--Penrose inverse}, and is denoted $x = a^\+$.  There is no simple formula for $a^\+$; its existence is established in \cite{Penrose1955} by using the linear dependence of the sequence $(a^\tr a)^k$, $k=1,2,3,\ldots$.  Nevertheless, many computational packages exist for working with the Moore--Penrose inverse, e.g.~Matlab \cite{Matlab}.  For example, with $a = \tmat{1&1&1\\1&0&0\\0&0&0} \in M_3(\RR)$, Matlab tells us that $a^\+ = \tmat{0&1&0\\\sfrac12&-\sfrac12&0\\\sfrac12&-\sfrac12&0}$.  One can readily check the identities \eqref{eq:AXA}, and we note that $aa^\+ = \tmat{1&0&0\\0&1&0\\0&0&0}$ and $a^\+a = \tmat{1&0&0\\0&\sfrac12&\sfrac12\\0&\sfrac12&\sfrac12}$ are symmetric idempotents, with $aa^\+ \R a \L a^\+a$; cf.~\eqref{eq:Green_MnR}.

\subsection[Regular $*$- and $*$-regular semigroups as DRC-semigroups]{\boldmath Regular $*$- and $*$-regular semigroups as DRC-semigroups}

The defining properties of the Moore--Penrose inverse in $\MnR$ led to the introduction in \cite{Drazin1979} of the class of \emph{$*$-regular semigroups}, to which we now turn.

First, a \emph{$*$-semigroup} is an algebra $(S,\cdot,{}^*)$, where $(S,\cdot)$ is a semigroup, and where ${}^*$ is an involution, i.e.~a unary operation $S\to S:a\mt a^*$ satisfying the laws
\[
a^{**} = a \AND (ab)^* = b^*a^* \qquad\text{for all $a,b\in S$.}
\]
For example, $\MnR$ becomes a $*$-semigroup with involution given by transpose: $a^*=a^\tr$.
The category of $*$-semigroups with $*$-morphisms (i.e.~the semigroup morphisms preserving ${}^*$) contains two important (full) subcategories, whose definitions we will shortly recall:
\bit
\item $\RSS$, the category of \emph{regular $*$-semigroups} \cite{NS1978}, and
\item $\SRS$, the category of \emph{$*$-regular semigroups} \cite{Drazin1979}.
\eit

Beginning with the latter, a \emph{regular $*$-semigroup} is a $*$-semigroup $(S,\cdot,{}^*)$ additionally satisfying
\[
a = aa^*a \qquad\text{for all $a\in S$.}
\]
It follows from the axioms that also $a^*aa^* = a^*$ in a regular $*$-semigroup, meaning that $a^*$ is a (semigroup) inverse of $a$.  

A \emph{$*$-regular semigroup} is a $*$-semigroup $(S,\cdot,{}^*)$ for which every element $a\in S$ has an inverse that is $\H$-related to $a^*$.  It follows from basic semigroup-theoretic facts that this inverse is unique; it is denoted by $a^\+$, and called the \emph{Moore--Penrose inverse} of $a$ \cite{Penrose1955}.  The terminology stems from the fact (see \cite[Corollary 1.2]{NP1985}) that $a^\+$ can also be characterised as the unique element $x$ satisfying
\[
a = axa \COMMA x = xax \COMMA (ax)^* = ax \AND (xa)^* = xa.
\]
In this way, a $*$-regular semigroup can be equivalently defined as an algebra $(S,\cdot,{}^*,{}^\+)$, where~$(S,\cdot)$ is a semigroup, and ${}^*$ and ${}^\+$ are unary operations satisfying:
\begin{equation}\label{eq:*reg}
a^{**} = a \COMMa
(ab)^* = b^*a^* \COMMa
aa^\+a = a \COMMa a^\+aa^\+ = a^\+ \COMMa
(aa^\+)^* = aa^\+ \ANd
(a^\+a)^* = a^\+a.
\end{equation}
Various other laws can be deduced as consequences (see for example \cite[Theorem 1.3]{NP1985}), including:
\begin{equation}\label{eq:*reg2}
a^{\+\+} = a \AND a^{*\+} = a^{\+*}.
\end{equation}
Regular $*$-semigroups are the $*$-regular semigroups in which~${a^\+=a^*}$ for all $a$.  The matrix monoid~$\MnR$ is $*$-regular (cf.~\eqref{eq:AXA}), but it is not a regular $*$-semigroup, as $a^\tr$ need not be an inverse of $a$, even when $n=1$.

It is well known (and easy to see) that any $*$-morphism $\phi:S\to T$ between $*$-regular semigroups preserves Moore--Penrose inverses, meaning that $(a^\+)\phi = (a\phi)^\+$ for all $a\in S$.  
Less obvious is the following:

\begin{thm}[Jones {\cite[Proposition 4.2]{Jones2021}}]\label{thm:Jones}
Any $*$-regular semigroup $S \equiv (S,\cdot,{}^*)$ gives rise to a DRC-semigroup $\bI(S) = (S,\cdot,D,R)$, under the unary operations
\[
D(a) = aa^\+ \AND R(a) = a^\+a \qquad\text{for $a\in S$.}  \epfreseq
\]
\end{thm}

\begin{rem}
The proof of Theorem \ref{thm:Jones} is completely routine in the special case that~$S$ is a regular $*$-semigroup.  We also note that the elements $D(a)=aa^\+$ and $R(a)=a^\+a$ appearing in the theorem are the unique projections of $S$ that are $\R$- and $\L$-related to $a$, respectively.  This was observed, for example, in \cite{CD2002}.
\end{rem}

We can think of the assignment $S\mt\bI(S)$ as an object map from $\SRS$ to $\DRC$, and this can be easily extended to a functor $\bI:\SRS\to\DRC$.
Given a $*$-morphism $\phi:S\to T$ in $\SRS$, we simply define $\bI(\phi)=\phi$.  Since $\phi$ preserves ${}^\+$, it also preserves~$D$ and~$R$.

The functor $\bI:\SRS\to\DRC$ is not full (surjective on morphisms), however.  Indeed, consider an abelian group $G$ with identity $1$, in which at least one element is not its own inverse (e.g.~the cyclic group of order $3$).  Then $S_1 = (G,\cdot,{}^{-1})$ and $S_2 = (G,\cdot,{}^\star)$ are both $*$-regular semigroups, where $a^\star= a$ for all $a\in G$.  Moreover, we have $\bI(S_1) = \bI(S_2)$, since $D(a) = R(a) = 1$ in both semigroups for all $a\in G$.  It follows that the identity map $\iota:G\to G$ is a DRC-morphism $\bI(S_1)\to \bI(S_2)$, but it is not a $*$-morphism, as $(a^{-1})\iota \not= (a\iota)^\star$ for any non-self-inverse $a\in G$.

On the other hand, the restriction of $\bI$ to the subcategory $\RSS$ of regular $*$-semigroups \emph{is} full, as follows from the next result (recall that ${}^\+$ coincides with ${}^*$ in a regular $*$-semigroup).

\begin{prop}
If $\phi:\bI(S)\to\bI(T)$ is a DRC-morphism, for $*$-regular semigroups $S$ and $T$, then $(a^\+)\phi = (a\phi)^\+$ for all $a\in S$.
\end{prop}

\pf
Since $\phi$ is a DRC-morphism, we have
\[
a\phi\cdot a^\+\phi = (aa^\+)\phi = D(a)\phi = D(a\phi) = a\phi\cdot (a\phi)^\+ \ANDSIM a^\+\phi\cdot a\phi = (a\phi)^\+\cdot a\phi.
\]
Combining these gives
\[
a^\+\phi = (a^\+\cdot a\cdot a^\+)\phi = a^\+\phi\cdot a\phi\cdot a^\+ \phi = a^\+\phi\cdot a\phi\cdot (a \phi)^\+ = (a\phi)^\+ \cdot a\phi\cdot (a \phi)^\+ = (a\phi)^\+.  \qedhere
\]
\epf

It follows that $\RSS$ is (isomorphic to) a full subcategory of $\DRC$, and is hence isomorphic to its image under the isomorphism $\bC:\DRC\to\CPC$ from Theorem \ref{thm:Cfunctor}.  One could then go on to show that this image is isomorphic to the category $\CPG$ of \emph{chained projection groupoids}, as defined in~\cite{EPA2024}, thereby proving the main result of that paper.  We will not give the full details here, as a similar deduction was considered in \cite{EPA2024} itself, where the isomorphism $\RSS\cong\CPG$ was used to reprove the Ehresmann--Nambooripad--Schein Theorem on inverse semigroups and inductive groupoids.

The situation for the more general $*$-regular semigroups, however, is more complicated.  Since the functor $\bI:\SRS\to\DRC$ is not an isomorphism onto a \emph{full} subcategory, we are left with what seems to be a very interesting and important open problem:

\begin{prob}\label{prob:*reg}
Obtain an `ESN-type' theorem for the category of $*$-regular semigroups.
\end{prob}

\subsection[Projection algebras of regular $*$- and $*$-regular semigroups]{\boldmath Projection algebras of regular $*$- and $*$-regular semigroups}

For the rest of this section we fix a $*$-regular semigroup $S \equiv (S,\cdot,{}^*,{}^\+)$, and denote the corresponding DRC-semigroup from Theorem~\ref{thm:Jones} by $S' = \bI(S) = (S,\cdot,D,R)$.  Using \eqref{eq:*reg} and \eqref{eq:*reg2} we see that
\begin{equation}\label{eq:Ra*}
R(a^*) = (a^*)^\+a^* = (a^\+)^*a^* = (aa^\+)^* = aa^\+ = D(a) \ANDSIM D(a^*) = R(a),
\end{equation}
meaning that $S'$ is an \emph{involutory DRC-semigroup} in the terminology of Jones \cite[Section 4]{Jones2021}.
Among other things, this leads to a simplification in the structure of the projection algebra $P = \bP(S')$.  First, we note that
\[
P = \set{p\in S}{p^2=p=p^*} = \set{p\in S}{p^2=p=p^\+}.
\]
Indeed, the first equality was shown in \cite[Proposition 4.2]{Jones2021}.  For the second, suppose first that $p^2=p=p^*$.  It follows that $p$ is an inverse of $p$ that is $\H$-related (indeed, \emph{equal}) to $p^*$, so that~${p=p^\+}$.  Conversely, if $p^2=p=p^\+$, then $p = pp = pp^\+ = (pp^\+)^* = p^*$.  Using \eqref{eq:Ra*}, it then follows that for any $p,q\in P$ we have
\[
q\th_p = R(qp) = 
R(q^*p^*) = 
R((pq)^*) = D(pq) = q\de_p.
\]
In other words, we have $\th_p = \de_p$ for all $p$, meaning that the two families of operations of~$P$ reduces to one, and $P$ is \emph{symmetric} in the terminology of Jones \cite[Section 11]{Jones2021}.

In fact, not only do the $\th_p$ and $\de_p$ operations coincide in a $*$-regular semigroup, but they can be defined equationally in terms of the ${}^\+$ operation.  To describe this, and for later use, we need the following basic result:

\begin{lemma}\label{lem:pa+}
If $S$ is a $*$-regular semigroup, then for any $a\in S$ and $p\in\bP(S)$ we have
\[
(pa)^\+ = (pa)^\+p \AND (ap)^\+ = p(ap)^\+.
\]
\end{lemma}

\pf
Since~$(pa)^\+ \H (pa)^* = a^*p$, it follows that $p$ is a right identity for $(pa)^\+$, which gives the first identity.  The second is dual.
\epf

\begin{prop}\label{prop:pthq}
If $S$ is a $*$-regular semigroup, then for any $p,q\in\bP(S)$ we have
\[
q\th_p = q\de_p = p(pq)^\+ = (qp)^\+p.
\]
\end{prop}

\pf
Using Lemma \ref{lem:pa+} we see that $q\th_p = R(qp) = (qp)^\+qp = (qp)^\+p$.  A symmetrical calculation gives ${q\de_p=p(pq)^\+}$, and we have already seen that $q\th_p=q\de_p$.
\epf

\begin{rem}\label{rem:pqp}
In the case that $S$ is a regular $*$-semigroup, we obtain the simpler expression
\begin{equation}\label{eq:pqp}
q\th_p = p(pq)^\+ = p(pq)^* = pqp.
\end{equation}
This does not hold in $*$-regular semigroups in general, however.  For example, with the projections $p = \tmat{\sfrac12&\sfrac12\\\sfrac12&\sfrac12}$ and $q = \tmatx{1&0\\0&0}$, both from $M_2(\RR)$, the matrix $pqp = \tmat{\sfrac14&\sfrac14\\\sfrac14&\sfrac14}$ is not an idempotent, and hence not a projection.
\end{rem}

The projection algebras of regular $*$-semigroups were axiomatised by Imaoka \cite{Imaoka1983}, as the algebras $P\equiv(P,\th)$, where $\th = \set{\th_p}{p\in P}$ is a set of unary operations satisfying
\begin{enumerate}[label=\textup{\textsf{(P\arabic*)$'$}},leftmargin=9mm]\bmc3
\item \label{P1'} $p\theta_p = p$,
\item \label{P2'} $\theta_p\theta_p=\theta_p$,
\item \label{P3'} $p\theta_q\theta_p=q\theta_p$,
\item \label{P4'} $\theta_p\theta_q\theta_p=\theta_{q\theta_p}$,
\item \label{P5'} $\theta_p\theta_q\theta_p\theta_q=\theta_p\theta_q$.
\item[] ~
\emc
\end{enumerate}
It is easy to see that these axioms hold in a regular $*$-semigroup (using \eqref{eq:pqp}), and that they imply axioms \ref{P1}--\ref{P5}, keeping $\th_p=\de_p$ in mind.  Note that \ref{P1'} and \ref{P2'} are exactly \ref{P1} and \ref{P6}, respectively, while \ref{P3'} is the $\th=\de$ version of \ref{P9}.  In particular, axioms \ref{P1'}--\ref{P3'} hold in the projection algebra of a $*$-regular semigroup.  Axioms \ref{P4'} and \ref{P5'} do not, however.  For example, in $M_3(\RR)$ we have
\[
t\theta_p\theta_q\theta_p \not= t\theta_{q\theta_p}
\AND
t\theta_p\theta_q\theta_p\theta_q \not= t\theta_p\theta_q
\]
for
\[
t = \tmat{\sfrac12&\sfrac12&0\\\sfrac12&\sfrac12&0\\0&0&1} \COMMA 
p = \tmat{\sfrac12&0&\sfrac12\\0&1&0\\\sfrac12&0&\sfrac12} \AND
q = \tmat{1&0&0\\0&\sfrac12&\sfrac12\\0&\sfrac12&\sfrac12}.
\]
(This example was found, and can be checked, with Matlab \cite{Matlab}, as with all the examples to come.)

As far as we are aware, there is no known axiomatisation for the projection algebras of $*$-regular semigroups.  Obtaining such an axiomatisation is a necessary first step in tackling Problem \ref{prob:*reg}.  We also believe it is of interest in its own right, as among other things it could also be useful in developing a theory of \emph{fundamental} $*$-regular semigroups; cf.~\cite{Imaoka1983,Jones2021,Jones2012,Yamada1981}.

\begin{prob}\label{prob:P}
Give an (abstract) axiomatisation for the class of projection algebras of $*$-regular semigroups.
\end{prob}

\subsection[Chained projection categories of regular $*$- and $*$-regular semigroups]{\boldmath Chained projection categories of regular $*$- and $*$-regular semigroups}\label{subsect:CPC*}

We continue to fix the  $*$-regular semigroup $S \equiv (S,\cdot,{}^*,{}^\+)$, and its associated DRC-semigroup $S' = \bI(S) = (S,\cdot,D,R)$ from Theorem \ref{thm:Jones}.  Also let $(P,\C,\ve) = \bC(S')$ be the chained projection category associated to~$S'$.  In particular, $P = \bP(S')$ is the projection algebra of  $S'$.  Note that for any $a\in\C\ (=S)$ we have
\[
\bd(a) = D(a) = aa^\+ \ANDSIM \br(a) = a^\+a.
\]

\begin{prop}
If $S$ is a $*$-regular semigroup, then $\C=\C(S')$ is a groupoid, in which $a^{-1} = a^\+$ for all $a\in\C$.
\end{prop}

\pf
Fix $a\in S$, and first note that \eqref{eq:*reg2} gives
\[
\bd(a^\+) = a^\+a^{\+\+} = a^\+a = \br(a) \ANDSIM \br(a^\+) = \bd(a).
\]
It follows that the compositions $a\circ a^\+$ and $a^\+\circ a$ exist in $\C$.  Specifically, we have
\[
\bd(a) = aa^\+ = a\circ a^\+ \ANDSIM \br(a) = a^\+\circ a,
\]
so that $a^\+$ is indeed the inverse of $a$ in $\C$.
\epf

Using \eqref{eq:*reg} and \eqref{eq:*reg2}, it is also worth noting that for $a\in S$ we have
\[
\bd(a^*) = a^*a^{*\+} = a^*a^{\+*} = (a^\+a)^* = a^\+a = \br(a) \ANDSIM \br(a^*) = \bd(a).
\]
In what follows, we will typically use the laws $\bd(a^\+)=\bd(a^*)=\br(a)$ and $\br(a^\+)=\br(a^*)=\bd(a)$ without explicit reference.

A further simplification occurs in the special case that $S$ is a regular $*$-semigroup (where we recall that ${}^*$ and ${}^\+$ coincide).

\begin{prop}
If $S$ is a regular $*$-semigroup, then the orders~$\leql$ and~$\leqr$ on $\C$ coincide.
\end{prop}

\pf
By symmetry it suffices to show that $a\leql b\implies a\leqr b$ for $a,b\in\C$, so suppose $a\leql b$.  This means that $a = {}_p\corest b = pb$, where $p = \bd(a) \leq \bd(b) = bb^*$.  It follows from the latter that $p = bb^*p$.  Write $q = \br(a) = a^*a$, and note that the biordered category axioms give $q\leq\br(b)$.  It then follows that
\[
b\rest_q = bq = ba^*a = b(pb)^*a = bb^*pa = pa = a,
\]
so that $a\leqr b$.  
\epf

\begin{rem}\label{rem:abpqe}
The orders $\leql$ and $\leqr$ are generally distinct for $*$-regular semigroups.  For example, consider the following matrices from $M_2(\RR)$:
\[
a = \tmatx{1&1\\0&0}  \COMMA b = \tmatx{1&1\\0&1} \COMMA p = \tmatx{1&0\\0&0} \COMMA q = \tmat{\sfrac12&\sfrac12\\\sfrac12&\sfrac12} \AND e = \tmatx{1&0\\0&1}.
\]
Then $p$, $q$ and $e$ are projections with $p,q\leq e$, and we have $D(a) = p$, $R(a) = q$ and ${D(b)=R(b)=e}$.  We also have ${}_{D(a)}\corest b = pb = a$, so that $a\leql b$.  But $a\not\leqr b$, since
$b\rest_{R(a)} = bq = \tmat{1&1\\\sfrac12&\sfrac12} \not= a$.
\end{rem}

Another special property enjoyed by the category $\C$ has to do with the $\vt/\Th$ and $\vd/\De$ maps from \eqref{eq:vtvd} and \eqref{eq:ThDe}.  We summarise this in Proposition \ref{prop:*3} below, but we first note that for $a\in\C$ we have
\begin{align}
\nonumber p\vt_a &= (pa)^\+a \text{ \ for $p\leq\bd(a)$,} & p\Th_a = (pa)^\+a \text{ \ for $p\in P$,}\\
\label{eq:Rpa*} p\vd_a &= a(ap)^\+ \text{ \ for $p\leq\br(a)$,} & p\De_a = a(ap)^\+ \text{ \ for $p\in P$.}
\end{align}
Indeed, if $p\leq\bd(a)$, then using \eqref{eq:Rpa} and Lemma \ref{lem:pa+} we have
\[
p\vt_a = R(pa) = (pa)^\+pa = (pa)^\+a.
\]
The other formulae are checked analogously.  If $S$ is a regular $*$-semigroup, then these reduce further; for example, $p\Th_a = a^*pa$ for $a\in\C$ and $p\in P$.  Such simplifications do not hold in general for $*$-regular semigroups; for example, for the projections $p,q\in M_2(\RR)$ from Remark~\ref{rem:pqp}, we have $q\Th_p = q\th_p \not= pqp = p^*qp$.

\newpage

\begin{prop}\label{prop:*3}
If $S$ is a $*$-regular semigroup, then for any $a\in\C$:
\ben
\item \label{*31} $\vt_a$ and $\vd_a$ are bijections, and $\vt_{a^\+} = \vt_a^{-1}$ and $\vd_{a^\+} = \vd_a^{-1}$,
\item \label{*32} $\vt_{a^*} = \vd_a$ and $\vd_{a^*}=\vt_a$,
\item \label{*33} $\Th_{a^*} = \De_a$ and $\De_{a^*} = \Th_a$.
\een
\end{prop}

\pf
For each item, it suffices by symmetry to prove the first assertion.

\pfitem{\ref{*31}}  It suffices by symmetry to show that $\vt_a\vt_{a^\+} = \id_{\bd(a)^\da}$, i.e.~that $p\vt_a\vt_{a^\+} = p$ for all $p\leq\bd(a)$.  To do so, fix some such $p$, and write $q = p\vt_a = \br({}_p\corest a)$.  We then have
\[
p = {}_p\corest \bd(a) = {}_p\corest(a\circ a^\+) = {}_p\corest a \circ {}_q\corest a^\+ ,
\]
from which it follows that $p = \br(p) = \br({}_p\corest a \circ {}_q\corest a^\+) = \br({}_q\corest a^\+) = q\vt_{a^\+} = p\vt_a\vt_{a^\+}$.

\pfitem{\ref{*32}}  
Let $p\leq\br(a)$.  Keeping in mind that $p\vd_a$ is a projection, we use \eqref{eq:Rpa*} and \eqref{eq:*reg2} to calculate
\[
p\vd_a = (p\vd_a)^* = (a(ap)^\+)^* = (ap)^{\+*}a^* = (ap)^{*\+}a^* = (pa^*)^\+a^* = p\vt_{a^*}.  
\]
\firstpfitem{\ref{*33}}  Using previously-established facts, we have $\Th_{a^*} = \th_{\bd(a^*)}\vt_{a^*} = \de_{\br(a)}\vd_a = \De_a$.
\epf

\begin{rem}
When $S$ is a regular $*$-semigroup, Proposition \ref{prop:*3}\ref{*31} says that $\vt_{a^*}$ is the inverse of $\vt_a$.  This was proved in \cite[Lemma 2.12]{EPA2024} using the identity ${}_p\corest a = a\rest_q$, where $p\leq\bd(a)$ and $q=\br({}_p\corest a)$.  This does not hold in a general $*$-regular semigroup.  For example, with $b,p,q\in M_2(\RR)$ as in Remark \ref{rem:abpqe}, we have $p\leq\bd(b)$ and $q = \br({}_p\corest b)$, but ${}_p\corest b \not= b\rest_q$.
\end{rem}

\begin{rem}\label{rem:G1}
Another point of difference concerning the $\vt_a$ and $\vd_a$ maps is that these are always projection algebra morphisms when $S$ is a regular $*$-semigroup, as shown in the proof \cite[Proposition 6.11]{EPA2024}.  This property was called \textsf{(G1d)} in \cite{EPA2024}, and was shown to be equivalent to certain other natural conditions in the context of the \emph{projection groupoids} of \cite{EPA2024}, including our~\ref{C1}.  However, the maps $\vt_a$ and $\vd_a$ need not be projection algebra morphims when~$S$ is $*$-regular.\footnote[4]{This is also the case in general for Ehresmann semigroups; see for example \cite[Section 4.2]{Gould2012}.}  For example, we have $(p\th_q)\vt_a \not= (p\vt_a)\th_{q\vt_a}$ for $a = \tmatx{1&1\\0&1}$, $p = \tmatx{1&0\\0&0}$ and $q = \tmatx{0&0\\0&1}$, all from $M_2(\RR)$.  It is worth noting that one of the other equivalent conditions in \textsf{(G1)} is:
\begin{equation}\label{eq:G1b}
\th_{p\Th_a} = \Th_{a^*}\th_p\Th_a \ (=\De_a\th_p\Th_a) \qquad\text{for all $p\in P$ and $a\in\C$.}
\end{equation}
This also need not hold for $*$-regular semigroups.  For example, for $p=\tmatx{1&0\\0&1}$, $t=\tmatx{1&0\\0&0}$ and $a=\tmatx{1&1\\0&1}$ from $M_2(\RR)$ we have $t\th_{p\Th_a} =\tmatx{1&0\\0&0}$ and $t\De_a\th_p\Th_a = \tmat{\sfrac12&\sfrac12\\\sfrac12&\sfrac12}$.
\end{rem}

The final issue we wish to address concerns axiom \ref{C2}, and is somewhat more technical.  Recall that this is quantified over all morphisms $b\in\C$ and all projections $p,s\in P$, and asserts the equality of the morphisms 
\begin{equation}\label{eq:lamrhoagain}
\lam(p,b,s) = \ve[e,e_1] \circ ({}_{p\th_q}\corest b)\rest_{f_1} \circ \ve[f_1,f]
\AND
\rho(p,b,s) = \ve[e,e_2] \circ {}_{e_2}\corest(b\rest_{s\de_r}) \circ \ve[f_2,f]
\end{equation}
in $\C$.  These morphisms were defined in terms of the projections $q=\bd(b)$, $r=\br(b)$, and $e,e_i,f,f_i$ given by
\begin{align}
\nonumber e(p,b,s) &= s\De_b\de_p , & e_1(p,b,s) &= s\De_b\de_{p\th_q} , & f_1(p,b,s) &= s\de_{p\Th_b} , & f(p,b,s) &= p\Th_b\th_s , \\
\label{eq:etc} & & e_2(p,b,s) &= p\th_{s\De_b} , & f_2(p,b,s) &= p\Th_b\th_{s\de_r}.
\end{align}
All of the above projections and morphisms are shown in Figure \ref{fig:lamrho}.  The idea behind axiom \ref{C2} is that if we have three morphisms $a,b,c\in\C$, with 
\[
p = \br(a) \COMMA q=\bd(b) \COMMA r=\br(b) \AND s=\bd(c),
\]
then
\[
(a\pr b)\pr c = a\rest_e \circ \lam(p,b,s) \circ {}_f\corest c
\AND
a\pr (b\pr c) = a\rest_e \circ \rho(p,b,s) \circ {}_f\corest c.
\]
(See the proof of Lemma \ref{lem:abc}.)

In the case of regular $*$-semigroups, the corresponding axiom \textsf{(G2)} from \cite{EPA2024} was quantified over all morphisms $b$, and all \emph{$b$-linked pairs} $(e,f)$.  Such a pair satisfies
\begin{equation}\label{eq:lp}
e\Th_b\th_f = f \AND f\De_b\de_e = e.
\end{equation}
(In the notation of \cite{EPA2024} the latter equation was $f\Th_{b^*}\th_e = e$, which is equivalent to the above since $\de_e=\th_e$ and $\De_b=\Th_{b^*}$.)  The idea here is that if $p,s\in P$ are arbitrary, then $e = e(p,b,s)$ and $f = f(p,b,s)$ from \eqref{eq:etc} are $b$-linked, and we have 
\begin{equation}\label{eq:pbsebf}
e(e,b,f) = e \COMMA e_i(e,b,f) = e_i(p,b,s) \COMMA f_i(e,b,f) = f_i(p,b,s) \AND f(e,b,f) = f,
\end{equation}
which in turn leads to $\lam(p,b,s) = \lam(e,b,f)$ and $\rho(p,b,s) = \rho(e,b,f)$.  For example, using \eqref{eq:G1b}, \eqref{eq:etc}, \eqref{eq:lp} and $\th_t=\de_t$, we have
\begin{align*}
f_1(e,b,f) = f\de_{e\Th_b} = f\th_{e\Th_b} = f\De_b\th_e\Th_b &= e\Th_b \\
&= (s\De_b\de_p)\Th_b = s\De_b\th_p\Th_b = s\th_{p\Th_b} = s\de_{p\Th_b} = f_1(p,b,s).
\end{align*}
In the general case of a $*$-regular semigroup, the equalities involving $e_i$ and $f_i$ in \eqref{eq:pbsebf} need not hold.  For example, consider the following matrices from $M_3(\RR)$:
\begin{equation}\label{eq:bps}
b = \tmat{1&1&0\\1&1&0\\0&0&1} \COMMA p = \tmat{1&0&0&\\0&\sfrac12&\sfrac12\\0&\sfrac12&\sfrac12} \AND s = \tmat{0&0&0\\0&0&0\\0&0&1}.
\end{equation}
Then $p$ and $s$ are projections, and we have
\begin{equation}\label{eq:efagain}
e = e(p,b,s) = \tmat{0&0&0&\\0&\sfrac12&\sfrac12\\0&\sfrac12&\sfrac12} \AND f = f(p,b,s) = s = \tmat{0&0&0\\0&0&0\\0&0&1}.
\end{equation}
Moreover, we have $e_i(p,b,s) = f_i(p,b,s) = s$ for $i=1,2$.  However, although ${e_2(e,b,f) = e_2(p,b,s)}$ and $f_2(e,b,f) = f_2(p,b,s)$, we have
\begin{equation}\label{eq:e1f1}
e_1(e,b,f) = \tmat{\sfrac16 & \sfrac16 & \sfrac13\\\sfrac16 & \sfrac16 & \sfrac13\\\sfrac13 & \sfrac13 & \sfrac23} \not= e_1(p,b,s)
\AND
f_1(e,b,f) = \tmat{\sfrac13&\sfrac13&\sfrac13\\\sfrac13&\sfrac13&\sfrac13\\\sfrac13&\sfrac13&\sfrac13} \not= f_1(p,b,s).
\end{equation}

Not only can axiom \ref{C2} be quantified over a smaller set of projections when $S$ is a regular $*$-semigroup, but the projections~$e_i$ and~$f_i$ in \eqref{eq:etc} take on much simpler forms.  Specifically, when $(e,f)$ is $b$-linked (and $S$ is a regular $*$-semigroup), we have
\begin{equation}\label{eq:easier}
e_1(e,b,f) = e\th_q \COMMA e_2(e,b,f) = f\De_b \COMMA f_1(e,b,f) = e\Th_b \AND f_2(e,b,f) = f\de_r,
\end{equation}
where again $q=\bd(b)$ and $r=\br(b)$.
For example, using $\th_t=\de_t$, \ref{P4'}, Lemma \ref{lem:Thpa}\ref{Thpa2} and \eqref{eq:lp}, we have
\[
e_1(e,b,f) = f\De_b\de_{e\th_q} = f\De_b\th_{e\th_q} = f\De_b\th_q\th_e \th_q = f\De_b\de_{\bd(b)}\th_e \th_q = f\De_b\th_e \th_q = e \th_q.
\]
The simplifications in \eqref{eq:easier} do not hold in general DRC-semigroups.  Indeed, small counterexamples can be readily constructed using Mace4 \cite{Mace4}, although we do not currently know of any \emph{$*$-regular} counterexamples.

As well as the simplification in the projections \eqref{eq:etc} and \eqref{eq:easier}, the `double restrictions' in~\eqref{eq:lamrhoagain} are just single restrictions for regular $*$-semigroups, as equality of the $\leql$ and $\leqr$ orders leads to the identities $({}_p\corest a)\rest_q = a\rest_q$ and ${}_r\corest (a\rest_s) = {}_r\corest a$ (for appropriate $p,q,r,s$).  As usual, these do not hold in general $*$-regular semigroups.

\section{Other special cases}\label{sect:other}

In this final section we consider another two important subclasses of DRC-semigroups, specifically the Ehresmann and DRC-restriction semigroups considered by Lawson \cite{Lawson1991} and Die and Wang~\cite{DW2024}.  Again, although we could in principle deduce the main results of \cite{Lawson1991,DW2024} from ours, we will not give all of the details.  Rather, we just outline the ways in which our constructions simplify in these cases.

\subsection{Ehresmann semigroups}\label{subsect:E}

We noted in Section \ref{subsect:DRC} that the \emph{Ehresmann semigroups} of Lawson \cite{Lawson1991} are precisely the DRC-semigroups $S$ satisfying the additional law
\[
D(a)D(b) = D(b)D(a) \AND R(a)R(b) = R(b)R(a) \qquad\text{for all $a,b\in S$.}
\]
Given \ref{DRC4}, either of the above identities implies the other, and of course they say that projections of $S$ commute.  This is a very strong property, and to unpack it we first prove the following:

\begin{prop}\label{prop:pqqp}
If $p,q\in\bP(S)$ for a DRC-semigroup $S$, then the following are equivalent:
\ben\bmc2
\item \label{pqqp1} $pq = qp$,
\item \label{pqqp2} $p\th_q = q\th_p$ and $p\de_q = q\de_p$,
\item \label{pqqp3} $p\th_q = q\th_p = p\de_q = q\de_p$,
\item \label{pqqp4} $p\th_q = q\th_p = p\de_q = q\de_p = pq = qp$.
\emc\een
\end{prop}

\pf
\firstpfitem{\ref{pqqp1}$\implies$\ref{pqqp4}}  Suppose first that $pq=qp$.  By Lemma \ref{lem:Dpa} we have $D(pq) = D(pq)p$, and since ${D(pq) = D(qp)}$, the same lemma gives $D(pq) = D(pq)q$.  Together with \ref{DRC1} these imply ${D(pq) = D(pq)pq = pq}$, which says $q\de_p = pq$.  The remaining equalities are analogous.

\pfitem{\ref{pqqp4}$\implies$\ref{pqqp3}$\implies$\ref{pqqp2}}  These implications are clear.

\pfitem{\ref{pqqp2}$\implies$\ref{pqqp1}}  Now suppose $p\th_q = q\th_p$ and $p\de_q = q\de_p$, which means that $R(pq)=R(qp)$ and $D(pq)=D(qp)$.  As in \eqref{eq:DpqRpq}, we have $pq = D(pq)R(pq)$, so also
\[
qp = D(qp)R(qp) = D(pq)R(pq) = pq.  \qedhere
\]
\epf

\begin{cor}\label{cor:pqqp}
If $p,q\in\bP(S)$ for a DRC-semigroup $S$, and if $pq=qp$, then 
\ben
\item \label{pqqp5} the meet $p\wedge q$ exists in the poset $(P,\leq)$, and $p\wedge q=pq=qp$,
\item \label{pqqp6} $p\F q \iff p = q$.
\een
\end{cor}

\pf
\firstpfitem{\ref{pqqp5}}  By Proposition \ref{prop:pqqp} we have $pq=p\th_q\in P$, and moreover $pq=p\th_q\leq q$ and similarly $pq\leq p$.  So $pq$ is a lower bound for $p$ and $q$.  If we had another lower bound $r\leq p,q$, then $r = prp = qrq$ (cf.~\eqref{eq:leqPS}); combined with $pq=qp$ this leads to $r = (pq)r(pq)$, i.e.~$r\leq pq$.

\pfitem{\ref{pqqp6}}  If $p\F q$, then $p = q\de_p$ and $q = p\th_q$, and since $q\de_p=p\th_q$ (by Proposition \ref{prop:pqqp}) it follows that $p=q$.
\epf

Corollary \ref{cor:pqqp} has some important simplifying consequences for an Ehresmann semigroup~$S$, since \emph{all} projections commute.  First, $P = \bP(S)$ is a meet semilattice, and hence a subsemigroup of $S$.  Second, ${\F} = \io_P$ is the equality relation.  Thus, the chain category $\CC(P)$ is simply a copy of~$P$, in which the only compositions are $p = p\circ p$ ($p\in P$), and in which ${\leql} = {\leqr} = {\leq}$ is the order from \eqref{eq:leqP}.  It follows that the chained projection category $(P,\C,\ve) = \bC(S)$ of the Ehresmann semigroup $S$ carries no more information than the projection category $(P,\C)$ itself.  
Moreover, since ${\F}=\io_P$ is trivial, axiom~\ref{C2} simplifies as well.  Indeed, we first note that since $e\F e_1,e_2$ and $f_1,f_2\F f$ by Lemma \ref{lem:ef}, it follows that $e=e_1=e_2$ and $f=f_1=f_2$.  Thus,~\ref{C2} asserts the equality of the double restrictions:
\begin{equation}
\label{eq:lamrhoE} ({}_{p\th_q}\corest b)\rest_f = {}_e\corest(b\rest_{s\de_r}) \qquad\text{for all $b\in\C(q,r)$ and $p,s\in P$, \quad where $e = s\De_b\de_p$ and $f = p\Th_b\th_s$.}
\end{equation}
This can be used to prove the following, which is one of the Ehresmann category axioms from~\cite{Lawson1991}.

\begin{prop}
In an Ehresmann semigroup $S$, we have ${\leql}\circ{\leqr} = {\leqr}\circ{\leql}$.
\end{prop}

\pf
By symmetry, it suffices to prove the inclusion ${\leql}\circ{\leqr} \sub {\leqr}\circ{\leql}$.  To do so, suppose $a,b\in \C = \C(S)$ are such that $a\leql u \leqr b$ for some $u\in\C$.  This means that
\[
a = {}_p\corest(b\rest_s) \qquad\text{for some $s\leq\br(b)$ and $p\leq\bd(b\rest_s) = s\vd_b$.}
\]
(So here $u = b\rest_s$.)  We can complete the proof by showing that
\begin{equation}\label{eq:pbt}
a = ({}_p\corest b)\rest_t, \WHERE t=\br(a).
\end{equation}
To do so, we now write $q=\bd(b)$ and $r=\br(b)$, and we work towards interpreting \eqref{eq:lamrhoE} in the current context.  First note that 
\[
s\leq\br(b) = r \AND p \leq \bd(b\rest_s) \leq \bd(b) = q,
\]
which imply that $p\th_q = p$ and $s\de_r = s$.  We also have
\[
e = s\De_b\de_p = s\de_r\vd_b\de_p = s\vd_b\de_p = p,
\]
because $p \leq s\vd_b$; cf.~\eqref{eq:leqP}.  Thus, \eqref{eq:lamrhoE} gives
\[
({}_p\corest b)\rest_f = {}_p\corest(b\rest_s),
\]
where $f = p\Th_b\th_s$.  Since $a = {}_p\corest(b\rest_s) = ({}_p\corest b)\rest_f$, it follows that $f = \br(a)$, completing the proof of~\eqref{eq:pbt}, and hence of the proposition.
\epf

Lawson in \cite{Lawson1991} denotes by $\leqe$ the composition ${\leql}\circ{\leqr} = {\leqr}\circ{\leql}$, and we note that the above proof shows that if $a \leqe b$, then we have
\[
a = ({}_p\corest b)\rest_t = {}_p\corest (b\rest_t) , \WHERE p = \bd(a) \ANd t = \br(a).
\]
This is \cite[Lemma 4.6]{Lawson1991}.

Finally, we note that ${\F} = \io_P$ has the consequence that the $\pr$ product on the category $\C$ from Definition \ref{defn:pr} takes on the simpler form
\[
a\pr b = a\rest_e \circ {}_e\corest b , \WHERE e = \br(a) \wedge \bd(b).
\]

\subsection{DRC-restriction semigroups}\label{subsect:rest}

We conclude the paper by briefly commenting on the \emph{DRC-restriction semigroups} considered by Die and Wang in~\cite{DW2024}.  These are the DRC-semigroups $S\equiv(S,\cdot,D,R)$ satisfying the additional laws
\[
D(ab)a = aD(R(a)b) \AND aR(ba) = R(bD(a))a \qquad\text{for all $a,b\in S$}.
\]
It was observed in \cite[Section 1]{DW2024} that these laws need only be quantified over $a\in S$ and $b\in\bP(S)$, given \ref{DRC2}.  The DRC-restriction laws might seem unwieldy, but it turns out that they are equivalent to the orders $\leql$ and $\leqr$ coinciding, as shown in \cite[Lemma 5.8]{DW2024}.  

The full subcategory $\DRCR$ of $\DRC$ consisting of all DRC-restriction semigroups (and DRC-morphisms) is isomorphic to its image under the isomorphism $\bC:\DRC\to\CPC$ from Theorem~\ref{thm:Sfunctor} (cf.~Theorem \ref{thm:iso}).  One could then use this to deduce the results from \cite{DW2024}.  Again we omit the details, but we do note that the constructions and results of \cite{DW2024} involve some (but not all) of the simplifications that we discussed in Section \ref{sect:*} for regular $*$-semigroups.  As we have already mentioned, the~$\leql$ and~$\leqr$ orders coincide, so one works with ordered rather than biordered categories.  However, the~$\th_p$ and~$\de_p$ operations need not be equal, i.e.~$\bP(S)$ need not be symmetric.  Similarly, the~$\vt/\Th$ and~$\vd/\De$ maps need not be related to each other; for example, $\vt_a$ ($a\in S$) need not be equal to any $\vd_b$ map (for $b\in S$).  On the other hand, $\vt_a$ is always a projection algebra morphism $\bd(a)^\da\to\br(a)^\da$, with a similar statement for $\vd_a$; see \cite[Lemmas 4.7 and 7.7]{DW2024}.  This leads to the various equivalent conditions in the analogue of \ref{C1} in \cite{DW2024}, including~\eqref{eq:G1b}; cf.~Remark \ref{rem:G1}.  Likewise, the simplifications in \ref{C2} discussed at the end of Section \ref{subsect:CPC*} all hold in DRC-restriction semigroups; specifically, the axiom can be quantified over $b$-linked pairs (cf.~\eqref{eq:lp}), and the projections $e_i,f_i$ take on the simpler forms in~\eqref{eq:easier}.

\footnotesize
\def\bibspacing{-1.1pt}
\bibliography{biblio}
\bibliographystyle{abbrv}

\end{document}